\documentclass[11pt]{article}

\usepackage[a4paper,margin=1in]{geometry}
\usepackage{amsmath,amssymb,amsthm,mathtools}
\usepackage{enumitem}

\usepackage{xcolor}
\usepackage{doi}
\usepackage{comment}
\usepackage[numbers]{natbib}
\usepackage{hyperref}
\hypersetup{hypertexnames=false}

\newtheorem{theorem}{Theorem}[section]
\newtheorem{lemma}{Lemma}[section]
\newtheorem{corollary}{Corollary}[section]
\newtheorem{remark}{Remark}[section]
\newtheorem{example}{Example}[section]
\newtheorem{proposition}{Proposition}[section]
\newtheorem{assumption}{Assumption}[section]

\newtheorem*{question}{Question}

\usepackage[normalem]{ulem}
\usepackage{subcaption}

\DeclareMathOperator{\supp}{supp}

\DeclareMathOperator{\tr}{Tr}

\title{Beyond weak convergence of signed spectral measures: Eigenvector overlaps and high-dimensional prediction risk}
\author{Thomas Stark\thanks{Department of Mathematics, Aarhus University, Denmark.
Email: \texttt{thomas.stark@math.au.dk}.},
Georg Köstenberger\thanks{Department of Statistics and Mathematics, Vienna University of Economics and Business, Vienna, Austria. Email: \texttt{georg.koestenberger@wu.ac.at}.}\ and Nina Dörnemann\thanks{Department of Mathematics, Aarhus University, Denmark.
Email: \texttt{ndoernemann@math.au.dk}.}
}
\date{}

\usepackage{titletoc}
\usepackage{graphicx}
\usepackage{float}

\begin{document}

\maketitle
\begin{abstract}
Many results in random matrix theory and high-dimensional statistics depend on weak convergence of spectral measures, which is usually established via pointwise convergence of Stieltjes transforms. 
This approach has recently been applied to signed measures. 
However, for signed measures, pointwise convergence of Stieltjes transforms alone no longer characterizes weak convergence. 
While this is unproblematic in many cases, in the high-dimensional regression setting, recently considered by Green and Romanov, the underlying signed \textit{eigenvector-overlap measures} can fail to converge weakly, despite pointwise convergence of their Stieltjes transforms.
 Using a novel approach to risk analysis, we show that weak convergence can be replaced by a simple convergence result for bilinear forms. This allows us to establish convergence of the prediction risk for estimators based on arbitrary continuous spectral filters, including (regularized) gradient descent and ridge regression, even when the underlying signed measures fail to converge weakly.
Under additional regularity conditions, we construct an explicit finite signed measure yielding an integral representation of the asymptotic prediction risk. 
This representation provides necessary and sufficient conditions
for minimum-norm interpolation and optimally tuned ridge regression
to attain the smallest asymptotic prediction risk among
estimators with continuous spectral filters.
Moreover, we provide characterizations of weak and vague convergence of finite signed measures in terms of pointwise convergence of their Stieltjes transforms.

\end{abstract}



\section{Introduction}

Weak convergence of empirical spectral measures is a central
topic in random matrix theory, going back to the seminal
works of \cite{wigner1955,MarchenkoPastur1967}; see also
\cite{Wachter1980,Silverstein1995empirical,
SilversteinChoi1995,Silverstein1995Strong,
BaiSilverstein1998}.
The resulting spectral theory has played an important role
in the development of high-dimensional statistics; see,
e.g., \cite{Johnstone2001,LedoitWolf2004,BBP2005,
ElKaroui2006,Hastie2022}. 
While the theory has primarily focused on probability measures,
recent contributions \cite{dobriban-sharp-detection,
LouartLiaoCouillet2018,GreenRomanov2025,
LedoitPeche2011,Sur2026} considered
Stieltjes transforms of signed measures. Similar to classical results for probability
measures, a central step in many of these arguments is the implication
\begin{equation}\label{eq:stieltjes-implication}
    \text{pointwise convergence of Stieltjes transforms}
    \implies
    \text{weak/vague convergence of measures.}
\end{equation}
While some popular references claim that this implication is valid for arbitrary sequences of finite signed measures (see, e.g., Theorem B.9 in \cite{BaiSilverstein2010}), it can indeed fail without further assumptions (see Examples~\ref{ex:stieltjes-not-imply-weak},\ref{ex:vague-convergence}, and Theorem \ref{thm:weak-convergence} below).
While this may introduce a formal gap in the corresponding arguments, it does not necessarily invalidate their conclusions. Indeed, in some cases, such as
\cite{LouartLiaoCouillet2018,LedoitPeche2011,Sur2026},
the additional conditions (see Section \ref{sec:weak}) required for Implication \eqref{eq:stieltjes-implication} are satisfied.
However, in other cases, such as \cite{GreenRomanov2025},
we show that the underlying signed measures
need not converge vaguely. Nevertheless, numerical evidence suggests that the
asymptotic prediction-risk formulas derived in
\cite{GreenRomanov2025} are valid.
This raises the question:
\begin{question}\label{question}
    Can the conclusions of these works be recovered when the underlying sequence of measures has neither a weak nor a vague limit?
\end{question}
To answer this question, one has to leave the well-trodden path of weak and vague convergence arguments.
From a mathematical point of view, this is based on the following observation: weak convergence tests against continuous bounded functions on { $\mathbb{R}$}, while the asymptotic analysis of the prediction risk frequently depends only on functions of the form
\begin{align*}
    (s,t)\mapsto f(s)g(t),
    \qquad f,g\in{ C_b(\mathbb{R})},
\end{align*}
cf. Proposition~\ref{prop:risk-decomposition}.
Convergence on this smaller class of test functions is closely related to the behavior of multivariate Stieltjes transforms, which in turn is covered by a surprisingly simple convergence result for bilinear forms.
This approach allows us to establish convergence of the
prediction risk for estimators with arbitrary continuous
spectral filters, including (regularized) gradient descent
and ridge regression, without requiring weak convergence
of the underlying signed measures. 
For this filter class,  our approach provides an alternative to the weak-convergence argument used in \cite{GreenRomanov2025} and yields new insights into the
structure of the asymptotic risk and the optimality of minimum-norm interpolation and ridge regression.

\paragraph*{Contribution}
Our contributions are fourfold:
\begin{enumerate}
    \item
    In Section~\ref{sec:weak}, we characterize weak and vague
    convergence of finite signed measures in terms of their
    Stieltjes transforms and provide counterexamples to both
    implications of \cite[Theorem~B.9]{BaiSilverstein2010}, which covers Implication \eqref{eq:stieltjes-implication}.
    For measures eventually supported on a common compact set,
    uniformly bounded total variation is a necessary and
    sufficient additional condition for
    Implication~\eqref{eq:stieltjes-implication}.

    \item
    In Section~\ref{sec:application}, we construct a Gaussian
    random-design model with two distinct population eigenvalues
    and a deterministic signal satisfying the assumptions of
    \cite{GreenRomanov2025}, whose underlying (signed) eigenvector-overlap measures
    have unbounded total variation and fail to converge weakly
    or vaguely, despite pointwise convergence of their
    Stieltjes transforms. Consequently, convergence of the prediction risk cannot,
in general, be justified by appealing to weak convergence
of the eigenvector-overlap measures.

    \item
  In Section~\ref{Out-of-sample prediction risk}, we establish
a convergence result for uniformly bounded bilinear forms. This yields asymptotic prediction-risk formulas for
    estimators with arbitrary continuous spectral filters,
    including (regularized) gradient descent and ridge regression,
    without requiring weak convergence of the eigenvector-overlap
    measures.

    \item 
    Under additional regularity conditions, we construct an
explicit finite signed Borel measure representing the asymptotic risk, even when the
eigenvector-overlap measures fail to converge weakly.
When the limiting signal-weighted and population spectral
measures coincide, this representation holds without
additional regularity conditions and reduces to the diagonal
pushforward of a positive one-dimensional measure. Based on this representation, we provide necessary and sufficient conditions
for minimum-norm interpolation and optimally tuned ridge
regression to be asymptotically optimal among all continuous
spectral filters, and identify the corresponding optimal
ridge regularization parameter. Among other things, this extends previous ridge-optimality results
of \cite{DobribanWager2018,Hastie2022,Sur2026} to deterministic signals under general population covariance structures.

\end{enumerate}

\paragraph*{Related work} 
We first explain why Implication~\eqref{eq:stieltjes-implication}
may seem plausible.
For tight sequences of probability measures,
pointwise convergence of Stieltjes transforms implies
weak convergence; see, e.g.,
\cite{BaiSilverstein2010,Fleermann2023,
tao-topic-in-random-matrix-theory}. Moreover, Theorem B.9 in \cite{BaiSilverstein2010} claims that pointwise convergence of Stieltjes transforms is equivalent to vague convergence of finite signed measures.
Since vague convergence and tightness together frequently imply weak convergence, this may suggest that Implication \eqref{eq:stieltjes-implication} also holds for tight families of signed measures. 
However, both implications of Theorem B.9 in \cite{BaiSilverstein2010}  can fail, and hence, Implication \eqref{eq:stieltjes-implication} does not hold in full generality for finite signed measures (see Example \ref{ex:stieltjes-not-imply-weak} and Example \ref{ex:vague-convergence} below).
We provide a full characterization of weak convergence of signed measures in terms of Stieltjes transforms in Theorem \ref{thm:weak-convergence} below, and give sufficient conditions for vague convergence in Theorem \ref{thm:vague-convergence}. 
Related convergence questions for signed measures, independent
of Stieltjes transforms, are studied in
\cite{Herdegen2022,Stanek2026}.
A closely related phenomenon for Laplace transforms is considered
in \cite{herdegen2023continuitytheoremgeneralisedsigned},
where additional conditions are likewise required when passing
from positive to signed measures.

Signed measures appear in a number of recent works.
A simple example occurs in the random matrix analysis of neural
networks in \cite{LouartLiaoCouillet2018}, where the relevant signed
measures are differences of two possibly random probability measures.
Their total variation norms are therefore bounded by two, and
Theorem~\ref{thm:vague-convergence} justifies the passage from
pointwise convergence of their Stieltjes transforms to vague
convergence.
Similarly, \cite[Lemma~6]{LedoitPeche2011} considers
eigenvector-weighted spectral measures associated with a bounded
real-valued function $g$. These measures need not be positive,
but their total variation norms are uniformly bounded by
$\sup_{x\in\mathbb R}|g(x)|$.
Thus, in both cases, the total-variation condition required by
our convergence criteria is satisfied.

In \cite{dobriban-sharp-detection}, the derivative of the
Marchenko--Pastur map is used to study the detection of weak
principal components in high-dimensional PCA.
The proof of the existence of this derivative as a weak limit
invokes Theorem~B.9 in \cite{BaiSilverstein2010} to deduce
weak convergence of finite signed measures from pointwise
convergence of their Stieltjes transforms.
Unlike differences of probability measures, the difference
quotients involved here do not automatically have uniformly
bounded total variation norms.

Signed measures also arise in the analysis of spectral shrinkage
and multi-step self-distillation under finite-spike covariance
models in \cite{Sur2026}.
In this case, \cite[Lemma~A.10]{Sur2026} establishes a uniform
total-variation bound for the relevant signed measures, and hence, the
corresponding vague-convergence argument falls within the scope
of our results.
By contrast, the assumptions imposed in
\cite{GreenRomanov2025} do not guarantee uniform
total-variation bounds for their eigenvector-overlap measures;
see Proposition~\ref{prop:GreenRomanov-TV-counterexample}
in Section~\ref{sec:application}.
For the prediction-risk functionals associated with continuous
spectral filters, we replace the affected weak-convergence step
by a bilinear-form argument in
Section~\ref{Out-of-sample prediction risk}.
Based on this argument we are able to establish convergence of the prediction risk, without requiring the underlying eigenvector-overlap measures to converge weakly. 
Subject to additional assumptions, we are able to characterize the asymptotic prediction risk in terms of a finite signed Borel measure. 
This risk representation yields novel insights into the optimality of minimum norm interpolation and ridge regression.                                            

Ridge optimality has been studied under several distinct
statistical formulations.
Under a Gaussian likelihood and a centered isotropic
Gaussian prior on the regression coefficients, suitably tuned
ridge regression is the posterior mean and hence Bayes-optimal
under squared-error loss; see, for example,
\cite[Remarks~3 and~6]{AliKolterTibshirani2019}.
Asymptotic prediction-risk formulas and optimal ridge tuning
under isotropic random effects and general population
covariance structures are derived in \cite{DobribanWager2018}.
In \cite{Hastie2022}, ridge and minimum-norm prediction
risks are analyzed for deterministic signals and general
population covariance structures, and the optimal ridge
parameter is identified in the isotropic case.
From a minimax perspective, asymptotic minimaxity of ridge
regression over Euclidean spheres under isotropic Gaussian
designs and Gaussian errors is established in
\cite{Dicker2016}.
More recently, optimality of tuned ridge regression among
spectral shrinkage estimators for isotropic covariances
is established in \cite[Lemma~3.2]{Sur2026}.
For finite-spike covariance models, the same work
characterizes optimal spectral filters and realizes
them through multi-step self-distillation.

Our results extend ridge-optimality characterizations
beyond isotropic population covariances and isotropic
random-effects assumptions.
Specifically, Theorem~\ref{cor:ridge-optimality} and Corollary~\ref{cor:ridge-optimality2} provide
necessary and sufficient conditions for ridge regression
to be optimal among all continuous spectral filters
for deterministic signals and general limiting population
spectra.
We show that ridge regression can remain optimal even for
anisotropic population covariances, provided the restrictions
of the limiting population and signal-weighted spectral
measures to $(0,\infty)$ are proportional.
The corresponding optimal regularization parameter is the
usual isotropic choice rescaled by the ratio of the masses
of the population and signal-weighted spectral measures
on $(0,\infty)$, respectively, thus accounting for their
potentially different point masses at zero.
Conversely, when these restrictions are not proportional,
no choice of the ridge parameter can attain the optimal
asymptotic risk among continuous spectral filters.

\section{Preliminaries}\label{sec:preliminaries}

Throughout this article, all measures are Borel measures. We write
\[
    \mathbb C^\pm
    \coloneqq
    \{z\in\mathbb C:\pm\operatorname{Im}(z)>0\}
\]
for the upper and lower complex half-plane, respectively. 
For a complex measure $\mu$ on $\mathbb R^d$, we denote its
total variation measure by $|\mu|$ and its total variation norm by $\|\mu\|_{\mathrm{TV}}\coloneqq|\mu|(\mathbb R^d)$. If $\mu$ is a finite signed measure, we denote its Jordan decomposition
by $\mu=\mu^+-\mu^-$, in which case $|\mu|=\mu^++\mu^-$.
For a finite complex measure $\mu$, we define its support by
\[
    \supp(\mu)\coloneqq\supp(|\mu|).
\]
In particular, for a finite signed measure,
\[
    \supp(\mu)
    =
    \supp(\mu^+)\cup\supp(\mu^-).
\]
For a Borel measure $\nu$ on $\mathbb R^d$ and Borel sets
$A,B\subseteq\mathbb R^d$, we denote by $\left.\nu\right|_A$
the restriction of $\nu$ to $A$, defined by
\[
    \left.\nu\right|_A(B)
    \coloneqq \nu(A\cap B).
\]
A family $\mathcal M$ of finite positive measures on $\mathbb R^d$ is
called tight if, for every $\varepsilon>0$, there exists a compact set
$K\subseteq\mathbb R^d$ such that
\[
    \sup_{\nu\in\mathcal M}\nu(K^c)<\varepsilon.
\]
For $\mathbb F\in\{\mathbb R,\mathbb C\}$ and
$A\subseteq\mathbb R^d$, we denote by $C(A;\mathbb F)$
the space of continuous functions $f:A\to\mathbb F$.
If $A$ is compact, this is a Banach space equipped with
the supremum norm $\|f\|_{A,\infty}\coloneqq \sup_{x\in A}|f(x)|$.
For $A\subseteq\mathbb R^d$, we write
$C(A)\coloneqq C(A;\mathbb R)$. Similarly, we denote by $C_b(\mathbb R^d;\mathbb F)$, $C_c(\mathbb R^d;\mathbb F)$, $C_0(\mathbb R^d;\mathbb F)$ the spaces of bounded continuous, compactly supported continuous, and
continuous functions vanishing at infinity, respectively. More
precisely, $f\in C_0(\mathbb R^d;\mathbb F)$ if, $f(x)\longrightarrow0$ as $\|x\|_2\longrightarrow\infty$. The spaces $C_b(\mathbb R^d;\mathbb F)$, $C_c(\mathbb R^d;\mathbb F)$, and
$C_0(\mathbb R^d;\mathbb F)$  are equipped with the supremum norm $\|f\|_\infty
    \coloneqq
    \sup_{x\in\mathbb R^d}|f(x)|$.
When $\mathbb F=\mathbb R$, we omit the range and simply write
$C_b(\mathbb R^d)$, $C_c(\mathbb R^d)$, and $C_0(\mathbb R^d)$.

We say that a sequence of finite signed measures $\mu_n$ on
$\mathbb R^d$ converges weakly to a finite signed measure $\mu$, and
write
\[
    \mu_n\Rightarrow\mu,
\]
if
\[
    \int_{\mathbb R^d}f\,d\mu_n
    \longrightarrow
    \int_{\mathbb R^d}f\,d\mu,\qquad \text{for every } f\in C_b(\mathbb R^d).
\]
If $\mu_n$ and $\mu$ are supported
on a compact set $K\subseteq\mathbb R^d$, then
$\mu_n\Rightarrow\mu$ if and only if
\[
    \int_K f\,d\mu_n
    \longrightarrow
    \int_K f\,d\mu,\qquad\text{for every }f\in C(K).
\]
We say that $\mu_n$ converges vaguely to $\mu$ if
\[
    \int_{\mathbb R^d}f\,d\mu_n
    \longrightarrow
    \int_{\mathbb R^d}f\,d\mu,\qquad\text{for every }f\in C_c(\mathbb R^d).
\]
For a finite complex measure $\mu$ on $\mathbb R^d$, we define its
multivariate Stieltjes transform by
\[
    m_\mu(z_1,\ldots,z_d)
    \coloneqq
    \int_{\mathbb R^d}
    \prod_{j=1}^d\frac{1}{x_j-z_j}
    \,d\mu(x),
    \qquad
    z_1,\ldots,z_d\in\mathbb C\setminus\mathbb R.
\]
If $\mu$ is compactly supported, then there exists $M>0$ such that $\supp(\mu)\subseteq[-M,M]^d$,
and, for such $M$, the same formula defines $m_\mu$ for $z_1,\ldots,z_d\in\mathbb C\setminus[-M,M]$. 
For vectors, $\|\cdot\|_2$ denotes the Euclidean norm, whereas for
matrices, $\|\cdot\|_2$ denotes the corresponding operator norm and
$\|\cdot\|_F$ the Frobenius norm. We write $I_p$ for the $p\times p$ identity matrix
and $\operatorname{spec}(A)$ for the spectrum of a square matrix $A$. For any normed space $X$, we denote its dual by $X^*$ and
write $\delta_x$ for the Dirac measure at $x$. Unless explicitly stated otherwise, all limits are taken as $n\to\infty$. Finally, we use $c,C>0$ to denote generic constants,
independent of $n$, whose values may change from one
occurrence to another.

\section{Stieltjes-transform criteria for weak and vague convergence}\label{sec:weak}

Pointwise convergence of Stieltjes transforms is a commonly used
criterion for weak or vague convergence; see, e.g.,
\cite{BaiSilverstein2010,Fleermann2023,tao-topic-in-random-matrix-theory}.
The following example shows that, for finite signed measures,
pointwise convergence of the Stieltjes transforms to zero need not
imply weak convergence to the zero measure, even when all measures
are supported on a common compact set.

\begin{example}\label{ex:stieltjes-not-imply-weak}
 Let \(K=[0,1]\) and, for each \(n\in\mathbb N\), define the finite signed
Borel measure
\[
        \mu_n
        \coloneqq
        \sqrt n\left(\delta_{1/n}-\delta_0\right)
\]
on \(K\). Since the positive and negative parts of \(\mu_n\) are given by
\[
        \mu_n^+=\sqrt n\,\delta_{1/n},
        \qquad
        \mu_n^-=\sqrt n\,\delta_0,
\]
the total variation norm of $\mu_{n}$ is given by
  \begin{equation*}
    \|\mu_{n}\|_{TV} = \mu_n^+(K)+\mu_n^-(K)= \sqrt{n}(\delta_{\frac{1}{n}}(K) + \delta_{0}(K)) = 2\sqrt{n},
  \end{equation*}
  and hence the total variation norms of the $\mu_{n}$'s are not uniformly bounded. 
  For any $z\in \mathbb{C}\setminus K$, we have
  \begin{equation*}
    m_{\mu_{n}}(z) = \int_{K} \frac{1}{x-z}d \mu_{n}(x) = \sqrt{n} \biggl( \frac{1}{\frac{1}{n} - z} + \frac{1}{z} \biggr)= \frac{1}{\sqrt{n}} \frac{1}{z(\frac{1}{n} - z)}\to 0. 
  \end{equation*}
  In words, the Stieltjes transforms of $\mu_{n}$ converge to the Stieltjes transform of the zero measure $\mu = 0$. However, $\mu_{n}$ does not converge to $\mu$ weakly. 
  To see this, let $f(x) = \sqrt{x}$, and note that 
  \begin{equation*}
    \int_{K} f(x) d \mu_{n}(x) = \sqrt{n} (\sqrt{1/n} - \sqrt{0}) = 1 \neq 0 = \int_K f(x) d \mu(x). 
  \end{equation*}
\end{example}

The obstruction to Implication \eqref{eq:stieltjes-implication} in the previous example is the diverging total variation:
the positive and negative contributions asymptotically cancel in
the Stieltjes transform, although their total masses diverge.
For measures supported on a common compact set, a uniform bound
on the total variation is a necessary and sufficient additional condition 
to deduce weak convergence from pointwise convergence of the
Stieltjes transforms.
The following proposition states this criterion for measures
eventually supported on a common compact set, a setting frequently
encountered in random matrix theory and high-dimensional statistics.

\begin{proposition}
\label{cor:weak-convergence-compact-support}
Let $\mu_n$, $n\geq1$, be finite signed Borel measures on
$\mathbb R^d$. Suppose that there exist a compact set
$K\subseteq\mathbb R^d$ and $n_0\in\mathbb N$ such that
\[
    \operatorname{supp}(\mu_n)\subseteq K,
    \qquad \text{for every } n\geq n_0.
\]
Then the following are equivalent:
\begin{enumerate}
    \item $\mu_n$ converges weakly to some finite signed Borel measure
    $\mu$.
    \item The following conditions hold:
    \begin{enumerate}
        \item $m_{\mu_n}(z_1,\ldots,z_d)$ converges for every
        $z_1,\ldots,z_d\in\mathbb C^+$;
        \item $\sup_{n\geq1}\|\mu_n\|_{\mathrm{TV}}<\infty.$
    \end{enumerate}
\end{enumerate}
Whenever these equivalent conditions hold, the weak limit $\mu$ is
supported on $K$ and is uniquely determined by
\[
    m_\mu(z_1,\ldots,z_d)
    =
    \lim_{n\to\infty}m_{\mu_n}(z_1,\ldots,z_d),
    \qquad
    z_1,\ldots,z_d\in\mathbb C^+.
\]
\end{proposition}
Without a common compact support, mass may escape to infinity,
so uniformly bounded total variation alone is no longer sufficient.
The following theorem gives a complete characterization
of weak convergence in terms of Stieltjes transforms.
Proposition~\ref{cor:weak-convergence-compact-support} follows
as a special case.
\begin{theorem}\label{thm:weak-convergence}
Let $\mu_n$, $n\geq 1$, be finite signed Borel measures on
$\mathbb R^d$. The following are equivalent:
\begin{enumerate}
    \item $\mu_n$ converges weakly to some finite signed Borel measure
    $\mu$.
    \item The following conditions hold:
    \begin{enumerate}
        \item $m_{\mu_n}(z_1,\ldots,z_d)$ converges for every
        $z_1,\ldots,z_d\in\mathbb C\setminus\mathbb R$;
        \item the family $\{|\mu_n|:n\geq1\}$ is tight;
        \item $\sup_{n\geq1}\|\mu_n\|_{\mathrm{TV}}<\infty.$
    \end{enumerate}
\end{enumerate}
Whenever these equivalent conditions hold, the weak limit $\mu$ is
uniquely determined by
\[
    m_\mu(z_1,\ldots,z_d)
    =
    \lim_{n\to\infty}m_{\mu_n}(z_1,\ldots,z_d),
    \qquad
    z_1,\ldots,z_d\in\mathbb C\setminus\mathbb R.
\]
\end{theorem}
We next turn to vague convergence. The role of tightness is fundamentally
different in this setting: mass may escape to infinity without affecting
integrals against compactly supported test functions. Consequently,
tightness is neither necessary nor expected. The role of total variation is more
subtle. Unlike weak convergence, vague convergence does not itself imply
a uniform bound on the total variation, as the following example shows.

\begin{example}\label{ex:vague-convergence}
Let $\mu_n=\delta_0+n\delta_n$, $\mu=\delta_0$. Then $\mu_n\to\mu$ vaguely. Indeed, if
$f\in C_c(\mathbb R)$, then $f(n)=0$ for all sufficiently large $n$,
and hence
\[
    \int_{\mathbb R} f\,d\mu_n
    =
    f(0)+nf(n)
    =
    f(0)
    =
    \int_{\mathbb R} f\,d\mu
\]
for all sufficiently large $n$. 
On the other hand, $\|\mu_n\|_{\mathrm{TV}}=n+1$,
and, for every $z\in\mathbb C\setminus\mathbb R$,
\[
    m_{\mu_n}(z)
    =
    -\frac{1}{z}
    +
    \frac{n}{n-z}
    \longrightarrow
    -\frac{1}{z}+1
    \neq
    -\frac{1}{z}
    =
    m_\mu(z).
\]

\end{example}

Thus, vague convergence need not imply either uniformly bounded total
variation or convergence of the corresponding Stieltjes transforms to the Stieltjes transform of the vague limit.
Example \ref{ex:stieltjes-not-imply-weak} and Example \ref{ex:vague-convergence} imply that both directions of Theorem B.9 in \cite{BaiSilverstein2010} can fail without a uniform bound on the total variation.

We also note that a uniform total-variation bound does not imply tightness, as
the sequence $\mu_n=\delta_n$ shows; see also
\cite[Proposition~2.4]{Herdegen2022}. This is consistent with the nature
of vague convergence: indeed, $\delta_n\to0$ vaguely even though the
family $\{\delta_n:n\geq1\}$ is not tight.
Nevertheless, once uniform total-variation control is imposed, pointwise
convergence of the multivariate Stieltjes transforms  characterizes
vague convergence. 

\begin{theorem}\label{thm:vague-convergence}
Let $\mu_n$, $n\geq1$, be finite signed Borel measures on
$\mathbb R^d$ satisfying
\[
    \sup_{n\geq1}\|\mu_n\|_{\mathrm{TV}}<\infty.
\]
Then the following are equivalent:
\begin{enumerate}
    \item $\mu_n$ converges vaguely to some finite signed Borel measure
    $\mu$.
    \item $m_{\mu_n}(z_1,\ldots,z_d)$ converges for every
    $z_1,\ldots,z_d\in\mathbb C\setminus\mathbb R$.
\end{enumerate}
Whenever these equivalent conditions hold,
\[
    \int_{\mathbb R^d}f\,d\mu_n
    \longrightarrow
    \int_{\mathbb R^d}f\,d\mu
\]
for every $f\in C_0(\mathbb R^d;\mathbb C)$, and
\[
    m_\mu(z_1,\ldots,z_d)
    =
    \lim_{n\to\infty}
    m_{\mu_n}(z_1,\ldots,z_d),
    \qquad
    z_1,\ldots,z_d\in\mathbb C\setminus\mathbb R.
\]
In particular, the vague limit $\mu$ is uniquely determined by the
pointwise limit of the Stieltjes transforms.
\end{theorem}

The proofs of Theorems~\ref{thm:weak-convergence}
and~\ref{thm:vague-convergence} follow the same subsequence argument.
For weak convergence, compactness follows from Prokhorov's theorem;
for vague convergence, it follows from the Banach--Alaoglu theorem.
The proof of Theorem~\ref{thm:weak-convergence} and Theorem~\ref{thm:vague-convergence} are deferred to Section~\ref{appendix_A} of the Supplementary Material. We would also like to emphasize that the Stieltjes transforms in Theorem \ref{thm:weak-convergence} and Theorem \ref{thm:vague-convergence} need to converge pointwise on $(\mathbb{C}\setminus \mathbb{R})^d$, since convergence on $(\mathbb{C}^+)^d$ is not enough to uniquely determine the limit in the multivariate case (for details, we refer to Proposition~\ref{prop:multivariate-stieltjes-inversion} in Section~\ref{appendix_A} of the Supplementary Material).

Next, we examine how these conditions arise for signed measures in
statistical random matrix theory.

\section{Signed measures in high-dimensional statistics:
Applications and obstructions}\label{sec:application}

Finite signed measures arise naturally in random matrix theory and its
applications to high-dimensional statistics. A simple example occurs in
the analysis of neural networks in \cite{LouartLiaoCouillet2018}, where
the relevant signed measures are differences of two possibly random
probability measures. Their total variation norms are therefore bounded
by two, and the convergence of their Stieltjes transforms established in
\cite[Proof of Theorem~2]{LouartLiaoCouillet2018} implies vague
convergence by Theorem~\ref{thm:vague-convergence}. A related
construction appears in \cite{dobriban-sharp-detection}, where scaled
differences of probability measures arise in the analysis of derivatives
of the Marchenko--Pastur map.
Another, structurally different, form of signed measures is provided
by bilinear spectral statistics arising from multi-resolvent
expressions, for example in the analysis of spectral shrinkage
estimators under spiked covariance models. In this setting, boundedness
of the total variation norms is no longer automatic. Nevertheless, a
uniform bound is established in \cite[Lemma~A.10]{Sur2026}.
 
The remainder of this section is concerned with the high-dimensional
principal component regression setting of \cite{GreenRomanov2025}.
We show that the assumptions imposed in \cite{GreenRomanov2025} do
not guarantee uniform control of the total variation norms of the
signed eigenvector overlap measures. Consequently, pointwise
convergence of their bivariate Stieltjes transforms is not sufficient
to conclude weak convergence. Proposition
\ref{prop:GreenRomanov-TV-counterexample} below gives a simple Gaussian
random-design example in which the total variation norms diverge.
Despite the simple form of the example, establishing this divergence
is technically involved, and hence, we defer the proof to
Section~\ref{sec_proofs}.

We begin by introducing the high-dimensional regression
framework and the signed eigenvector overlap measures that arise in \cite{GreenRomanov2025}.
Let $\{(x_i,y_i)\}_{i=1}^n$ be i.i.d.\ training observations taking values in $\mathbb R^p\times\mathbb R$, and let $(x_0,y_0)$ be an independent test observation with the same distribution as the training observations. We assume that 
\[ 
y_i=x_i^\top\beta_n^\star+\varepsilon_i, \qquad i=0,\ldots,n, 
\] 
where $\beta_n^\star\in\mathbb R^p$, $\beta_n^\star \neq 0$, $\mathbb E[x_i]=0$, $\Sigma_n = \mathbb E[x_ix_i^\top]$ is positive semidefinite, $\varepsilon_0,\ldots,\varepsilon_n$ are i.i.d., independent of the covariates, and satisfy $\mathbb E[\varepsilon_i]=0$, $\mathbb E[\varepsilon_i^2]=\sigma^2$ for some $\sigma\geq0$. Write $X\coloneqq(x_1,\ldots,x_n)^\top$, $\varepsilon \coloneqq(\varepsilon_1,\ldots,\varepsilon_n)^\top$ and $y\coloneqq(y_1,\ldots,y_n)^\top$. 
The sample covariance matrix $S_n= \frac{1}{n} X^\top X$ admits the eigendecomposition
\begin{equation}\label{eq:sam-cov-decomposition}
    S_n = \sum_{i=1}^{p}\lambda_i(S_n)w_iw_i^\top.
\end{equation}
For $f\in C([0,\infty))$, we consider the spectral estimator 
\begin{align} \label{eq:spectral-estimator}
\hat\beta(f) \coloneqq f(S_n)\frac{X^\top y}{n}, 
\end{align} where $f(S_n)$ is defined by functional calculus, i.e.,
\begin{equation*}
    f(S_n) = \sum_{i=1}^p f(\lambda_i(S_n)) w_i w_i^\top.
\end{equation*}
This class includes, among others, ridge regression and finite-step (regularized) gradient descent. In particular, ridge regression with regularization parameter $\lambda>0$ corresponds to the filter 
\[ f_\lambda(s)=\frac{1}{s+\lambda}, \] whereas $t$ steps of (regularized) gradient descent with
regularization parameter $\lambda\geq0$, step size $\eta>0$, and
initialization $\hat\beta_0=0$ correspond to the filter
\[
    f_{t,\lambda,\eta}(s)
    =
    \eta
    \sum_{j=0}^{t-1}
    \{1-\eta(s+\lambda)\}^j.
\]
see, e.g., \cite[Section~2]{StarkSteinberger2025}.
Its conditional (out-of-sample) prediction error is given by 
\[ 
\mathbb E\left[ \bigl(x_0^\top\hat\beta(f)-y_0\bigr)^2 \,\middle|\,X \right] = \mathbb E\left[ \bigl(\hat\beta(f)-\beta_n^\star\bigr)^\top \Sigma_n \bigl(\hat\beta(f)-\beta_n^\star\bigr) \,\middle|\,X \right] +\sigma^2.
\] 
Since the irreducible noise term $\sigma^2$ does not depend on the estimator, we consider the excess (out-of-sample) prediction risk 
\begin{align} \label{eq:def_R_n}
\mathrm R_{\Sigma_n}\bigl(\hat\beta(f)\bigr) \coloneqq \mathbb E\left[ \bigl(\hat\beta(f)-\beta_n^\star\bigr)^\top \Sigma_n \bigl(\hat\beta(f)-\beta_n^\star\bigr) \,\middle|\,X \right]. 
\end{align}

The excess prediction risk depends on two measures; the eigenvector overlap measure $\widehat{\mathsf{K}}_n^{\beta^\star_n}$, and a measure corresponding to the variance term, $\mathsf{Q}_n$, which are given by 
\begin{align} \label{eq:def_K}
\widehat{\mathsf K}_n^{\beta_n^\star}
=
\frac{1}{\|\beta_n^\star\|_2^2}
\sum_{i,j=1}^{p}
\langle\beta_n^\star,w_i\rangle
\langle w_i,\Sigma_n w_j\rangle
\langle w_j,\beta_n^\star\rangle
\,
\delta_{(\lambda_i(S_n),\lambda_j(S_n))}
\end{align}
and 
\begin{align}\label{eq:definition_Q_n}
    \mathsf Q_n = \frac1p
    \sum_{i=1}^p
    w_i^\top\Sigma_n w_i\,\delta_{\lambda_i(S_n)}
\end{align}
respectively. 
Note that $\widehat{\mathsf{K}}_n^{\beta_n^\star}$ is a signed measure, while $\mathsf{Q}_n$ is nonnegative.  
The next proposition expresses the risk of the estimator $\hat{\beta}(f)$  in terms of the measures $\widehat{\mathsf{K}}_n^{\beta^\star_n}$ and $\mathsf{Q}_n$.
\begin{proposition} \label{prop:risk-decomposition} Consider the preceding data model. Then, for every $f\in C([0,\infty))$, 
\begin{align*} 
\mathrm R_{\Sigma_n}\bigl(\hat\beta(f)\bigr) = \|\beta_n^\star\|_2^2 \int_{[0,\infty)^2} \bigl(1-sf(s)\bigr) \bigl(1-tf(t)\bigr) \,d\widehat{\mathsf K}^{\beta_n^\star}_n(s,t)+\sigma^2\gamma_n \int_{[0,\infty)} s f(s)^2\,d\mathsf Q_n(s),
\end{align*}
where $\gamma_n\coloneqq p/n$.
\end{proposition}
The proof is deferred to Section~\ref{appendix_B} of the Supplementary Material.
The most important feature of this representation is the fact that the risk depends on $\widehat{\mathsf{K}}_n^{\beta_n^\star}$ only through a function of the form
\begin{equation*}
    (s,t)\mapsto f(s)g(t),\qquad f,g\in C([0,\infty)).
\end{equation*}
This property will allow us to replace weak convergence by convergence of certain bilinear forms in Section \ref{Out-of-sample prediction risk}.
This is necessary, since the measures $\widehat{\mathsf{K}}_n^{\beta_n^\star}$ may not have uniformly bounded total variation given the assumptions imposed in \cite{GreenRomanov2025}, and hence, may not converge weakly by Theorem \ref{thm:weak-convergence}; see Proposition \ref{prop:GreenRomanov-TV-counterexample} below. 
Indeed, a crude upper bound is given by 
\begin{align*}\label{eq:sig_K_bound}
\begin{split}
\|\widehat{\mathsf K}^{\beta_n^\star}_n\|_{\mathrm{TV}}
&\leq
\frac{1}{\|\beta_n^\star\|_2^2}
\sum_{i,j=1}^{p}
\left|
\langle\beta_n^\star,w_i\rangle
\right|
\left|
\langle w_i,\Sigma_n w_j\rangle
\right|
\left|
\langle w_j,\beta_n^\star\rangle
\right|
\\
&\leq
\frac{1}{\|\beta_n^\star\|_2^2}
\left(
\sum_{i,j=1}^{p}
\langle\beta_n^\star,w_i\rangle^2
\langle w_j,\beta_n^\star\rangle^2
\right)^{1/2}
\left(
\sum_{i,j=1}^{p}
\langle w_i,\Sigma_n w_j\rangle^2
\right)^{1/2}
=
\|\Sigma_n\|_F.
\end{split}
\end{align*}
However, the assumptions of \cite{GreenRomanov2025} do not imply that the Forbenius norm of $\Sigma_n$ is uniformly bounded, they only imply a uniform bound on the spectral norm of $\Sigma_n$. 
Thus, under a uniform bound on \(\|\Sigma_n\|_2\), we have $\|\Sigma_n\|_F\leq
\sqrt{p}\,\|\Sigma_n\|_2 =O(\sqrt{p})$, which yields, almost surely,
\[
\|\widehat{\mathsf K}^{\beta_n^\star}_n\|_{\mathrm{TV}}
=
O(\sqrt p). 
\]
In what follows, we will see that this bound is essentially sharp in its dependence on $p$. In particular, we construct a simple Gaussian random-design model satisfying the assumptions of \cite{GreenRomanov2025} for which the sequence 
$\{\|\widehat{\mathsf K}^{\beta_n^\star}_n\|_{\mathrm{TV}}\}_{n\geq 1}$ is not uniformly bounded.

\begin{proposition}[Failure of uniform total-variation control in \cite{GreenRomanov2025}]
\label{prop:GreenRomanov-TV-counterexample}
Consider the preceding data model and let $n\to\infty$ through multiples of $4$, put
$p=\frac n2$,
and let $Z_n$ be the $n\times p$ upper left-block of a doubly infinite array of independent standard
Gaussian entries. Then, there exists a sufficiently small constant
$\kappa>0$, independent of $n$, such that, for
\[
P_n
=
\begin{pmatrix}
0&0\\
0&I_{p/2}
\end{pmatrix},
\qquad
\Sigma_n
=
I_p+\kappa P_n
=
\begin{pmatrix}
I_{p/2}&0\\
0&(1+\kappa)I_{p/2}
\end{pmatrix},\qquad X_n\coloneqq Z_n\Sigma_n^{1/2},
\]
and any sequence of deterministic vectors
$\beta_n^\star\in\mathbb R^p$ satisfying $\|\beta_n^\star\|_2=1, \|P_n\beta_n^\star\|_2^2=\frac12$,
the following holds: for every fixed
$\varepsilon\in(0,1/2)$, there exists a constant $c>0$ such that
\begin{align}\label{eq:statement_failure}
\mathbb P\left(
\|\widehat{\mathsf K}^{\beta_n^\star}_n\|_{\mathrm{TV}}
\geq
c p^{1/2-\varepsilon}
\right)
\longrightarrow1.
\end{align}
Furthermore, this construction satisfies 
Assumptions~1--6 of \cite{GreenRomanov2025}.
\end{proposition}

\begin{remark} 
The heuristic discussion preceding
\cite[Equation~(99)]{GreenRomanov2025} suggests that
the off-diagonal overlaps
$\langle w_i,\Sigma_n w_j\rangle$, $i\neq j$,
are typically of order $p^{-1/2}$, while there are
of order $p^2$ such terms. In the proof of Proposition~\ref{prop:GreenRomanov-TV-counterexample}
we show that their cumulative absolute contribution can lead to diverging total variation norms.
\end{remark}

Since \eqref{eq:statement_failure} implies that
$\|\widehat{\mathsf K}^{\beta_n^\star}_n\|_{\mathrm{TV}}$
diverges in probability, Theorem~\ref{thm:weak-convergence} rules out
almost-sure weak convergence of
$\widehat{\mathsf K}^{\beta_n^\star}_n$. Hence, the passage from
pointwise convergence of the Stieltjes transforms to weak convergence
in \cite{GreenRomanov2025} requires additional justification.
In the next section, we discuss one remedy based on restricting the
class of test functions. The resulting convergence is sufficient for
the prediction-risk functionals considered below. Moreover, we derive
additional sufficient conditions under which the limiting bilinear
form admits a representation by a finite signed Borel measure.

\section{Prediction risk for continuous spectral filters}
\label{Out-of-sample prediction risk}

The preceding section shows that the signed measures $\widehat{\mathsf K}^{\beta_n^\star}_n$ need not be uniformly bounded in total variation, even under the assumptions imposed in \cite{GreenRomanov2025}. Consequently, pointwise convergence of their bivariate Stieltjes transforms does not, in general, imply  weak convergence by the results of Section~\ref{sec:weak}. 

Nevertheless, the relevant prediction-risk functionals appearing in the weak-convergence step of \cite{GreenRomanov2025} do not require convergence against arbitrary functions in $C([0,M]^2)$.
For the continuous spectral filters considered here, the relevant
test functions have the form
\[ 
(s,u)\mapsto g(s)h(u), \qquad g,h\in C([0,M]; \mathbb{C}), \]
see Proposition \ref{prop:risk-decomposition} above. 
Integrals against this smaller set of test functions are bilinear forms that can converge, even if the total variation of the signed measures in question is unbounded.
To illustrate this distinction, suppose that the eigenvalues of $S_n$ are eventually contained in $[0,M]$, and associate with $\widehat{\mathsf K}^{\beta_n^\star}_n$ the bilinear form \begin{align*} 
B_n(g,h) &\coloneqq \int_{[0,M]}\int_{[0,M]} g(s)h(u)\, d\widehat{\mathsf K}^{\beta_n^\star}_n(s,u) 
= \frac{1}{\|\beta_n^\star\|_2^2} (\beta_n^\star)^\top g(S_n)\Sigma_n h(S_n) \beta_n^\star, 
\end{align*}
for $g,h\in C([0,M];\mathbb C)$. Using the elementary bounds, $\|g(S_n)\|_2\leq \|g\|_{\infty},  \|h(S_n)\|_2\leq \|h\|_{\infty}, $
we have  \begin{align} \label{eq:bound_B_n}
|B_n(g,h)| \leq \|\Sigma_n\|_2 \|g\|_{\infty} \|h\|_{\infty}. \end{align}
Thus, a uniform spectral-norm bound on $\Sigma_n$ yields a uniform
bound on the operator norms of the bilinear forms $B_n$, although it
does not yield a uniform bound on
$\|\widehat{\mathsf K}^{\beta_n^\star}_n\|_{\mathrm{TV}}$.
The bivariate Stieltjes transform of
$\widehat{\mathsf K}^{\beta_n^\star}_n$ is precisely the evaluation
$B_n(r_\zeta,r_\xi)$, where
\[
    r_z(s)\coloneqq\frac1{s-z},
    \qquad s\in[0,M],\quad z\in\mathbb C^+.
\]
Convergence on these kernels, together with uniform boundedness of
the bilinear forms, determines the limit on all pairs of continuous
functions. The corresponding result is stated in
Lemma~\ref{lem:abstract-bilinear-extension} in Section~\ref{sec_proofs} and is
used in the proof of Theorem~\ref{lem:GR-bilinear-convergence} below.

We now combine this observation with the single- and two-resolvent
deterministic equivalents of \cite{DoernemannStark} to obtain
asymptotic risk formulas for continuous spectral filters.
We work under the following assumptions.

\begin{assumption}
\label{ass:GR-basic} In the following, all random variables are defined on a common probability space $(\Omega, \mathcal A, \mathbb P)$. We assume that:
\begin{enumerate}[label=(A\arabic*)]
    \item\label{ass:GR-basic-1}
    The dimension $p=p_n$ satisfies
    \[
        \gamma_n\coloneqq\frac{p}{n}
        \longrightarrow\gamma\in(0,\infty).
    \]

    \item\label{ass:GR-basic-2}
    The data vectors have the form
    $x_i=\Sigma_n^{1/2}z_i$, $1\leq i\leq n$.
    For each $n\in\mathbb N$, the entries of
    $Z_n=(z_1,\ldots,z_n)^\top
    =(z_{ij})_{1\leq i\leq n,\,1\leq j\leq p}$
    are real-valued independent random variables satisfying $\mathbb E(z_{ij})=0$, $\mathbb E(z_{ij}^2)=1$,
    and,
    \[
        \sup_{n\geq1}
        \max_{1\leq i\leq n,\,1\leq j\leq p_n}
        \mathbb E\bigl(|z_{ij}|^{8}\bigr)
        <\infty.
    \]

    \item\label{ass:GR-basic-3}
    The matrices $\Sigma_n$ are deterministic, symmetric,
    nonnegative definite, and satisfy
    \[
        \sup_{n\geq1}\|\Sigma_n\|_2<\infty.
    \]

    \item\label{ass:GR-basic-4}
    The empirical spectral measures
    \[
        H_n\coloneqq
        \frac1{p}\sum_{j=1}^{p}\delta_{\lambda_j(\Sigma_n)},
    \]
    where $\lambda_1(\Sigma_n),\ldots,\lambda_p(\Sigma_n)$ are the eigenvalues
    of $\Sigma_n$, converge weakly to a Borel probability
    measure $H$ on $[0,\infty)$:
    \[
        H_n\Rightarrow H, \qquad \text{with} \qquad H((0,\infty))>0.
    \]

    \item\label{ass:GR-basic-5}
    The signal vectors
    $\beta_n^\star\in\mathbb R^{p}\setminus\{0\}$
    are deterministic. Let $\Sigma_n
        =\sum_{j=1}^{p} \lambda_j(\Sigma_n)q_{j}q_{j}^\top$
    be an orthonormal eigendecomposition. The signal-weighted
    spectral measures
    \[
        G_n\coloneqq
\frac1{\|\beta_n^\star\|_2^2}
        \sum_{j=1}^{p}
        |\langle\beta_n^\star,q_{j}\rangle|^2
        \delta_{\lambda_j(\Sigma_n)}
    \]
    converge weakly to a Borel probability measure $G$ on
    $[0,\infty)$:
    \[
        G_n\Rightarrow G.
    \]
\end{enumerate}
\end{assumption}

We write,
\begin{align*}
    F_{S_n}
    &\coloneqq
    \frac1p
    \sum_{j=1}^p
    \delta_{\lambda_j(n^{-1}X_n^\top X_n)},
    &
    \overline F_{S_n}
    &\coloneqq
    \frac1n
    \sum_{j=1}^n
    \delta_{\lambda_j(n^{-1}X_nX_n^\top)}.
\end{align*}
Under Assumptions~\ref{ass:GR-basic-1}--\ref{ass:GR-basic-4},
\cite{Pan2010} gives, almost surely,
\begin{align}\label{eq:GR-MP-weak-convergence}
\begin{split}
    \overline F_{S_n}
    &\Rightarrow
    \overline F_{\gamma,H},
    \quad\text{and consequently}\\
    F_{S_n}
    &=
    \left(1-\frac1{\gamma_n}\right)\delta_0
    +
    \frac1{\gamma_n}\overline F_{S_n}
    \Rightarrow
    \left(1-\frac1{\gamma}\right)\delta_0
    +
    \frac1{\gamma}\overline F_{\gamma,H}
    \eqqcolon F_{\gamma,H}.
\end{split}
\end{align}
We denote the Stieltjes transform associated with the
Marchenko--Pastur measure $\overline F_{\gamma,H}$ by
\begin{align}\label{eq:companion_stielt}
    v(z)
    \coloneqq
    \int_{[0,\infty)}\frac1{s-z}\,d\overline F_{\gamma,H}(s),
    \qquad z\in\mathbb C^+,
\end{align}
and set
\[
    \Delta(\zeta,\xi)
    \coloneqq
    \begin{cases}
        \displaystyle
        \frac{v(\zeta)-v(\xi)}{\zeta-\xi},
        &\zeta\neq\xi,\\[3mm]
        \displaystyle
        \frac{d}{d\zeta}v(\zeta),
        &\zeta=\xi.
    \end{cases}
\]
The following lemma records the limiting Stieltjes transforms needed
for the prediction-risk formula.
\begin{lemma}
\label{lem:GR-Stieltjes-convergence}
Suppose that Assumption~\ref{ass:GR-basic} holds.
Then almost surely, for every $\zeta,\xi\in\mathbb C^+$,
\begin{align}
\label{eq:stielt-K}
    m_{\widehat{\mathsf K}^{\beta_n^\star}_n}(\zeta,\xi)
    \longrightarrow
    L(\zeta,\xi)
    &\coloneqq
    \frac{1}{\zeta\xi}
    \frac{\Delta(\zeta,\xi)}{v(\zeta)v(\xi)}
    \int_0^\infty
    \frac{t}{
        \{1+tv(\zeta)\}
        \{1+tv(\xi)\}}
    \,dG(t),
\end{align}
and
\begin{align}
\label{eq:Q:stielt:def}
    m_{\mathsf Q_n}(\zeta)
    \longrightarrow
    m_{\mathsf Q}(\zeta)
    \coloneqq
    -\frac1\gamma
    \left(
        \frac{1+\zeta v(\zeta)}{\zeta v(\zeta)}
    \right).
\end{align}
Moreover, almost surely, $\mathsf Q_n\Rightarrow\mathsf Q$, where
$\mathsf Q=\mathsf Q_{\gamma,H}$ denotes the finite positive Borel
measure characterized by the Stieltjes transform in
\eqref{eq:Q:stielt:def}.
\end{lemma}

\begin{proof}
The convergences in \eqref{eq:stielt-K} and \eqref{eq:Q:stielt:def}
follow directly from results of 
\cite{DoernemannStark}; see Proposition~2.1 and Remark 2.2 as well as Theorem~3.2 and Remark 3.1.
Finally, $\mathsf Q_n\Rightarrow\mathsf Q$ follows from \eqref{eq:Q:stielt:def},
Lemma~\ref{lem:eventual-spectral-confinement}, and
Proposition~\ref{cor:weak-convergence-compact-support}.
\end{proof}

The uniform bound on the bilinear forms allows us to pass from the
two-resolvent limit in Lemma~\ref{lem:GR-Stieltjes-convergence}
to arbitrary pairs of continuous functions. Combining this with
Proposition~\ref{prop:risk-decomposition} yields the following
asymptotic prediction-risk formula.

\begin{theorem}
\label{lem:GR-bilinear-convergence}
Suppose that Assumption~\ref{ass:GR-basic} holds and that
$\|\beta_n^\star\|_2\longrightarrow\tau$ for some $\tau>0$.
Then there exists a deterministic constant $M>0$ such that:
\begin{enumerate}
    \item Almost surely,
    \[
        \operatorname{spec}(S_n)
        \subseteq
        [0,M]
    \]
    for all sufficiently large $n$.

    \item There exists a unique deterministic bounded bilinear form
    $B=B_{\gamma,H,G}$ on $C([0,M];\mathbb C)$ satisfying
    \[
        B(r_\zeta,r_\xi)
        =
        L(\zeta,\xi),
        \qquad
        \zeta,\xi\in\mathbb C^+,
    \]
    where $L=L_{\gamma,H,G}$ is defined in
    \eqref{eq:stielt-K}. Moreover, almost surely,
    \[
        \int_{[0,M]^2}
        g(s)h(t)\,
        d\widehat{\mathsf K}^{\beta_n^\star}_n(s,t)
        \longrightarrow
        B(g,h)
    \]
    for every $g,h\in C([0,M];\mathbb C)$.

    \item Almost surely, for every
    $f\in C([0,\infty))$,
    \begin{align}
    \label{eq:continuous-filter-risk-limit}
    \begin{split}
        \mathrm R_{\Sigma_n}
        \bigl(
            \hat\beta(f)
        \bigr)
        \longrightarrow{}&
        \tau^2
        B
        \bigl(
            s\mapsto1-sf(s),
            s\mapsto1-sf(s)
        \bigr)
        \\
        &\quad+
        \sigma^2\gamma
        \int_{[0,M]}
        sf(s)^2
        \,d\mathsf Q(s),
    \end{split}
    \end{align}
    where $\mathsf Q$ is the finite positive Borel measure from
    Lemma~\ref{lem:GR-Stieltjes-convergence}.
\end{enumerate}
\end{theorem}

The convergence in Theorem~\ref{lem:GR-bilinear-convergence}
is weaker than weak convergence of
$\widehat{\mathsf K}^{\beta_n^\star}_n$: it concerns functions of
product form and their finite linear combinations. Although these
functions are dense in $C([0,M]^2)$, extending the convergence to
arbitrary continuous functions by uniform approximation would require
a uniform total-variation bound. 
Nevertheless, convergence on products suffices for
\eqref{eq:continuous-filter-risk-limit}. 
This raises the question whether the limiting bilinear form $B$
can be represented by a finite signed Borel measure on $[0,M]^2$.
The answer is affirmative when $G=H$, where $G$ and $H$ are the weak limits defined in
Assumptions~\ref{ass:GR-basic-5} and~\ref{ass:GR-basic-4}.
In this case, the measure representation holds under
Assumption~\ref{ass:GR-basic} alone, without imposing the
additional regularity conditions of Theorem~\ref{thm:GR-fix} below.
More
precisely, the representing measure reduces to the finite positive
measure supported on the diagonal,
\[
    \mathsf K
    =
    D_\#\mathsf Q,
    \qquad
    D(s)=(s,s).
\]
Consequently, the limiting prediction risk takes a particularly simple
form, as recorded in the following corollary.
Its proof is deferred to Section~\ref{sec_proofs}.

\begin{corollary}\label{cor:collapse-G-equals-H}
    If we have $G = H$ in addition to the assumptions of Theorem~\ref{lem:GR-bilinear-convergence}, then 
    \begin{equation*}
        B(g,h)
        =
        \int_{[0,M]} g(s)h(s)\,d\mathsf Q(s),
    \end{equation*}
    and, almost surely, for every $f\in C([0,\infty))$, we have
    \begin{align}
\label{eq:risk-G-equals-H}
\mathrm R_{\Sigma_n}(\hat\beta(f))
\longrightarrow
\tau^2 \int_{[0,M]} \bigl(1-sf(s)\bigr)^2\,d\mathsf Q(s)
+
\sigma^2\gamma \int_{[0,M]} s f(s)^2\,d\mathsf Q(s).
\end{align}
\end{corollary}

One may consider the first term in \eqref{eq:risk-G-equals-H} as the asymptotic bias and the second as the asymptotic variance term. Thus, when the limiting measures $G$ and $H$
coincide, both the asymptotic bias and the asymptotic variance are governed
by the same one-dimensional positive measure, and the signed
two-dimensional overlap structure disappears from the final risk formula.
Notably, the condition $G=H$ also holds in the anisotropic two-block
covariance model of
Proposition~\ref{prop:GreenRomanov-TV-counterexample}, even though the
total variation of the corresponding bivariate eigenvector-overlap
measures diverges with probability tending to one. Hence, the relevant
asymptotic risk functionals can still be characterized even when weak
convergence of the underlying signed overlap measures is unavailable.

\begin{remark}[The optimal filter when $G=H$]
\label{rem:optimal-filter-G-equals-H}
Assume that $G=H$ and $\sigma>0$.
Pointwise minimization of the combined integrand in
\eqref{eq:risk-G-equals-H} shows that
\[
    f_\star(s)=\frac{1}{s+\lambda_\star},
    \qquad
    \lambda_\star\coloneqq\frac{\sigma^2\gamma}{\tau^2},
\]
minimizes the limiting prediction risk over all 
 filters $f\in C([0,\infty))$.

The same optimal regularization parameter appears in
\cite[Theorem~2.1]{DobribanWager2018} under an isotropic
random-effects model with general population covariance, and in
\cite[Theorem~5]{Hastie2022} for a deterministic signal with
isotropic covariance. Both results concern optimization within
the ridge family. 
More recently, \cite[Lemma~3.2]{Sur2026} showed that
optimally tuned ridge regression is optimal among all
spectral-shrinkage estimators in the isotropic setting;
see the discussion following Theorem~\ref{cor:ridge-optimality} and Corollary~\ref{cor:ridge-optimality2}.
Here, the condition $G=H$ allows a general limiting population
spectral distribution while retaining optimality over all
continuous spectral filters.
\end{remark}

Moreover, beyond this special case, we introduce additional sufficient conditions that allow us to represent the limiting bounded bilinear form $B$ on $C([0,M];\mathbb C)^2$ by a finite signed Borel measure on $[0,M]^2$. This representation, in turn, enables us to characterize when minimum-norm interpolation and optimally tuned ridge regression are optimal over all continuous filter functions.

First, we assume that the positive part of the asymptotic population spectrum is
separated from zero.
\begin{assumption}\label{eq:GR-main-positive-population-gap}
We assume that
\begin{equation*} 
    \tau_-
    \coloneqq
    \inf\bigl(
        \operatorname{supp}(H)\setminus\{0\}
    \bigr)
    >0,
    \qquad \text{and} \qquad
    \gamma[1-H(\{0\})]\neq1.
\end{equation*}     
\end{assumption}
As shown in
Lemma~\ref{lem:positive-spectrum-separated-zero}, these conditions
imply that the positive limiting spectrum
$\mathcal S^+\coloneqq
\mathcal S\cap(0,\infty)$, where
$\mathcal S\coloneqq\operatorname{supp}(F_{\gamma,H})$,
is separated from zero. This rules out the critical hard-edge regime in which zero may be a boundary point of the positive limiting spectrum.

Next, recall that $v$ denotes the Stieltjes
transform of $\overline F_{\gamma,H}$ defined in~\eqref{eq:companion_stielt}. By
\cite[Theorem~1.1]{SilversteinChoi1995}, the limiting
Marchenko--Pastur distribution $F_{\gamma,H}$ admits a continuous
density $f_{\gamma,H}$ on $(0,\infty)$. Moreover, the
Stieltjes transform $v$ of the companion distribution
$\overline F_{\gamma,H}$ has continuous boundary values
\[
    v_0(s)
    \coloneqq
    \lim_{\substack{z\to s\\ z\in\mathbb C^+}} v(z),
    \qquad s>0.
\]
\begin{assumption}\label{eq:GR-main-zero-density-regularity}
    We assume that
\begin{equation*}
    -\frac1{v_0(s)}
    \notin
    \operatorname{supp}(H)
\end{equation*}
for every $s\in\mathcal{S}^+$ satisfying $f_{\gamma,H}(s)=0$.
\end{assumption}
Assumption~\ref{eq:GR-main-zero-density-regularity} applies to every
$s\in\mathcal S^+$ satisfying $f_{\gamma,H}(s)=0$ and therefore covers
both nonzero boundary points of $\mathcal S$ and interior zero-density
points, or so-called \emph{cusps}. This condition is related to, but weaker than, the notion of a regular
cusp in \cite{Hachem2016}, which additionally imposes a
finite-dimensional separation condition involving the population
spectral measures $H_n$.\footnote{Such a finite-dimensional condition
is needed in \cite{Hachem2016} to control the local eigenvalue behavior
near the cusp. In the present work, we only require regularity of the
limiting Marchenko--Pastur equation, in particular to establish the
properties of $v_0$ and its H\"older continuity on $\mathcal S^+$.
Therefore, a condition formulated directly in terms of the limiting
measure $H$ is sufficient.}

Finally, we impose an analogous condition on the support of the
limiting signal-weighted spectral measure $G$.
\begin{assumption}
\label{eq:GR-main-no-critical-signal}
We assume that
\[
    -\frac1{v_0(s)}
    \notin
    \operatorname{supp}(G)
\]
for every $s\in\mathcal S^+$ satisfying
$f_{\gamma,H}(s)=0$.
\end{assumption}

We do not require the spiked-covariance structure of
\cite{GreenRomanov2025}, but instead impose the noncriticality
condition directly on $\operatorname{supp}(G)$ through
Assumption~\ref{eq:GR-main-no-critical-signal}.
Under the spiked model of \cite{GreenRomanov2025},
\[
    \operatorname{supp}(G)
    \subseteq
    \operatorname{supp}(H)
    \cup
    \{\tau_1,\ldots,\tau_{k_0}\},
\]
where $\tau_1,\ldots,\tau_{k_0}$ denote the population spikes.
Since Assumption~\ref{eq:GR-main-zero-density-regularity}
already excludes $-1/v_0(s)$ from $\operatorname{supp}(H)$
at every zero-density point $s\in\mathcal S^+$,
a violation of Assumption~\ref{eq:GR-main-no-critical-signal}
could then only occur at a population spike belonging to
$\operatorname{supp}(G)$.

Before stating the main theorem, we introduce
the following notation. Let $M$ be the deterministic constant
from Theorem~\ref{lem:GR-bilinear-convergence}, and define
\begin{align}\label{eq:def_measure_Cn}
    \mathsf C_n
    \coloneqq
    \frac{1}{\|\beta_n^\star\|_2^2}
    \sum_{i=1}^p
    \left\langle
        \Sigma_n^{1/2}\beta_n^\star,w_i
    \right\rangle^2
    \delta_{\lambda_i(S_n)}.
\end{align}
By Lemma~\ref{lem:GR-auxiliary-positive-measure} in Section~\ref{appendix_C} of the Supplementary Material,
$\mathsf C_n\Rightarrow\mathsf C$ almost surely, where
$\mathsf C$ is a deterministic finite positive measure
on $[0,M]$, with Stieltjes transform
\begin{align}
\label{eq:GR-C-transform_0}
    m_{\mathsf C}(z)
    =
    -\frac1z
    \int_{[0,\infty)}
    \frac{t}{1+t v(z)}
    \,dG(t),
    \qquad z\in\mathbb C^+.
\end{align}
With $\mathsf Q$ from
Lemma~\ref{lem:GR-Stieltjes-convergence}, set
\begin{equation}\label{eq:def_mu_main}
    \mu\coloneqq\mathsf C+\mathsf Q,
    \qquad
    \rho\coloneqq\frac{d\mathsf C}{d\mu}.
\end{equation}
We choose the $\mu$-a.e. representative of
$\rho:[0,M]\to[0,1]$ defined pointwise in the proof of Lemma~\ref{lem:GR-radon-nikodym-energy} in Section~\ref{appendix_C} of the Supplementary Material. Moreover, we write
\begin{equation}\label{eq:kappa_density_def_main}
    \kappa(s,t)\coloneqq
    \begin{cases}
        \displaystyle\frac{\rho(s)-\rho(t)}{s-t},
        & s\neq t,\\[2mm]
        0, & s=t,
    \end{cases}
    \qquad
    \ell(s)\coloneqq\int_{[0,M]}\kappa(s,t)\,d\mu(t).
\end{equation}
By Lemma~\ref{lem:GR-radon-nikodym-energy} in Section~\ref{appendix_C} of the Supplementary Material,
$\kappa\in L^1(\mu\otimes\mu)$. Hence, $\ell$ exists $\mu$-almost everywhere and $\ell\in L^1(\mu)$.

We define the candidate representing measure by
\begin{equation}\label{eq:GR-candidate-main}
    \mathsf K
    \coloneqq
    \gamma\kappa\,(\mu\otimes\mu)
    +D_\#\bigl((\rho-\gamma\ell)\mu\bigr),
\end{equation}
where $D(s)=(s,s)$, $D_\#$ denotes pushforward under $D$,
and $\kappa\,(\mu\otimes\mu)$ and $\ell\mu$ denote the signed
measures with densities $\kappa$ and $\ell$ with respect to
$\mu\otimes\mu$ and $\mu$, respectively. Thus, $\mathsf K$ is a
finite signed Borel measure.

The next theorem shows that, under the additional Assumptions
\ref{eq:GR-main-positive-population-gap} --
\ref{eq:GR-main-no-critical-signal}, the limiting expression in
\eqref{eq:continuous-filter-risk-limit} admits an integral
representation with respect to the finite signed Borel measure
$\mathsf K$ defined in
\eqref{eq:GR-candidate-main}.

\begin{theorem}
\label{thm:GR-fix}
Suppose that Assumption~\ref{ass:GR-basic} and Assumptions
\ref{eq:GR-main-positive-population-gap}--
\ref{eq:GR-main-no-critical-signal} hold. Suppose furthermore that
$\|\beta_n^\star\|_2\longrightarrow\tau$ for some $\tau>0$.
Then the measure $\mathsf{K} = \mathsf K_{\gamma,H,G}$ defined in
\eqref{eq:GR-candidate-main} is a finite signed Borel measure
supported on $[0,M]^2$, where $M$ is the constant from
Theorem~\ref{lem:GR-bilinear-convergence}, and satisfies
\begin{align}\label{eq:resolvent_kernel_int_representation}
    m_{\mathsf K}(\zeta,\xi)
    \coloneqq \int_{[0,M]^2}
    \frac{1}{(s-\zeta)(u-\xi)}
    \,d\mathsf K(s,u)=
    L(\zeta,\xi),
    \qquad
    \zeta,\xi\in\mathbb C^+,
\end{align}
where $L = L_{\gamma,H,G}$ is defined in \eqref{eq:stielt-K}. 
In particular,
\[
    B(g,h)
    =
    \int_{[0,M]^2}g(s)h(t)\,d\mathsf K(s,t),
    \qquad g,h\in C([0,M];\mathbb C).
\]
Moreover,
almost surely, for every $f\in C([0,\infty))$,
\begin{align}
\label{eq:continuous-filter-risk-limit-representation}
   \begin{split}
    \mathrm R_{\Sigma_n}
    \bigl(
        \hat\beta(f)
    \bigr)
    \longrightarrow
    \tau^2
    \int_{[0,M]^2}
    \bigl(1-sf(s)\bigr)
    \bigl(1-tf(t)\bigr)
    \,d\mathsf K(s,t)
    +
    \sigma^2\gamma
    \int_{[0,M]}
    sf(s)^2
    \,d\mathsf Q(s).
    \end{split}
\end{align}
\end{theorem}

The representation~\eqref{eq:GR-candidate-main}
splits $\mathsf K$ into an off-diagonal component
$\gamma\kappa\,(\mu\otimes\mu)$ and a diagonal component
$D_\#((\rho-\gamma\ell)\mu)$.
The former vanishes on every Borel subset of the diagonal
$\Delta=\{(s,s):s\in[0,M]\}$, whereas the latter is concentrated on
$\Delta$.
The correction term $-\gamma D_\#(\ell\mu)$ ensures that both
marginals of $\mathsf K$ coincide with $\mathsf C$.

\begin{remark}[Comparison with \cite{GreenRomanov2025}]
In \cite{GreenRomanov2025}, the authors incorrectly infer weak
convergence of the empirical overlap measure $\widehat{\mathsf K}_n^{\beta_n^\star}$ to a finite
signed measure $\mathbb K$ from pointwise convergence of their
bivariate Stieltjes transforms.
In contrast, Theorem~\ref{thm:GR-fix} establishes convergence only
against a smaller class of test functions, namely functions
$(s,u)\mapsto g(s)h(u)$ with $g,h\in C([0,M])$. This suffices to identify the limiting bilinear
form and obtain the integral representation
in~\eqref{eq:continuous-filter-risk-limit-representation}
with respect to $\mathsf K$ defined in
\eqref{eq:GR-candidate-main}. Moreover, the bivariate Stieltjes
transform $m_{\mathsf K}$ coincides with the proposed limiting
transform $\mathsf m_{\mathbb K}$ in
\cite[Lemma~7]{GreenRomanov2025}. 
Thus, although weak convergence of the empirical overlap
measures need not hold, we recover the same limiting
Stieltjes transform. Moreover, our representation
\eqref{eq:GR-candidate-main} provides a simpler description
of the corresponding finite signed measure than the
componentwise representation in
\cite[Lemma~11]{GreenRomanov2025}.
\end{remark}

For the following results, we denote the limiting prediction
risk in \eqref{eq:continuous-filter-risk-limit-representation}
by $R(f)$. That is, for $f\in C([0,\infty))$,
\[
R(f)\coloneqq
\tau^2
\int_{[0,M]^2}
\bigl(1-sf(s)\bigr)
\bigl(1-tf(t)\bigr)
\,d\mathsf K(s,t)
+
\sigma^2\gamma
\int_{[0,M]}
sf(s)^2
\,d\mathsf Q(s).
\]
We now use Theorem~\ref{thm:GR-fix} to characterize the optimality of the minimum-norm interpolator. It can be written in the spectral form
\[
    \widehat\beta_{\mathrm{MN}}
    \coloneqq
    S_n^\dagger\frac{X_n^\top y_n}{n}
    =
    X_n^\dagger y_n
    =
    \hat\beta(f_{\mathrm{MN}}),
    \qquad
    f_{\mathrm{MN}}(s)\coloneqq
    \begin{cases}
        s^{-1},&s>0,\\
        0,&s=0.
    \end{cases}
\]
The filter $f_{\mathrm{MN}}$ is not continuous at zero and therefore does
not directly fall within the class of continuous filters considered in
Theorem~\ref{thm:GR-fix}. However, the assumptions below ensure that
the positive sample eigenvalues are almost surely separated from zero
for all sufficiently large $n$. Consequently, $f_{\mathrm{MN}}$ can
then be replaced by a continuous filter that agrees with it on the
empirical spectrum. This allows us to identify the limiting risk of
$\widehat\beta_{\mathrm{MN}}$ and to characterize when it is
asymptotically optimal among continuous spectral filters.
\begin{theorem}
\label{thm:min-norm-optimality}
Suppose that the assumptions of
Theorem~\ref{thm:GR-fix} hold.
Suppose further that
\[
    \gamma>1
    \qquad\text{and}\qquad
    \inf_{n\geq1}\lambda_p(\Sigma_n)>0,
\]
where $\lambda_1(\Sigma_n)\geq\cdots\geq\lambda_p(\Sigma_n)>0$
denote the ordered eigenvalues of $\Sigma_n$.
In addition, $Z_n$ is the $n\times p$ upper-left block
of a single doubly infinite array $(z_{ij})_{i,j\geq1}$
of i.i.d.\ random variables satisfying the moment conditions in Assumption~\ref{ass:GR-basic-2}.
Then, almost surely,
\begin{equation}
\label{eq:MN-risk-limit}
    \mathrm R_{\Sigma_n}
    \bigl(\widehat\beta_{\mathrm{MN}}\bigr)
    \longrightarrow
    R(f_{\mathrm{MN}})
    \coloneqq
    \tau^2\mathsf K\bigl(\{(0,0)\}\bigr)
    +
    \sigma^2\gamma
    \int_{(0,M]}\frac1s\,d\mathsf Q(s)
    <\infty.
\end{equation}
Moreover,
\begin{equation}
\label{eq:MN-optimality-characterization}
    R(f_{\mathrm{MN}})
    =
    \inf_{f\in C([0,\infty))}R(f)
    \quad\Longleftrightarrow\quad
    \sigma^2=0
    \ \text{and}\
    G=H.
\end{equation}
\end{theorem}
Similarly, we characterize when optimally tuned ridge regression is
optimal among all continuous spectral filters, without relying on the
stronger distributional assumptions imposed in
Theorem~\ref{thm:min-norm-optimality}.

\begin{theorem}
\label{cor:ridge-optimality}
Suppose that the assumptions of
Theorem~\ref{thm:GR-fix} hold, that
$G(\{0\})<1$, and that $\sigma>0$.
For every $\lambda>0$, the ridge filter
$f_\lambda(s)\coloneqq(s+\lambda)^{-1}, s\geq 0$, satisfies
\begin{equation}
\label{eq:ridge-optimality-characterization}
    R(f_\lambda)
    =
    \inf_{f\in C([0,\infty))}R(f)
    \quad\Longleftrightarrow\quad
    \frac{\left.G\right|_{(0,\infty)}}{1-G(\{0\})}
    =
    \frac{\left.H\right|_{(0,\infty)}}{1-H(\{0\})}
    \quad\text{and}\quad
    \lambda
    =
    \frac{\sigma^2\gamma\{1-H(\{0\})\}}
    {\tau^2\{1-G(\{0\})\}}.
\end{equation}
\end{theorem}
Finally, while our analysis considers the prediction risk conditional
on $X$, the recent work of \cite{Sur2026} considers the
prediction risk conditional on the entire training sample
$(X,y)$. For $f\in C([0,\infty))$, this risk is given by 
\[
\mathbb E\left[
\bigl(x_0^\top\hat\beta(f)-y_0\bigr)^2
\,\middle|\,X,y
\right]
=
(\hat\beta(f)-\beta_n^\star)^\top
\Sigma_n(\hat\beta(f)-\beta_n^\star)
+\sigma^2,
\]
and the corresponding excess prediction risk by
\begin{align} \label{eq:def-Ln}
\mathcal L_n(f)
\coloneqq
(\hat\beta(f)-\beta_n^\star)^\top
\Sigma_n(\hat\beta(f)-\beta_n^\star).
\end{align}
The following corollary shows that, for every
$f\in C([0,\infty))$, $\mathcal L_n(f)$ converges
almost surely to the same deterministic limit $R(f)$
as the prediction risk $\mathrm R_{\Sigma_n}(\hat\beta(f))$, which is conditional only on $X$.
\begin{corollary}\label{cor:ridge-optimality2}
Suppose that the assumptions of
Theorem~\ref{cor:ridge-optimality} hold.
Assume additionally that
$\mathbb E|\varepsilon_1|^{4+\delta}<\infty$
for some $\delta>0$. Then,  for every $f\in C([0,\infty))$, 
\begin{align}\label{eq:risk_cond_XY}
 \mathcal{L}_n(f) \longrightarrow R(f)\qquad \text{almost surely}.
\end{align}
In particular, the almost-sure limit of $\mathcal L_n(f)$
satisfies the ridge-optimality characterization
\eqref{eq:ridge-optimality-characterization}.
\end{corollary}

Theorem~\ref{thm:min-norm-optimality} provides necessary and
sufficient conditions for minimum-norm interpolation to be
optimal among all continuous spectral filters. Under its
assumptions, interpolation is strictly suboptimal whenever
$\sigma^2>0$, whereas in the noiseless case $\sigma^2=0$, it is optimal
if and only if $G=H$.
Theorem~\ref{cor:ridge-optimality} establishes the
corresponding characterization for ridge regression. Corollary~\ref{cor:ridge-optimality2} extends this
characterization to the prediction risk conditional on the
entire training data $(X,y)$.
In particular, ridge regression can attain the continuous-filter
optimum even for non-isotropic population covariances, provided
the normalized restrictions of $G$ and $H$ to $(0,\infty)$
coincide. Their respective masses at zero determine the
optimal regularization parameter.
When this condition fails, no choice of the ridge parameter
can attain the optimal limiting risk among continuous spectral
filters.

In a related work, \cite{Sur2026} establish the optimality of optimally tuned ridge regression
among all continuous spectral shrinkage estimators in the isotropic case
(cf. Lemma~3.2 therein). Under our assumptions, this optimality statement
is recovered as the special case $G=H=\delta_\theta$, $\theta>0$,
of Corollary~\ref{cor:ridge-optimality2}. Moreover, \cite{Sur2026} show that $\mathcal L_n(f)$
converges almost surely to a deterministic limit,
which we denote by $\mathcal L(f)$, for every continuous
filter $f$, and that self-distillation with suitably
chosen weights approximates $\inf_f \mathcal L(f)$
(cf. Theorem 3.8 therein). Going beyond this existence result, we give in \eqref{eq:continuous-filter-risk-limit-representation}, under additional
regularity assumptions, an explicit integral formula for $\mathcal L(f)=R(f)$
in terms of the measures $\mathsf K$ and $\mathsf Q$.

Our results should also be distinguished from benign overfitting
\cite{BartlettEtAl2020,TsiglerBartlett2023}. An interpolating
estimator may achieve good predictive performance without
minimizing the limiting risk among continuous spectral filters.

A natural direction for future research is to quantify the
improvements achievable through spectral filtering and to
develop data-driven procedures approaching the corresponding
oracle risk. This could be pursued by estimating the relevant
spectral measures $H$ and $G$ together with $\sigma^2/\tau^2$,
or by selecting spectral filters directly through estimates of their prediction risk.

\section{Proofs of the main results in Sections~\ref{sec:application}
and~\ref{Out-of-sample prediction risk}}\label{sec_proofs}

In this section, we prove Proposition~\ref{prop:GreenRomanov-TV-counterexample},
Theorems~\ref{lem:GR-bilinear-convergence},
\ref{thm:GR-fix}, and~\ref{thm:min-norm-optimality},
and Corollaries~\ref{cor:collapse-G-equals-H}
and~\ref{cor:ridge-optimality}. Technical details for Section~\ref{sec:proof-counterexample}
are given in Section~\ref{appendix_B},
for Section~\ref{sec:proof-bilinear}
in Sections~\ref{appendix_C} and~\ref{appendix_D},
and for Section~\ref{sec:proof-overlap-measure}
in Section~\ref{appendix_C} of the Supplementary Material.

\subsection{Proof of Proposition~\ref{prop:GreenRomanov-TV-counterexample}} \label{sec:proof-counterexample}

The proof of Proposition \ref{prop:GreenRomanov-TV-counterexample} mainly relies on
Lemmas~\ref{lem:Gaussian-block-spherical}--\ref{lem:block-spherical-deterministic-transfer} stated below. The first lemma gives the desired lower bound on the total variation
of $\widehat{\mathsf K}_n^{\overline u_n}$, where $\overline u_n$
is a block-spherical random signal. The second lemma shows that replacing $\widehat{\mathsf K}_n^{\overline u_n}$
by $\widehat{\mathsf K}_n^{\beta_n^\star}$ preserves the distribution
of the total variation norm.
\begin{lemma}\label{lem:Gaussian-block-spherical}
Under the assumptions of
Proposition~\ref{prop:GreenRomanov-TV-counterexample}, let
$g_n\sim\mathcal N(0,I_p)$ be independent of $Z_n$.
Write $g_n$ in terms of its first and last $p/2$ coordinates as
\[
g_n
=
\begin{pmatrix}
g_n^{(1)}
\\
g_n^{(2)}
\end{pmatrix},
\qquad
g_n^{(1)}\sim\mathcal{N}(0,I_{p/2}), \,g_n^{(2)}\sim\mathcal{N}(0,I_{p/2}),
\]
and define
\[
\overline u_n
=
\begin{pmatrix}
\displaystyle
\frac1{\sqrt2}
\frac{g_n^{(1)}}{\|g_n^{(1)}\|_2}
\\[3mm]
\displaystyle
\frac1{\sqrt2}
\frac{g_n^{(2)}}{\|g_n^{(2)}\|_2}
\end{pmatrix}.
\]
Then, for every fixed $\varepsilon\in(0,1/2)$, there exists a
constant $c_4>0$ such that
\begin{align*}\label{eq:block-spherical-TV-bound}
\mathbb P\left(
\|\widehat{\mathsf K}_n^{\overline u_n}\|_{\mathrm{TV}}
\geq
c_4 p^{1/2-\varepsilon}
\right)
\longrightarrow1.
\end{align*}
\end{lemma}

\begin{lemma}
\label{lem:block-spherical-deterministic-transfer}
Let $\overline u_n$ be the block-spherical random vector defined in
Lemma~\ref{lem:Gaussian-block-spherical}. Then
\[
    \|\widehat{\mathsf K}_n^{\beta_n^\star}\|_{\mathrm{TV}}
    \overset{d}=
    \|\widehat{\mathsf K}_n^{\overline u_n}\|_{\mathrm{TV}},
\]
where $\overset{d}{=}$ denotes equality in distribution.
\end{lemma}

\begin{proof}[Proof of Proposition~\ref{prop:GreenRomanov-TV-counterexample}]
By Lemma~\ref{lem:verification_GR_assumptions} in Section~\ref{appendix_B} of the Supplementary Material,
Assumptions~1--6 of \cite{GreenRomanov2025} are satisfied.
Fix $\varepsilon\in(0,1/2)$ and let $\overline u_n$
be the random vector defined in
Lemma~\ref{lem:Gaussian-block-spherical}.
By that lemma, there exists a constant $c>0$ such that
\[
    \mathbb P\left(
    \|\widehat{\mathsf K}_n^{\overline u_n}\|_{\mathrm{TV}}
    \geq c p^{1/2-\varepsilon}
    \right)
    \longrightarrow1.
\]
Lemma~\ref{lem:block-spherical-deterministic-transfer} gives
\[
    \|\widehat{\mathsf K}_n^{\beta_n^\star}\|_{\mathrm{TV}}
    \overset{d}=
    \|\widehat{\mathsf K}_n^{\overline u_n}\|_{\mathrm{TV}},
\]
which proves the claim.
\end{proof}

\begin{remark}[Relation to the equivariance argument in
\cite{Dicker2016}]
The block-orthogonal invariance argument in
Lemma~\ref{lem:block-spherical-deterministic-transfer}
is closely related to \cite[Proposition~3]{Dicker2016}. In the isotropic Gaussian model considered there, invariance under the
full orthogonal group $\mathsf O(p)$, together with orthogonal
equivariance of the estimator, implies that its risk is rotationally
invariant in the signal parameter, and hence, constant on each Euclidean
sphere. This relates the corresponding minimax problem to a Bayes
problem under the uniform spherical prior.

In the present anisotropic model, the symmetry group is reduced to
$\mathsf O(p/2)\times\mathsf O(p/2)$, since only block-diagonal
orthogonal matrices commute with $\Sigma_n$. The orbit of
$\beta_n^\star$ under this group is the product of
two spheres of radius $1/\sqrt2$, and the Haar distribution on this
orbit is the distribution of $\overline u_n$. 
\end{remark}
The proofs of Lemmas~\ref{lem:Gaussian-block-spherical}
and~\ref{lem:block-spherical-deterministic-transfer} are given in
Section~\ref{appendix_B} of the Supplementary Material. We outline the main idea of the proof of
Lemma~\ref{lem:Gaussian-block-spherical}. First, we establish a lower
bound on $\|\widehat{\mathsf K}_n^{u_n}\|_{\mathrm{TV}}$, where
$u_n=g_n/\|g_n\|_2$ is uniformly distributed on the unit sphere.
Working with $u_n$ is convenient in our setting. To see why, note that
$\|\widehat{\mathsf K}_n^{u_n}\|_{\mathrm{TV}}$ depends on $u_n$
through the coordinates
$\langle u_n,w_i\rangle=\langle g_n,w_i\rangle/\|g_n\|_2$.
Conditional on $S_n$, the numerators $\langle g_n,w_i\rangle$ are
independent standard Gaussian variables.

Second, we show that the total variation norms of
$\widehat{\mathsf K}_n^{u_n}$ and $\widehat{\mathsf K}_n^{\overline u_n}$
differ by $O_{\mathbb P}(1)$. Thus, the lower bound of order
$p^{1/2-\varepsilon}$ obtained in the first step carries over to
$\|\widehat{\mathsf K}_n^{\overline u_n}\|_{\mathrm{TV}}$.

To explain the first step, we state two auxiliary lemmas. The first
controls the terms $\langle w_i,P_nw_j\rangle$ for distinct indices
$i,j$ away from the spectral edges. These terms appear in the expression
for $\|\widehat{\mathsf K}_n^{u_n}\|_{\mathrm{TV}}$. The second lemma
uses these estimates to obtain the lower bound on this norm.

\begin{lemma}\label{lem:bulk-overlap-control}
Under the assumptions of
Proposition~\ref{prop:GreenRomanov-TV-counterexample}, consider the
spectral decomposition in \eqref{eq:sam-cov-decomposition} and set $c_{ij}
\coloneqq
\langle w_i,P_nw_j\rangle$, $1\leq i,j\leq p$.
There exist constants $\eta\in(0,1/2)$ and $c_1>0$ such that, with
\[
J_n
\coloneqq
\left\{
\lceil\eta p\rceil,\ldots,
\lfloor(1-\eta)p\rfloor
\right\},
\]
we have
\begin{align}\label{eq:definition_A_n}
\mathbb P\left(
\sum_{\substack{i,j\in J_n\\i\neq j}}
c_{ij}^2
\geq
c_1p
\right)
\longrightarrow1.
\end{align}
Moreover, for every fixed \(\varepsilon\in(0,1/2)\),
\begin{align}
\label{eq:outline-step2}
\mathbb P\left(
\max_{\substack{i,j\in J_n\\i\neq j}}
|c_{ij}|
\leq
p^{-1/2+\varepsilon}
\right)
\longrightarrow1.
\end{align}
\end{lemma}

\begin{lemma}\label{lem:Gaussian-retention}
Under the assumptions and notation of
Lemma~\ref{lem:bulk-overlap-control}, let $g_n$ be the Gaussian vector
from Lemma~\ref{lem:Gaussian-block-spherical}, and set
\[
u_n
\coloneqq
\frac{g_n}{\|g_n\|_2}.
\]
Then, for every fixed $\varepsilon\in(0,1/2)$, there exists a
constant $c_3>0$ such that
\begin{align*}\label{eq:outline-step4}
\mathbb P\left(
\|\widehat{\mathsf K}_n^{u_n}\|_{\mathrm{TV}}
\geq
c_3 p^{1/2-\varepsilon}
\right)
\longrightarrow1.
\end{align*}
\end{lemma}
The proofs of Lemmas \ref{lem:bulk-overlap-control} and \ref{lem:Gaussian-retention} 
are postponed to Section~\ref{appendix_B} of the Supplementary Material.

\subsection{Proof of Theorem~\ref{lem:GR-bilinear-convergence}} \label{sec:proof-bilinear}

We first record an auxiliary result on the convergence of bilinear
forms. Recall that the norm of a bilinear form $B:X^2\to\mathbb F$
is given by
\begin{equation*}
    \|B\| = \sup_{\substack{x,y\in X\\\|x\|,\|y\|\leq1}} |B(x,y)|,
\end{equation*}
and that $B$ is continuous if and only if $\|B\|<\infty$.

\begin{lemma}\label{lem:abstract-bilinear-extension}
    Let $M>0$, and let 
    \begin{equation*}
        B_n: C([0,M]; \mathbb{C})^2\to \mathbb{C}
    \end{equation*}
    for $n\geq 1$, be continuous bilinear forms with $C_B \coloneqq\sup_{n\geq 1 } \|B_n\| < \infty$. 
    If $B_n(r_\zeta, r_\xi)$ converges to $L(\zeta,\xi)$ for every $\zeta, \xi\in \mathbb{C}^+$, where $r_z(s) = (s-z)^{-1}$, $s\in [0,M]$ denotes the resolvent kernel of $z\in \mathbb{C}^+$, then there exists a unique bilinear form $B:C([0,M]; \mathbb{C})^2\to \mathbb{C}$ such that $\|B\|\leq  C_B$, 
    \[
        B(r_\zeta,r_\xi)
        =
        L(\zeta,\xi),
        \qquad
        \zeta,\xi\in\mathbb C^+,
    \]
    and 
    \begin{equation*}
        B_n(g,h)\to B(g,h)
    \end{equation*}
    for all $g,h \in C([0,M];\mathbb{C})$. 
\end{lemma}
The proof of Lemma~\ref{lem:abstract-bilinear-extension} is straightforward, and hence deferred to Section~\ref{appendix_D} of the Supplementary Material.

\begin{proof}[Proof of Theorem~\ref{lem:GR-bilinear-convergence}]
By Lemma~\ref{lem:eventual-spectral-confinement}, there exists a
deterministic $M>0$ such that, almost surely,
\begin{equation}
\label{eq:GR-bilinear-support-bound}
    \operatorname{spec}(S_n)\subseteq[0,M]
\end{equation}
for all sufficiently large $n$. Since
$\widehat{\mathsf K}^{\beta_n^\star}_n$ is supported on
$\operatorname{spec}(S_n)^2$, it follows that
\[
    \operatorname{supp}
    \bigl(
        \widehat{\mathsf K}^{\beta_n^\star}_n
    \bigr)
    \subseteq
    [0,M]^2
\]
eventually almost surely. This proves the first assertion.

In what follows, we work on the intersection of the probability-one
event on which \eqref{eq:GR-bilinear-support-bound} holds eventually
and the probability-one event from
Lemma~\ref{lem:GR-Stieltjes-convergence} on which, simultaneously
for every $\zeta,\xi\in\mathbb C^+$,
\[
    m_{\widehat{\mathsf K}^{\beta_n^\star}_n}(\zeta,\xi)
    \longrightarrow
    L(\zeta,\xi)
    \qquad\textnormal{and}\qquad
    \mathsf Q_n\Rightarrow\mathsf Q.
\]
Denote this event by $\Omega_0$ and fix $\omega\in\Omega_0$.
There exists $N=N_\omega$ such that
\[
    \operatorname{spec}(S_n)\subseteq[0,M],
    \qquad
    n\geq N.
\]
Let $n\geq N$.
For $g,h\in C([0,M];\mathbb C)$, define
\begin{align}
\label{eq:def_B_n}
    B_n(g,h)
    \coloneqq
    \int_{[0,M]^2}
    g(s)h(t)\,
    d\widehat{\mathsf K}^{\beta_n^\star}_n(s,t)
    =
    \frac1{\|\beta_n^\star\|_2^2}
    (\beta_n^\star)^\top
    g(S_n)\Sigma_n h(S_n)
    \beta_n^\star,
\end{align}
where the last equality follows from the definition of
$\widehat{\mathsf K}^{\beta_n^\star}_n$.
By \eqref{eq:bound_B_n} and
Assumption~\ref{ass:GR-basic-3},
\[
    \sup_{n\geq N}\|B_n\|
    \leq \sup_{n\geq1}\|\Sigma_n\|_2
    <\infty.
\]

For $\zeta,\xi\in\mathbb C^+$, the spectral functional calculus gives
\[
    r_\zeta(S_n)
    =
    (S_n-\zeta I_p)^{-1}
    =
    R(\zeta),
\]
and hence, by \eqref{eq:def_B_n},
\[
\begin{aligned}
    B_n
    \bigl(
        r_\zeta,
        r_\xi
    \bigr)
    &=
    \frac1{\|\beta_n^\star\|_2^2}
    (\beta_n^\star)^\top
    R(\zeta)\Sigma_nR(\xi)
    \beta_n^\star
    =
    m_{\widehat{\mathsf K}^{\beta_n^\star}_n}(\zeta,\xi).
\end{aligned}
\]
Lemma~\ref{lem:GR-Stieltjes-convergence} therefore yields
\[
    B_n
    \bigl(
        r_\zeta,
        r_\xi
    \bigr)
    \longrightarrow
    L(\zeta,\xi)
\]
for every $\zeta,\xi\in\mathbb C^+$.
Applying Lemma~\ref{lem:abstract-bilinear-extension} to
$\{B_n:n\geq N\}$ yields a unique bounded bilinear form
$B_{\gamma,H,G}$ such that
\[
    B_n(g,h)
    \longrightarrow
    B_{\gamma,H,G}(g,h)
\]
for every $g,h\in C([0,M];\mathbb C)$ and
\[
    B_{\gamma,H,G}
    \bigl(
        r_\zeta,
        r_\xi
    \bigr)
    =
    L(\zeta,\xi).
\]
Since $L$ is deterministic, the uniqueness assertion in
Lemma~\ref{lem:abstract-bilinear-extension} shows that the resulting
bilinear form does not depend on the choice of
$\omega\in\Omega_0$. In what follows, we write $B\coloneqq B_{\gamma,H,G}$.
This concludes the proof of the second assertion.

For the third one, fix a continuous function
$f:[0,\infty)\to\mathbb R$. Since
$s\mapsto1-sf(s)$ is continuous on $[0,M]$, the preceding
convergence of the bilinear forms gives
\begin{align}
\label{eq:GR-bias-product-convergence}
\int_{[0,M]^2}
    \bigl(1-sf(s)\bigr)
    \bigl(1-tf(t)\bigr)
    \,d\widehat{\mathsf K}^{\beta_n^\star}_n(s,t)
\longrightarrow
    B
    \bigl(
        s\mapsto1-sf(s),
        s\mapsto1-sf(s)
    \bigr).
\end{align}
Since $\mathsf Q_n$ is supported on $[0,M]$ for all
sufficiently large $n$, its weak limit $\mathsf Q$ from
Lemma~\ref{lem:GR-Stieltjes-convergence} is also supported on
$[0,M]$. Since $s\mapsto sf(s)^2$ is continuous on
$[0,M]$, the weak convergence in
Lemma~\ref{lem:GR-Stieltjes-convergence} therefore gives
\begin{align}
\label{eq:GR-variance-measure-convergence}
    \int_{[0,M]}
    sf(s)^2
    \,d\mathsf Q_n(s)
    \longrightarrow
    \int_{[0,M]}
    sf(s)^2
    \,d\mathsf Q(s).
\end{align}
By Proposition~\ref{prop:risk-decomposition}, for every sufficiently
large $n$,
\begin{align*}
    \mathrm R_{\Sigma_n}
    \bigl(
        \hat\beta(f)
    \bigr)
    &=
    \|\beta_n^\star\|_2^2
    \int_{[0,M]^2}
    \bigl(1-sf(s)\bigr)
    \bigl(1-tf(t)\bigr)
    \,d\widehat{\mathsf K}^{\beta_n^\star}_n(s,t)
    \\
    &\quad+
    \sigma^2\gamma_n
    \int_{[0,M]}
    sf(s)^2
    \,d\mathsf Q_n(s).
\end{align*}
Since $\|\beta_n^\star\|_2^2
    \longrightarrow
    \tau^2$ and $\gamma_n
    \longrightarrow
    \gamma$,
combining
\eqref{eq:GR-bias-product-convergence} and
\eqref{eq:GR-variance-measure-convergence} yields
\begin{align*}
    \mathrm R_{\Sigma_n}
    \bigl(
        \hat\beta(f)
    \bigr)
    \longrightarrow{}&
    \tau^2
    B
    \left(
        s\mapsto1-sf(s),
        s\mapsto1-sf(s)
    \right)
    +
    \sigma^2\gamma
    \int_{[0,M]}
    sf(s)^2
    \,d\mathsf Q(s).
\end{align*}
Since $\omega\in\Omega_0$ was arbitrary and all the preceding
convergences hold on the same probability-one event for every
continuous $f:[0,\infty)\to\mathbb R$, this proves
\eqref{eq:continuous-filter-risk-limit} and completes the proof.
\end{proof}

\subsection{Proof of Corollary~\ref{cor:collapse-G-equals-H}}

\begin{proof}[Proof of Corollary~\ref{cor:collapse-G-equals-H}]
By Lemma~\ref{lem:GR-Stieltjes-convergence},
$\mathsf Q$ is a finite positive Borel measure supported on $[0,M]$.
Define
\[
    \widetilde B(g,h)
    \coloneqq
    \int_{[0,M]^2}
    g(s)h(t)\,
    d(D_\#\mathsf Q)(s,t)
    =
    \int_{[0,M]}
    g(s)h(s)\,d\mathsf Q(s).
\]
Then $\widetilde B$ is a bounded bilinear form, since
\[
    |\widetilde B(g,h)|
    \leq
    \mathsf Q([0,M])
    \|g\|_{\infty}
    \|h\|_{\infty}.
\]
We show that $B$ from Theorem~\ref{lem:GR-bilinear-convergence} and $\widetilde B$ agree on every pair of resolvent
kernels. For $z\in\mathbb C^+$, the fixed-point equation for $v$,
\cite[Equation~(1.3)]{SilversteinChoi1995}, gives
\[
    v(z)
    =
    \left\{
        -z
        +
        \gamma
        \int_0^\infty
        \frac{t}{1+tv(z)}
        \,dH(t)
    \right\}^{-1}.
\]
Since $H$ is a probability measure, it follows that
\[
    \int_0^\infty
    \frac1{1+tv(z)}
    \,dH(t)
    =
    1-\frac{1+zv(z)}{\gamma}.
\]
Hence, for $\zeta,\xi\in\mathbb C^+$ with $\zeta\neq\xi$,
\begin{align*}
    &\{v(\zeta)-v(\xi)\}
    \int_0^\infty
    \frac{t}{
        \{1+tv(\zeta)\}
        \{1+tv(\xi)\}}
    \,dH(t)
    \\
    &\qquad=
    \int_0^\infty
    \left\{
        \frac1{1+tv(\xi)}
        -
        \frac1{1+tv(\zeta)}
    \right\}
    \,dH(t)
    =
    \frac{
        \zeta v(\zeta)-\xi v(\xi)
    }{\gamma}.
\end{align*}
Using this identity, \eqref{eq:stielt-K},
\eqref{eq:Q:stielt:def}, and 
\[
    \frac1{(s-\zeta)(s-\xi)}
    =
    \frac1{\zeta-\xi}
    \left(
        \frac1{s-\zeta}
        -
        \frac1{s-\xi}
    \right),
\]
we obtain
\begin{align*}
    B
    \bigl(
        r_\zeta,
        r_\xi
    \bigr)
    =
    L(\zeta,\xi)=
    \frac{
        \zeta v(\zeta)-\xi v(\xi)
    }{
        \gamma\zeta\xi
        v(\zeta)v(\xi)
        (\zeta-\xi)
    }
    =
    \frac{
        m_{\mathsf Q}(\zeta)
        -
        m_{\mathsf Q}(\xi)
    }{\zeta-\xi}
    =
    \widetilde B
    \bigl(
        r_\zeta,
        r_\xi
    \bigr).
\end{align*}
The map
\[
    \zeta
    \mapsto
    r_\zeta
\]
is continuous from $\mathbb C^+$ to
$C([0,M];\mathbb C)$ with respect to the supremum norm. Since
$B$ and $\widetilde B$ are bounded bilinear forms, the preceding
identity extends by continuity to $\zeta=\xi$. Thus, the two bilinear
forms agree on every pair of resolvent kernels. By the uniqueness
assertion in Theorem~\ref{lem:GR-bilinear-convergence},
$B=\widetilde B$. Consequently, for every
$g,h\in C([0,M];\mathbb C)$,
\begin{align}
\label{eq:B-HH-diagonal-representation}
    B(g,h)
    =
    \int_{[0,M]^2}
    g(s)h(t)\,
    d(D_\#\mathsf Q)(s,t)
    =
    \int_{[0,M]}
    g(s)h(s)\,d\mathsf Q(s).
\end{align}
Finally, applying Theorem~\ref{lem:GR-bilinear-convergence}
and \eqref{eq:B-HH-diagonal-representation} with
$g(s)=h(s)=1-sf(s)$ yields, almost surely,
\begin{align*}
    \mathrm R_{\Sigma_n}
    \bigl(
        \hat\beta(f)
    \bigr)
    \longrightarrow
    \tau^2
    \int_{[0,M]}
    \bigl(1-sf(s)\bigr)^2
    \,d\mathsf Q(s)
    +
    \sigma^2\gamma
    \int_{[0,M]}
    sf(s)^2
    \,d\mathsf Q(s).
\end{align*}
This proves \eqref{eq:risk-G-equals-H}.
\end{proof}

\subsection{Proof of Theorem~\ref{thm:GR-fix}} \label{sec:proof-overlap-measure}

\begin{lemma}
\label{lem:GR-candidate-measure}
Suppose that Assumption~\ref{ass:GR-basic} and
Assumptions~\ref{eq:GR-main-positive-population-gap},
\ref{eq:GR-main-zero-density-regularity}, and
\ref{eq:GR-main-no-critical-signal} hold.
Then $\mathsf K$ is a finite signed Borel measure
supported on $[0,M]^2$ satisfying
\begin{align}
\label{eq:GR-candidate-TV}
    \|\mathsf K\|_{\mathrm{TV}}
    \leq
    \mathsf C([0,M])
    +
    2\gamma\mathcal E
    <\infty,
\end{align}
where $\mathcal E \coloneqq \int_{[0,M]^2}
    |\kappa(s,t)|
    \,d\mu(s)d\mu(t).$
Moreover,
\begin{equation}
\label{eq:GR-candidate-transform}
    m_{\mathsf K}(\zeta,\xi)
    =
    L(\zeta,\xi),
    \qquad
    \zeta,\xi\in\mathbb C^+,
\end{equation}
where $L$ is defined in
Lemma~\ref{lem:GR-Stieltjes-convergence},
and the limiting bilinear form $B$ from
Theorem~\ref{lem:GR-bilinear-convergence}
admits the representation
\begin{equation}
\label{eq:GR-bilinear-measure-representation}
    B(g,h)
    =
    \int_{[0,M]^2}
    g(s)h(t)
    \,d\mathsf K(s,t)
\end{equation}
for every $g,h\in C([0,M];\mathbb C)$.
\end{lemma}
The proof is deferred to Section~\ref{appendix_C} of the Supplementary Material.

\begin{proof}[Proof of Theorem~\ref{thm:GR-fix}]
By Lemma~\ref{lem:GR-candidate-measure}, $\mathsf K$ is a
finite signed Borel measure supported on $[0,M]^2$ satisfying
\eqref{eq:resolvent_kernel_int_representation}, and the limiting
bilinear form $B$ admits the representation
\eqref{eq:GR-bilinear-measure-representation}.
Substituting this representation into
\eqref{eq:continuous-filter-risk-limit} from
Theorem~\ref{lem:GR-bilinear-convergence} yields
\eqref{eq:continuous-filter-risk-limit-representation}.
\end{proof}

\subsection{Proof of Theorem~\ref{thm:min-norm-optimality}}

\begin{proof}[Proof of Theorem~\ref{thm:min-norm-optimality}]
First, we establish the asymptotic risk of $\widehat\beta_{\mathrm{MN}}$.
Set
\[
    c_\Sigma
    \coloneqq
    \inf_{n\geq1}\lambda_{\min}(\Sigma_n)>0.
\]
By Assumption~\ref{ass:GR-basic-2} and the additional i.i.d.\
assumption of the theorem, \cite[Theorem~2]{BaiYin1993} yields
\[
    \lambda_{\min}\left(\frac1n Z_nZ_n^\top\right)
    \overset{\mathrm{a.s.}}{\longrightarrow}
    (1-\sqrt{\gamma})^2>0.
\]
Moreover,
\[
    \lambda_{\min}\left(\frac1nX_nX_n^\top\right)
    =
    \lambda_{\min}\left(\frac1nZ_n\Sigma_n Z_n^\top\right)
    \geq
    c_\Sigma
    \lambda_{\min}\left(\frac1nZ_nZ_n^\top\right).
\]
Hence, almost surely, $X_n$ has full row rank for all sufficiently
large $n$, and the nonzero eigenvalues of
$S_n=n^{-1}X_n^\top X_n$ are uniformly bounded away from zero.
Together with the upper spectral bound from
Theorem~\ref{lem:GR-bilinear-convergence}, this yields, for all
sufficiently large $n$,
\begin{equation}
\label{eq:MN-empirical-spectral-gap}
    \operatorname{spec}(S_n)
    \subseteq
    E
    \coloneqq
    \{0\}\cup[\delta,M],
\end{equation}
for some deterministic $\delta>0$.

We next transfer this restriction to $\mathsf Q$ and $\mathsf C$.
By \eqref{eq:def_measure_Cn} and
\eqref{eq:definition_Q_n}, both $\mathsf C_n$ and
$\mathsf Q_n$ are supported on $\operatorname{spec}(S_n)$.
Thus, almost surely,
\[
    \mathsf C_n((0,\delta))
    =
    \mathsf Q_n((0,\delta))
    =
    0,
\]
for all sufficiently large $n$. Moreover,
Lemma~\ref{lem:GR-auxiliary-positive-measure} and
Lemma~\ref{lem:GR-Stieltjes-convergence} give
$\mathsf C_n\Rightarrow\mathsf C$ and
$\mathsf Q_n\Rightarrow\mathsf Q$ almost surely.
Applying the Portmanteau theorem to the open interval
$(0,\delta)$ yields
\[
    0
    \leq
    \mathsf C((0,\delta))
    \leq
    \liminf_{n\to\infty}\mathsf C_n((0,\delta))
    =
    0,
\]
and the same argument gives $\mathsf Q((0,\delta))=0$.
Since both limiting measures are supported on $[0,M]$,
it follows from \eqref{eq:def_mu_main} that
\begin{equation}
\label{eq:MN-limiting-support}
    \operatorname{supp}(\mathsf C)
    \cup
    \operatorname{supp}(\mathsf Q)
    \subseteq E,
    \qquad
    \operatorname{supp}(\mu)\subseteq E.
\end{equation}
Consequently,
$\kappa\,(\mu\otimes\mu)$ is concentrated on $E\times E$,
whereas $(\rho-\gamma\ell)\mu$ is concentrated on $E$.
Its pushforward under $D$ is therefore concentrated on
$D(E)\subseteq E\times E$.
Definition \eqref{eq:GR-candidate-main} now implies
\begin{equation}
\label{eq:MN-K-support}
    \operatorname{supp}(\mathsf K)
    \subseteq
    E\times E.
\end{equation}

Define
\[
    \widetilde f_{\mathrm{MN}}(s)
    \coloneqq
    \begin{cases}
        s/\delta^2, & 0\leq s\leq\delta,\\
        s^{-1}, & s\geq\delta,
    \end{cases}
\]
and note that this function is continuous on $[0,\infty)$.
By \eqref{eq:MN-empirical-spectral-gap},
\[
    \widetilde f_{\mathrm{MN}}(S_n)=S_n^\dagger,
\]
for all sufficiently large $n$, almost surely. Moreover,
\[
    \widehat\beta(\widetilde f_{\mathrm{MN}})
    =
    \widetilde f_{\mathrm{MN}}(S_n)
    \frac{X_n^\top y_n}{n}
    =
    \widehat\beta_{\mathrm{MN}}
\]
for all sufficiently large $n$, almost surely.
We may therefore apply
\eqref{eq:continuous-filter-risk-limit-representation}
to $\widetilde f_{\mathrm{MN}}$.
Write $R(f)$ for the limiting functional in
\eqref{eq:continuous-filter-risk-limit}.
Since this functional depends only on the restriction of
$f$ to $[0,M]$, we also use it for
$f\in C([0,M])$. To evaluate the bias term, set
\[
    \varphi_{\mathrm{MN}}(s)
    \coloneqq
    1-s\widetilde f_{\mathrm{MN}}(s),
    \qquad s\in[0,M].
\]
We have $\varphi_{\mathrm{MN}}(0)=1$, whereas
$\varphi_{\mathrm{MN}}(s)=1-s/s=0$ for $s\in[\delta,M]$.
Thus,
\begin{equation}
\label{eq:MN-residual-on-support}
    \varphi_{\mathrm{MN}}(s)
    =\mathbf 1_{\{0\}}(s),
    \qquad s\in E.
\end{equation}
Using the integral representation of $B$, together with
\eqref{eq:MN-K-support} and
\eqref{eq:MN-residual-on-support}, we obtain
\begin{align*}
    B(\varphi_{\mathrm{MN}},\varphi_{\mathrm{MN}})
    &=
    \int_{E\times E}
    \varphi_{\mathrm{MN}}(s)\varphi_{\mathrm{MN}}(t)
    \,d\mathsf K(s,t)=
    \mathsf K\bigl(\{(0,0)\}\bigr).
\end{align*}
For the variance term, the integrand is zero at $s=0$
and equals $1/s$ on $[\delta,M]$. By
\eqref{eq:MN-limiting-support},
\begin{align*}
    \int_{[0,M]}
    s\widetilde f_{\mathrm{MN}}(s)^2\,d\mathsf Q(s)
    =
    \int_{[\delta,M]}
    s\left(\frac1s\right)^2\,d\mathsf Q(s)
    =
    \int_{(0,M]}\frac1s\,d\mathsf Q(s).
\end{align*}
Consequently, almost surely,
\[
    \mathrm R_{\Sigma_n}
    \bigl(\widehat\beta_{\mathrm{MN}}\bigr)
    \longrightarrow
    \tau^2\mathsf K\bigl(\{(0,0)\}\bigr)
    +
    \sigma^2\gamma
    \int_{(0,M]}\frac1s\,d\mathsf Q(s).
\]
The limit is finite, since $\mathsf K$ is a finite signed Borel measure
and
\[
    \int_{(0,M]}\frac1s\,d\mathsf Q(s)
    =
    \int_{[\delta,M]}\frac1s\,d\mathsf Q(s)
    \leq
    \frac1\delta\mathsf Q([0,M])
    <
    \infty.
\]
This proves \eqref{eq:MN-risk-limit}. In particular,
\begin{equation}
\label{eq:MN-continuous-extension-risk}
    R(\widetilde f_{\mathrm{MN}})
    =
    R(f_{\mathrm{MN}}).
\end{equation}

\medskip

To establish \eqref{eq:MN-optimality-characterization}, we first show
that $\mathsf Q$ has a strictly positive atom at zero. By
\eqref{eq:MN-empirical-spectral-gap} and
\eqref{eq:MN-residual-on-support},
\[
    \mathsf Q_n(\{0\})
    =
    \int_{[0,M]}
    \varphi_{\mathrm{MN}}(s)\,d\mathsf Q_n(s),
\]
for all sufficiently large $n$, almost surely. Since $\varphi_{\mathrm{MN}}$ is continuous on $[0,M]$
and agrees with $\mathbf 1_{\{0\}}$ on $E$ by
\eqref{eq:MN-residual-on-support}, weak convergence of
$\mathsf Q_n$, together with \eqref{eq:MN-limiting-support},
gives
\[
    \mathsf Q_n(\{0\})
    \longrightarrow
    \int_{[0,M]}
    \varphi_{\mathrm{MN}}(s)\,d\mathsf Q(s)
    =
    \mathsf Q(\{0\})
    \qquad\text{almost surely}.
\]
On the other hand, by \eqref{eq:definition_Q_n} and since $X_n$ has
full row rank, for all sufficiently large $n$, we get
\[
    \mathsf Q_n(\{0\})
    =
    \frac1p
    \sum_{\substack{1\leq i\leq p\\\lambda_i(S_n)=0}}
    \langle w_i,\Sigma_nw_i\rangle
    \geq
    c_\Sigma\frac{p-n}{p}.
\]
Passing to the limit yields
\[
    \mathsf Q(\{0\})
    \geq
    c_\Sigma\left(1-\frac1\gamma\right)>0.
\]
By \eqref{eq:def_mu_main},
\begin{equation}
\label{eq:MN-rho-at-zero}
    \mu(\{0\})
    =
    \mathsf C(\{0\})+\mathsf Q(\{0\})
    >0,
    \qquad
    \mathsf C(\{0\})
    =
    \int_{\{0\}}\rho(s)\,d\mu(s)
    =
    \mu(\{0\})\rho(0).
\end{equation}
Note that by \eqref{eq:def_B_n}, $B_n(g,h)=B_n(h,g)$ and $B_n(g,g)\geq0$, for $g,h\in C([0,M])$.
Passing to the limit gives the same properties for $B$.

Let $h\in C([0,M])$ and
$\varepsilon\in\mathbb R$.
Here $sh$ denotes the function $s\mapsto sh(s)$.
Expanding the two terms in
\eqref{eq:continuous-filter-risk-limit}, using symmetry of
$B$ and the identity
$s\widetilde f_{\mathrm{MN}}(s)=1$ on the positive support
of $\mathsf Q$, gives
\begin{align}
\label{eq:MN-risk-expansion}
\begin{split}
R(\widetilde f_{\mathrm{MN}}+\varepsilon h)
-
R(\widetilde f_{\mathrm{MN}})
&=
2\varepsilon
\left\{
    -\tau^2B(\varphi_{\mathrm{MN}},sh)
    +
    \sigma^2\gamma
    \int_{(0,M]}h(t)\,d\mathsf Q(t)
\right\}
\\
&\qquad+
\varepsilon^2
\left\{
    \tau^2B(sh,sh)
    +
    \sigma^2\gamma
    \int_{[0,M]}s h(s)^2\,d\mathsf Q(s)
\right\}.
\end{split}
\end{align}
We next evaluate $B(\varphi_{\mathrm{MN}},sh)$ using
\eqref{eq:GR-candidate-main}.
The integral representation of $B$ and the transformation
formula for the pushforward give
\begin{align*}
    B(\varphi_{\mathrm{MN}},sh)
    &=
    \gamma
    \int_{E\times E}
    \varphi_{\mathrm{MN}}(s)\,t h(t)\kappa(s,t)
    \,d\mu(s)\,d\mu(t)
    \\
    &\quad+
    \int_E
    \varphi_{\mathrm{MN}}(s)\,s h(s)
    \{\rho(s)-\gamma\ell(s)\}
    \,d\mu(s).
\end{align*}
The second integral vanishes: at $s=0$ because of the factor $s$, and
on $E\setminus\{0\}$ because $\varphi_{\mathrm{MN}}(s)=0$ by
\eqref{eq:MN-residual-on-support}.
Fubini's theorem is applicable to the remaining integral because
$\kappa\in L^1(\mu\otimes\mu)$ by
Lemma~\ref{lem:GR-radon-nikodym-energy} and the other
factors are bounded. Consequently,
\begin{align*}
    B(\varphi_{\mathrm{MN}},sh)
    &=
    \gamma
    \int_E\int_E
    \varphi_{\mathrm{MN}}(s)\,t h(t)\kappa(s,t)
    \,d\mu(t)\,d\mu(s)
    \\
    &=
    \gamma\mu(\{0\})
    \int_{(0,M]}
    t h(t)\kappa(0,t)\,d\mu(t).
\end{align*}
In the second equality, the integral in $s$ is restricted
to $\{0\}$ by \eqref{eq:MN-residual-on-support};
the contribution at $t=0$ is zero.
For $t>0$, \eqref{eq:kappa_density_def_main} gives
\[
    t\kappa(0,t)
    =
    t\frac{\rho(0)-\rho(t)}{-t}
    =
    \rho(t)-\rho(0).
\]
Using this identity, \eqref{eq:MN-rho-at-zero}, and
\eqref{eq:def_mu_main}, we obtain
\begin{align}\label{eq:B_rho_sh}
\begin{split}
    B(\varphi_{\mathrm{MN}},sh)
    &=
    \gamma
    \int_{(0,M]}
    h(t)
    \left\{
        \mu(\{0\})\rho(t)
        -
        \mathsf C(\{0\})
    \right\}
    \,d\mu(t)
    \\
    &=
    \gamma
    \left\{
        \bigl(
            \mathsf C(\{0\})
            +
            \mathsf Q(\{0\})
        \bigr)
        \int_{(0,M]}h(t)\,d\mathsf C(t)
        -
        \mathsf C(\{0\})
        \int_{(0,M]}h(t)\,d\mu(t)
    \right\}
    \\
    &=
    \gamma
    \left\{
        \mathsf Q(\{0\})
        \int_{(0,M]}h(t)\,d\mathsf C(t)
        -
        \mathsf C(\{0\})
        \int_{(0,M]}h(t)\,d\mathsf Q(t)
    \right\}.
\end{split}
\end{align}

Suppose now that
\[
    R(f_{\mathrm{MN}})
    =
    \inf_{f\in C([0,\infty))}R(f).
\]
By \eqref{eq:MN-continuous-extension-risk},
$\widetilde f_{\mathrm{MN}}$ is a minimizer in
$C([0,\infty))$.
Fix $h\in C([0,M])$ and extend it continuously to
$[0,\infty)$, for instance by setting $h(s)=h(M)$
for $s>M$.
Since $\widetilde f_{\mathrm{MN}}+\varepsilon h
\in C([0,\infty))$ for every $\varepsilon\in\mathbb R$,
the function
\[
    \varepsilon
    \mapsto
    R(\widetilde f_{\mathrm{MN}}+\varepsilon h)
\]
has a minimum at $\varepsilon=0$. Hence, its derivative at zero
vanishes. By \eqref{eq:MN-risk-expansion}, this means that the
coefficient of $2\varepsilon$ is zero. Substituting
\eqref{eq:B_rho_sh} into the coefficient of $2\varepsilon$ in \eqref{eq:MN-risk-expansion} gives
\begin{align*}
0
&=
-\tau^2\gamma
\left\{
\mathsf Q(\{0\})
\int_{(0,M]}h(t)\,d\mathsf C(t)
-
\mathsf C(\{0\})
\int_{(0,M]}h(t)\,d\mathsf Q(t)
\right\}
\\
&\qquad
+
\sigma^2\gamma
\int_{(0,M]}h(t)\,d\mathsf Q(t).
\end{align*}
Dividing by $\gamma$ and rearranging yields
\begin{align}\label{eq:ident_for_proof}
\tau^2\mathsf Q(\{0\})
\int_{(0,M]}h(t)\,d\mathsf C(t)
=
\left\{
\tau^2\mathsf C(\{0\})+\sigma^2
\right\}
\int_{(0,M]}h(t)\,d\mathsf Q(t).
\end{align}
Since $\mathsf Q(\{0\})>0$, define
\[
    a\coloneqq
    \frac{\mathsf C(\{0\})+\sigma^2/\tau^2}{\mathsf Q(\{0\})}
    \geq0.
\]
By \eqref{eq:ident_for_proof}, we have, for every
$h\in C([0,M])$,
\begin{align*}
    \int_{[0,M]}h(t)\,d(\mathsf C-a\mathsf Q)(t)
    =
    \{\mathsf C(\{0\})-a\mathsf Q(\{0\})\}h(0)
    =
    -\frac{\sigma^2}{\tau^2}h(0).
\end{align*}
Since finite Borel measures on $[0,M]$ are determined by their
integrals against continuous functions, we obtain
\begin{equation}
\label{eq:MN-measure-identity}
    \mathsf C
    =
    a\mathsf Q-\frac{\sigma^2}{\tau^2}\delta_0.
\end{equation}

By \eqref{eq:Q:stielt:def} and
\eqref{eq:MP-Silverstein-equation},
\begin{align}\label{eq:def_Q_limit}
    m_{\mathsf Q}(z)
    =
    -\frac1{\gamma z}\left(z+\frac1{v(z)}\right)
    =
    -\frac1z
    \int_{[0,\infty)}\frac{t}{1+tv(z)}\,dH(t),
    \qquad z\in\mathbb C^+.
\end{align}
Next, we construct the Stieltjes transform of the measure in
\eqref{eq:MN-measure-identity}, using $m_{\delta_0}(z)=-1/z$
and \eqref{eq:GR-C-transform}.
Multiplying by $-z$ therefore gives
\begin{equation}
\label{eq:MN-transform-identity}
    \int_{[0,\infty)}
    \frac{t}{1+tv(z)}\,d(G-aH)(t)
    =
    -\frac{\sigma^2}{\tau^2},
    \qquad z\in\mathbb C^+.
\end{equation}
The function
\[
    w\mapsto
    \int_{[0,\infty)}\frac{t}{1+tw}\,d(G-aH)(t)
\]
is holomorphic on $\mathbb C^+$. Indeed, the bounds
\[
    \left|\frac{t}{1+tw}\right|
    \leq\frac1{\operatorname{Im}(w)},
    \qquad
    \left|\frac{t^2}{(1+tw)^2}\right|
    \leq\frac1{\operatorname{Im}(w)^2},
    \qquad t\geq0,
\]
justify differentiation under the integral sign on compact
subsets of $\mathbb C^+$, since $G+aH$ is finite.
Moreover, $v$ is holomorphic and nonconstant by
\eqref{eq:MP-Silverstein-equation}. Thus
$\{v(z):z\in\mathbb C^+\}$ is a nonempty open subset of
$\mathbb C^+$ by the open mapping theorem.
Equation~\eqref{eq:MN-transform-identity} and the identity theorem
now imply
\begin{align}\label{eq:int_ident_1}
    \int_{[0,\infty)}
    \frac{t}{1+tw}\,d(G-aH)(t)
    =
    -\frac{\sigma^2}{\tau^2},
    \qquad w\in\mathbb C^+.
\end{align}
Taking $w=iy$, $y>0$, and using
$|t/(1+ity)|\leq1/y$, we obtain
\[
    \frac{\sigma^2}{\tau^2}
    =
    \left|
        \int_{[0,\infty)}
        \frac{t}{1+ity}\,d(G-aH)(t)
    \right|
    \leq
    \frac{1+a}{y}.
\]
Letting $y\to\infty$ proves that $\sigma^2=0$.
Hence, \eqref{eq:int_ident_1} reduces to
\begin{equation}
\label{eq:int_ident_2}
    \int_{[0,\infty)}
    \frac{t}{1+tw}\,d(G-aH)(t)
    =
    0,
    \qquad w\in\mathbb C^+.
\end{equation}
For $z\in\mathbb C^+$, substituting $w=-1/z\in\mathbb C^+$
into \eqref{eq:int_ident_2} and dividing by $-z$ then gives
\[
    \int_{[0,\infty)}
    \frac{t}{t-z}\,d(G-aH)(t)=0.
\]
By Assumption~\ref{ass:GR-basic-3} and the additional lower
eigenvalue bound of the theorem, $H_n$ and $G_n$ are supported on
\[
    [c_\Sigma,C_\Sigma],
    \qquad\text{where}\qquad
    C_\Sigma\coloneqq
    \sup_{n\geq1}\|\Sigma_n\|_2<\infty.
\]
Together with Assumptions~\ref{ass:GR-basic-4} and
\ref{ass:GR-basic-5}, this implies
\[
    \operatorname{supp}(H)
    \cup
    \operatorname{supp}(G)
    \subseteq
    [c_\Sigma,C_\Sigma].
\]
Hence, $A\mapsto \int_A t\,d(G-aH)(t)$ defines
a compactly supported finite signed measure whose Stieltjes transform
vanishes on $\mathbb C^+$. Lemma~\ref{lemma:injectivity-of-stieltjes-transform} yields
\[
    t\,dG(t)=a t\,dH(t).
\]
Integrating against $t^{-1}$ on arbitrary Borel subsets of
$[c_\Sigma,C_\Sigma]$ gives $G=aH$, since $c_\Sigma>0$.
As both measures have total mass one, $a=1$ and $G=H$.
This proves necessity.

Conversely, suppose that $\sigma^2=0$ and $G=H$.
By Corollary~\ref{cor:collapse-G-equals-H}, for every
$f\in C([0,\infty))$,
\[
    R(f)
    =
    \tau^2
    \int_{[0,M]}
    \{1-sf(s)\}^2\,d\mathsf Q(s).
\]
Since $\mathsf Q$ is supported on
$E=\{0\}\cup[\delta,M]$ and
$\widetilde f_{\mathrm{MN}}(s)=s^{-1}$ on $[\delta,M]$,
\[
    R(\widetilde f_{\mathrm{MN}})
    =
    \tau^2\mathsf Q(\{0\}).
\]
Moreover, for every $f\in C([0,\infty))$,
\[
    R(f)-R(\widetilde f_{\mathrm{MN}})
    =
    \tau^2
    \int_{(0,M]}
    \{1-sf(s)\}^2\,d\mathsf Q(s)
    \geq0.
\]
Together with
\eqref{eq:MN-continuous-extension-risk}, this proves the converse.
\end{proof}

\subsection{Proof of Theorem~\ref{cor:ridge-optimality}}

\begin{proof}[Proof of Theorem~\ref{cor:ridge-optimality}]
Fix $\lambda>0$ and set
\[
    \varphi_\lambda(s)
    \coloneqq
    1-sf_\lambda(s)
    =
    \frac{\lambda}{s+\lambda},
    \qquad s\in[0,M].
\]
As in \eqref{eq:MN-risk-expansion}, if $f_\lambda$ minimizes $R$,
then the coefficient of $\varepsilon$ in
$R(f_\lambda+\varepsilon h)-R(f_\lambda)$ must vanish. Hence,
\begin{equation}
\label{eq:ridge-first-variation}
    \tau^2B(\varphi_\lambda,sh)
    =
    \sigma^2\gamma
    \int_{[0,M]}
    \frac{s}{s+\lambda}h(s)\,d\mathsf Q(s),
    \qquad h\in C([0,M]).
\end{equation}

We evaluate the left-hand side using
\eqref{eq:GR-candidate-main}.
Theorem \ref{thm:GR-fix} and
and the transformation formula
for the pushforward give
\begin{align*}
    B(\varphi_\lambda,sh)
    &=
    \gamma
    \int_{[0,M]}\int_{[0,M]}
    \varphi_\lambda(s)\,t h(t)\kappa(s,t)
    \,d\mu(s)\,d\mu(t)
    \\
    &\quad+
    \int_{[0,M]}
    \varphi_\lambda(t)\,t h(t)
    \{\rho(t)-\gamma\ell(t)\}
    \,d\mu(t).
\end{align*}
By the definition of $\rho$ in \eqref{eq:def_mu_main},
\[
    \int_{[0,M]}
    \varphi_\lambda(t)\,t h(t)\rho(t)
    \,d\mu(t)
    =
    \int_{[0,M]}
    \varphi_\lambda(t)\,t h(t)
    \,d\mathsf C(t).
\]
Moreover, by the symmetry of $\kappa$ and the
definition of $\ell$ in \eqref{eq:kappa_density_def_main},
\[
    \ell(t)
    =
    \int_{[0,M]}
    \kappa(s,t)\,d\mu(s)
\]
for $\mu$-almost every $t$. Hence,
\begin{align*}
    B(\varphi_\lambda,sh)
    &=
    \int_{[0,M]}
    \varphi_\lambda(t)\,t h(t)\,d\mathsf C(t)
    \\
    &\quad+
    \gamma
    \int_{[0,M]}\int_{[0,M]}
    \{\varphi_\lambda(s)-\varphi_\lambda(t)\}
    t h(t)\kappa(s,t)
    \,d\mu(s)\,d\mu(t).
\end{align*}
Here Fubini's theorem applies since
$\kappa\in L^1(\mu\otimes\mu)$ and the other factors are bounded.
Since
\[
    \varphi_\lambda(s)-\varphi_\lambda(t)
    =
    -\frac{\lambda(s-t)}
    {(s+\lambda)(t+\lambda)},
\]
Equation~\eqref{eq:kappa_density_def_main} gives
\[
    \{\varphi_\lambda(s)-\varphi_\lambda(t)\}\kappa(s,t)
    =
    \frac{
        \lambda\{\rho(t)-\rho(s)\}
    }{
        (s+\lambda)(t+\lambda)
    }.
\]
This identity also holds for $s=t$, since both sides vanish.
Moreover, by \eqref{eq:def_mu_main},
\[
    \int_{[0,M]}
    \frac{1}{s+\lambda}\,d\mu(s)
    =
    \int_{[0,M]}
    \frac{1}{s+\lambda}\,d\mathsf C(s)
    +
    \int_{[0,M]}
    \frac{1}{s+\lambda}\,d\mathsf Q(s),
\]
whereas the definition of $\rho$ gives
\[
    \int_{[0,M]}
    \frac{\rho(s)}{s+\lambda}\,d\mu(s)
    =
    \int_{[0,M]}
    \frac{1}{s+\lambda}\,d\mathsf C(s).
\]
Consequently,
\begin{align}
\label{eq:ridge-bilinear-identity}
\begin{split}
    B(\varphi_\lambda,sh)
    &=
    \lambda
    \left\{
        1+
        \gamma
        \int_{[0,M]}
        \frac{1}{s+\lambda}\,d\mathsf Q(s)
    \right\}
    \int_{[0,M]}
    \frac{t}{t+\lambda}h(t)\,d\mathsf C(t)
    \\
    &\quad-
    \gamma\lambda
    \left\{
        \int_{[0,M]}
        \frac{1}{s+\lambda}\,d\mathsf C(s)
    \right\}
    \int_{[0,M]}
    \frac{t}{t+\lambda}h(t)\,d\mathsf Q(t).
\end{split}
\end{align}
Suppose that
$R(f_\lambda)=\inf_{f\in C([0,\infty))}R(f)$.
Substituting \eqref{eq:ridge-bilinear-identity} into
\eqref{eq:ridge-first-variation} yields
\begin{align}\label{eq:measure_cont_rep}
\begin{split}
&\tau^2\lambda
\left\{
    1+
    \gamma
    \int_{[0,M]}
    \frac{1}{s+\lambda}\,d\mathsf Q(s)
\right\}
\int_{[0,M]}
\frac{t}{t+\lambda}h(t)\,d\mathsf C(t)
\\
&\qquad=
\gamma
\left\{
    \tau^2\lambda
    \int_{[0,M]}
    \frac{1}{s+\lambda}\,d\mathsf C(s)
    +
    \sigma^2
\right\}
\int_{[0,M]}
\frac{t}{t+\lambda}h(t)\,d\mathsf Q(t)
\end{split}
\end{align}
for every $h\in C([0,M])$.
Define
\begin{align} \label{eq:def-a}
    a
    \coloneqq
    \frac{
        \gamma
        \left\{
            \tau^2\lambda
            \displaystyle\int_{[0,M]}
            \frac{1}{s+\lambda}\,d\mathsf C(s)
            +
            \sigma^2
        \right\}
    }{
        \tau^2\lambda
        \left\{
            1+
            \gamma
            \displaystyle\int_{[0,M]}
            \frac{1}{s+\lambda}\,d\mathsf Q(s)
        \right\}
    }
    >0.
\end{align}
Since finite Borel measures on $[0,M]$ are determined by their
integrals against continuous functions, \eqref{eq:measure_cont_rep} implies
\[
    \frac{t}{t+\lambda}\,d\mathsf C(t)
    =
    a
    \frac{t}{t+\lambda}\,d\mathsf Q(t).
\]
On each interval $[1/k,M]$,
$(t+\lambda)/t$ is bounded. Integrating against this term on Borel
subsets of $[1/k,M]$ and then letting $k\to\infty$ yields
\[
    \left.\mathsf C\right|_{(0,M]}
    =
    a\left.\mathsf Q\right|_{(0,M]}.
\]
Including the masses at zero gives
\begin{equation}
\label{eq:ridge-measure-identity}
    \mathsf C
    =
    a\mathsf Q
    +
    \{\mathsf C(\{0\})-a\mathsf Q(\{0\})\}\delta_0.
\end{equation}
Taking Stieltjes transforms and using
\eqref{eq:GR-C-transform_0} and \eqref{eq:def_Q_limit}, we obtain
\[
    \int_{[0,\infty)}
    \frac{t}{1+tv(z)}\,d(G-aH)(t)
    =
    \mathsf C(\{0\})-a\mathsf Q(\{0\}),
    \qquad z\in\mathbb C^+.
\]
The same holomorphic continuation argument used in the proof of
Theorem~\ref{thm:min-norm-optimality} gives
\[
    \int_{[0,\infty)}
    \frac{t}{1+tw}\,d(G-aH)(t)
    =
    \mathsf C(\{0\})-a\mathsf Q(\{0\}),
    \qquad w\in\mathbb C^+.
\]
Taking $w=iy$, $y>0$, and using
$|t/(1+ity)|\leq1/y$, we get
\[
    |\mathsf C(\{0\})-a\mathsf Q(\{0\})|
    \leq
    \frac{1+a}{y}.
\]
Letting $y\to\infty$ shows that
$\mathsf C(\{0\})=a\mathsf Q(\{0\})$, and hence, by
\eqref{eq:ridge-measure-identity},
$\mathsf C=a\mathsf Q$.
The preceding integral identity therefore reduces to
\[
    \int_{[0,\infty)}
    \frac{t}{1+tw}\,d(G-aH)(t)
    =
    0,
    \qquad w\in\mathbb C^+.
\]
Substituting $w=-1/z$, with $z\in\mathbb C^+$, gives
\[
    \int_{[0,\infty)}
    \frac{t}{t-z}\,d(G-aH)(t)
    =
    0,
    \qquad z\in\mathbb C^+.
\]
By Assumptions~\ref{ass:GR-basic-3}--\ref{ass:GR-basic-5},
both $G$ and $H$ are supported on $[0,C_\Sigma]$, where
\[
    C_\Sigma
    \coloneqq
    \sup_{n\geq1}\|\Sigma_n\|_2
    <\infty.
\]
Thus,  $A\mapsto \int_At\,d(G-aH)(t)$ induces a compactly supported  finite signed measure.
Lemma~\ref{lemma:injectivity-of-stieltjes-transform} yields
\[
    t\,dG(t)
    =
    a t\,dH(t).
\]
Integrating against $t^{-1}$ on Borel subsets of
$[1/k,C_\Sigma]$ and letting $k\to\infty$ gives
\[
    \left.G\right|_{(0,\infty)}
    =
    a\left.H\right|_{(0,\infty)}.
\]
In particular,
\begin{align}\label{eq:ratio_H_G}
    1-G(\{0\})
    =
    a\{1-H(\{0\})\}.
\end{align}
Since $H(\{0\})<1$, we have
\[
    a=\frac{1-G(\{0\})}{1-H(\{0\})}
    \qquad\text{and}\qquad
    \left.G\right|_{(0,\infty)}
    =
    \frac{1-G(\{0\})}{1-H(\{0\})}\left.H\right|_{(0,\infty)}.
\]
Since $\mathsf C=a\mathsf Q$, the definition of $a$ in \eqref{eq:def-a} gives
\[
    a\tau^2\lambda
    \left\{
        1+
        \gamma
        \int_{[0,M]}
        \frac{1}{s+\lambda}\,d\mathsf Q(s)
    \right\}
    =
    \gamma
    \left\{
        a\tau^2\lambda
        \int_{[0,M]}
        \frac{1}{s+\lambda}\,d\mathsf Q(s)
        +
        \sigma^2
    \right\}.
\]
The integral terms cancel, and using \eqref{eq:ratio_H_G} gives \[
\lambda=\frac{\sigma^2\gamma(1-H(\{0\}))}{\tau^2(1-G(\{0\}))}.
\]
This proves necessity.

Conversely, suppose that
\begin{align} \label{eq:lambda}
    \left.G\right|_{(0,\infty)}
    =
    \frac{1-G(\{0\})}{1-H(\{0\})}\left.H\right|_{(0,\infty)}
    \qquad\text{and}\qquad
    \lambda=\frac{\sigma^2\gamma(1-H(\{0\}))}{\tau^2(1-G(\{0\}))}>0.
\end{align}
Define $c\coloneqq  (1-G(\{0\}))/(1-H(\{0\}))$.
By \eqref{eq:GR-C-transform_0} and \eqref{eq:def_Q_limit},
\[
    m_{\mathsf C}(z)
    =
    -\frac1z
    \int_{[0,\infty)}
    \frac{t}{1+tv(z)}\,dG(t)
    =
    -\frac cz
    \int_{[0,\infty)}
    \frac{t}{1+tv(z)}\,dH(t)
    =
    c\,m_{\mathsf Q}(z),
    \qquad z\in\mathbb C^+,
\]
since the integrands vanish at $t=0$.
Lemma~\ref{lemma:injectivity-of-stieltjes-transform} therefore gives
$\mathsf C=c\mathsf Q$.
Consequently, $\mu=(1+c)\mathsf Q$ and
\[
    \rho(s)=\frac{c}{1+c}
\]
for $\mu$-almost every $s$.
By \eqref{eq:kappa_density_def_main}, it follows that
$\kappa=0$ $(\mu\otimes\mu)$-almost everywhere and
$\ell=0$ $\mu$-almost everywhere.
Thus \eqref{eq:GR-candidate-main} and the integral representation of $B$ in Theorem~\ref{thm:GR-fix} yield
\begin{align} \label{eq:rep_B}
    \mathsf K=cD_\#\mathsf Q,
    \qquad
    B(g,h)
    =
    c
    \int_{[0,M]}g(s)h(s)\,d\mathsf Q(s),
    \qquad
    g,h\in C([0,M];\mathbb C).
\end{align}
Consequently, using \eqref{eq:continuous-filter-risk-limit}, \eqref{eq:lambda} and \eqref{eq:rep_B}, completing the square gives
\[
    R(f)-R(f_\lambda)
    =
    \tau^2c
    \int_{[0,M]}
    s(s+\lambda)
    \left\{
        f(s)-\frac1{s+\lambda}
    \right\}^2
    \,d\mathsf Q(s)
    \geq0
\]
for every $f\in C([0,\infty))$.
Since $f_\lambda\in C([0,\infty))$, this proves
\eqref{eq:ridge-optimality-characterization}.
\end{proof}

\section*{Acknowledgements}
The work of TS and ND is partially supported by  the Aarhus University Research
Foundation (AUFF),  project numbers 47221 and 47388.
ND gratefully acknowledges the hospitality of the Department of
Statistics at Harvard University during the preparation of this
manuscript. ND's stay is supported by a RECRUIT Sabbatical
Grant from the Novo Nordisk Foundation, project number 52629.

\section*{AI use declaration}
Thomas Stark and Nina Dörnemann have made use of ChatGPT (GPT-6) in the course of this project. Use has included text editing, spell-checking, proof-reading, and proof editing.
 The examples in Section~\ref{sec:weak} were generated with AI assistance.
All other results and proofs are based on the authors' own ideas. AI assisted in the development of some of these arguments. The gaps in Theorem~B.9 of \cite{BaiSilverstein2010} and \cite{GreenRomanov2025} were identified by the authors themselves.
The authors reviewed the AI-assisted material and take full responsibility for the final manuscript.

\bibliographystyle{abbrvnat}
\bibliography{main}

\appendix
\clearpage
\pagenumbering{arabic}
\setcounter{page}{1}
\setcounter{section}{0}
\setcounter{subsection}{0}
\setcounter{subsubsection}{0}
\renewcommand{\thesection}{\Alph{section}}
\renewcommand{\theHsection}{supp.\Alph{section}}
\setcounter{equation}{0}
\setcounter{theorem}{0}
\setcounter{lemma}{0}
\setcounter{proposition}{0}
\setcounter{corollary}{0}
\setcounter{definition}{0}
\setcounter{example}{0}
\setcounter{remark}{0}
\setcounter{assumption}{0}
\setcounter{footnote}{0}
\numberwithin{equation}{section}
\numberwithin{figure}{section}
\numberwithin{table}{section}
\setcounter{figure}{0}
\setcounter{table}{0}
\markboth{T. Stark, G. Köstenberger and N. Dörnemann}{Supplementary Material}
\section*{Supplementary Material to "Beyond weak convergence of signed spectral measures: Eigenvector overlaps and high-dimensional prediction risk"}
\startcontents[supplement]

\printcontents[supplement]{}{1}{
}

\section{Proofs of results in  Section~\ref{sec:weak}}\label{appendix_A}
Before presenting the characterizations of weak and vague convergence in terms of Stieltjes transforms convergence, we would like to discuss another somewhat subtle point: When does a Stieltjes transform uniquely determine a measure? 
In the one-dimensional case, for a finite real measure $\mu$, it is enough to consider $m_\mu(z)$ for $z\in \mathbb{C}^+$, since this determines $m_\mu$ for all $z\in \mathbb{C}\setminus \mathbb{R}$, via $m_\mu(\overline{z}) = \overline{m_\mu(z)}$, and the measure can be reconstructed using the inversion formula (see, e.g., Theorem B.8 in \cite{BaiSilverstein2010}). 
In the multivariate case, this is no longer true. 
The values of $m_\mu(z_1,z_2)$ for $z_1,z_2\in \mathbb{C}^+$ only determine the values of $m_\mu(w_1,w_2)$ for $w_1,w_2\in \mathbb{C}^-$, but they do not determine the behavior of $m_\mu$ on $\mathbb{C}^+\times \mathbb{C}^-$ or $\mathbb{C}^-\times \mathbb{C}^+$. 
There exist probability measures $P$ and $Q$, such that $m_P(z_1,z_2) = m_Q(z_1,z_2)$ for all $z_1,z_2\in \mathbb{C}^+$, but $P\neq Q$ (see Example \ref{ex:noninjectivity-upper-poly-halfplane} below).
However, if the measures in question are compactly supported, the maximal domain of the Stieltjes transforms is open and connected, and hence, they are uniquely determined by their values on $(\mathbb{C}^+)^d$ by the identity theorem (see Lemma \ref{lemma:injectivity-of-stieltjes-transform} below). 
Proposition \ref{prop:multivariate-stieltjes-inversion} below constitutes an inversion formula, similar to the inversion formulas for Cauchy-transforms considered in \cite{HuangWang2016, Nedic2023}. 
An immediate consequence of the inversion formula is the injectivity of the Stieltjes transform, which we record in Lemma \ref{lemma:injectivity-of-stieltjes-transform}.

\subsection{Proposition~\ref{prop:multivariate-stieltjes-inversion} and its proof}

\begin{proposition}
\label{prop:multivariate-stieltjes-inversion}
Let $\mu$ be a finite complex Borel measure on $\mathbb R^d$. Then,
for every $f\in C_0(\mathbb R^d;\mathbb C)$,
\[
\int_{\mathbb R^d} f(t)\,d\mu(t)
=
\lim_{\eta\downarrow0}
\frac{1}{(2\pi i)^d}
\int_{\mathbb R^d}
f(x)
\sum_{\varepsilon\in\{-1,1\}^d}
\left(\prod_{j=1}^d\varepsilon_j\right)
m_\mu
\left(
x_1+i\varepsilon_1\eta,\ldots,
x_d+i\varepsilon_d\eta
\right)
\,dx.
\]
\end{proposition}

\begin{proof}[Proof of Proposition~\ref{prop:multivariate-stieltjes-inversion}]
For $\eta>0$, define the Poisson kernel
\[
    P_\eta(u)
    \coloneqq
    \frac{1}{\pi}\frac{\eta}{u^2+\eta^2},
    \qquad
    u\in\mathbb R.
\]
For every $t,x\in\mathbb R$,
\begin{align}
\label{eq:poisson-resolvent}
    \frac{1}{2\pi i}
    \left\{
        \frac{1}{t-(x+i\eta)}
        -
        \frac{1}{t-(x-i\eta)}
    \right\}
    =
    P_\eta(t-x).
\end{align}
Applying \eqref{eq:poisson-resolvent} in each coordinate and using
the definition of the multivariate Stieltjes transform gives
\begin{align}
\label{eq:multiv-poisson}
\frac{1}{(2\pi i)^d}
\sum_{\varepsilon\in\{-1,1\}^d}
\left(\prod_{j=1}^d\varepsilon_j\right)
m_\mu
\left(
x_1+i\varepsilon_1\eta,\ldots,
x_d+i\varepsilon_d\eta
\right)
=
\int_{\mathbb R^d}
\prod_{j=1}^d
P_\eta(t_j-x_j)
\,d\mu(t).
\end{align}

Define
\[
    Q_\eta(u)
    \coloneqq
    \prod_{j=1}^dP_\eta(u_j),
    \qquad
    u\in\mathbb R^d.
\]
We verify that $\{Q_\eta\}_{\eta>0}$ is an approximate identity on
$\mathbb R^d$ in the sense of
\cite[Section~7.6.1]{MakarovPodkorytov2013}. Since $P_\eta\geq0$, we
have
\[
    Q_\eta(u)\geq0,
    \qquad
    u\in\mathbb R^d.
\]
Moreover, by Tonelli's theorem,
\[
    \int_{\mathbb R^d}Q_\eta(u)\,du
    =
    \prod_{j=1}^d
    \int_{\mathbb R}P_\eta(u_j)\,du_j
    =
    1.
\]
Finally,
\[
    Q_\eta(u)
    =
    \eta^{-d}Q_1(u/\eta).
\]
Hence, for every $\delta>0$, the change of variables $y=u/\eta$
gives
\[
    \int_{\{\|u\|_2>\delta\}}
    Q_\eta(u)\,du
    =
    \int_{\{\|y\|_2>\delta/\eta\}}
    Q_1(y)\,dy
    \longrightarrow0
    \qquad\text{as }\eta\downarrow0,
\]
since $Q_1\in L^1(\mathbb R^d)$. Thus
$\{Q_\eta\}_{\eta>0}$ is an approximate identity on
$\mathbb R^d$; cf.~\cite[Section~7.6]{MakarovPodkorytov2013}.

Now let $f\in C_0(\mathbb R^d;\mathbb C)$. Since
\[
    \int_{\mathbb R^d}
    \int_{\mathbb R^d}
    |f(x)|Q_\eta(t-x)
    \,dx\,d|\mu|(t)
    \leq
    \|f\|_\infty\|\mu\|_{\mathrm{TV}}
    <\infty,
\]
Fubini's theorem is applicable. Multiplying
\eqref{eq:multiv-poisson} by $f(x)$ and integrating with respect to
$x$ therefore yields
\begin{align*}
\frac{1}{(2\pi i)^d}
\int_{\mathbb R^d}
f(x)
\sum_{\varepsilon\in\{-1,1\}^d}
\left(\prod_{j=1}^d\varepsilon_j\right)
m_\mu
\left(
x_1+i\varepsilon_1\eta,\ldots,
x_d+i\varepsilon_d\eta
\right)
\,dx
=
\int_{\mathbb R^d}
(Q_\eta*f)(t)
\,d\mu(t),
\end{align*}
where
$(Q_\eta*f)(t)
\coloneqq
\int_{\mathbb R^d}f(x)Q_\eta(t-x)\,dx$
denotes the convolution of $Q_\eta$ and $f$.
Since $f\in C_0(\mathbb R^d;\mathbb C)$ is bounded and uniformly
continuous, \cite[Theorem~7.6.3 and the subsequent remark]
{MakarovPodkorytov2013} yields
\[
    \|Q_\eta*f-f\|_\infty
    \longrightarrow0
    \qquad\text{as }\eta\downarrow0.
\]
Therefore,
\[
    \left|
        \int_{\mathbb R^d}
        (Q_\eta*f-f)(t)
        \,d\mu(t)
    \right|
    \leq
    \|\mu\|_{\mathrm{TV}}
    \|Q_\eta*f-f\|_\infty
    \longrightarrow0,
\]
which proves the claim.
\end{proof}

\subsection{Proof of Theorem \ref{thm:weak-convergence}}

\begin{proof}[Proof of Theorem \ref{thm:weak-convergence}] 
Suppose first $\mu_n\Rightarrow \mu$ for some finite signed Borel measure $\mu$ on $\mathbb R^d$. 
For fixed $z\in\mathbb C\setminus\mathbb R$, write $z=a+ib$, where
$a\in\mathbb R$ and $b\neq0$. Then
\[
    \frac{1}{x-z}
    =
    \frac{x-a}{(x-a)^2+b^2}
    +
    i\,\frac{b}{(x-a)^2+b^2}.
\]
Hence, the real and imaginary parts of the function
$x\mapsto(x-z)^{-1}$ are bounded and continuous on $\mathbb R$.
Similarly, for $z_1,\ldots,z_d\in\mathbb C\setminus\mathbb R$, the real and imaginary
parts of the function
\[
    (x_1,\ldots,x_d)
    \mapsto
    \frac{1}{(x_1-z_1)\cdots(x_d-z_d)}
\]
are bounded continuous functions on $\mathbb R^d$, since they are finite sums and products of bounded continuous functions. Applying weak convergence to the real and imaginary parts separately therefore yields
\[
    m_{\mu_n}(z_1,\ldots,z_d)
    \longrightarrow
    m_\mu(z_1,\ldots,z_d)
\]
for every $z_1,\ldots,z_d\in\mathbb C\setminus \mathbb R$.
This proves condition~$2$(a).

Conditions $2$(b) and $2$(c) are a consequence of Prokhorov's theorem \cite[Theorem 8.6.2]{bogachev-measure-theory}, since every sequence in $\{\mu_n\}_{n\geq 1}$ either contains a constant subsequence, or contains a subsequence of $(\mu_n)_{n\geq 1}$, it contains a weakly convergent subsequence.

Conversely, suppose that Conditions~$2$(a)-(c) are satisfied. 
By \cite[Theorem 8.6.2]{bogachev-measure-theory} every subsequence $\mu_{n_k}$ of $\mu_n$ has a weakly convergence subsequence $\mu_{n_{k_l}}$.
We denote this weak limit by $\mu$, and note that
\begin{equation*}
    m_\mu(z_1,\dots, z_d) = \lim_{l\to \infty} m_{\mu_{n_{k_l}}}(z_1,\dots, z_d)
\end{equation*}
for all $z_1,\dots, z_d\in \mathbb{C}\setminus \mathbb{R}$, by the arguments used in the first half of this proof. 
If $\mu_{n_{m}}$ is a different subsequence of $\mu_n$ containing a different weakly convergent subsequence $\mu_{n_{m_{r}}}\Rightarrow \nu$, we get 
\begin{equation*}
    m_\mu(z_1,\dots, z_d) = \lim_{l\to \infty} m_{\mu_{n_{k_l}}}(z_1,\dots, z_d)= \lim_{r\to \infty} m_{\mu_{n_{m_{r}}}}(z_1,\dots, z_d)  = m_\nu(z_1,\dots, z_d)
\end{equation*}
for all $z_1,\dots, z_d\in \mathbb{C}\setminus\mathbb{R}$. 
Lemma \ref{lemma:injectivity-of-stieltjes-transform} implies that $\nu = \mu$, and hence, every subsequence of $\mu_n$ contains a subsequence which converges to $\mu$.
This implies that $\mu_n\Rightarrow \mu$, completing the proof.

\end{proof}

\subsection{Proof of Proposition~\ref{cor:weak-convergence-compact-support}}

\begin{proof}[Proof of Proposition~\ref{cor:weak-convergence-compact-support}]
The implication $1\Rightarrow2$ follows directly from
Theorem~\ref{thm:weak-convergence}. Conversely, the eventual common
compact support implies tightness of $\{|\mu_n|:n\geq1\}$, and together
with the uniform total-variation bound yields weakly convergent
subsequences; pointwise convergence of the Stieltjes transforms on
$(\mathbb C^+)^d$ and
Lemma~\ref{lemma:injectivity-of-stieltjes-transform} show that all
subsequential limits coincide, and hence $\mu_n\Rightarrow\mu$, with
$\operatorname{supp}(\mu)\subseteq K$.
\end{proof}

\subsection{Lemma~\ref{lemma:injectivity-of-stieltjes-transform} and its proof}

\begin{lemma}\label{lemma:injectivity-of-stieltjes-transform}
    Let $\mu$ and $\nu$ be two finite complex Borel measures on $\mathbb{R}^d$. Then $\mu = \nu$ if and only if
    \begin{equation*}
        m_{\mu}(z_1,\dots, z_d) = m_\nu(z_1,\dots, z_d),\qquad \text{for every } z_1,\dots, z_d \in \mathbb C\setminus \mathbb R.
    \end{equation*} 
    Moreover, if $\mu$ and $\nu$ are compactly supported, then $\mu=\nu$ if and only if
\[
m_\mu(z_1,\ldots,z_d)=m_\nu(z_1,\ldots,z_d),\qquad\text{for every } z_1,\dots, z_d \in \mathbb C^+.
\]
\end{lemma}
\begin{proof}[Proof of Lemma \ref{lemma:injectivity-of-stieltjes-transform}]
The implication from $\mu=\nu$ to equality of the Stieltjes transforms
is immediate.

Conversely, suppose that
\[
    m_\mu(z_1,\ldots,z_d)
    =
    m_\nu(z_1,\ldots,z_d),
    \qquad
    \text{for every }
    z_1,\ldots,z_d\in\mathbb C\setminus\mathbb R.
\]
Set $\rho\coloneqq\mu-\nu.$
Then $\rho$ is a finite complex Borel measure on $\mathbb R^d$ and
\[
    m_\rho(z_1,\ldots,z_d)=0,
    \qquad
    \text{for every }
    z_1,\ldots,z_d\in\mathbb C\setminus\mathbb R.
\]
Let $f\in C_0(\mathbb R^d;\mathbb C)$. By
Proposition~\ref{prop:multivariate-stieltjes-inversion},
\begin{align*}
    \int_{\mathbb R^d}f(t)\,d\rho(t)
    =
    \lim_{\eta\downarrow0}
    \frac{1}{(2\pi i)^d}
    \int_{\mathbb R^d}
    f(x)
    \sum_{\varepsilon\in\{-1,1\}^d}
    \left(\prod_{j=1}^d\varepsilon_j\right)
    m_\rho
    \left(
        x_1+i\varepsilon_1\eta,\ldots,
        x_d+i\varepsilon_d\eta
    \right)
    \,dx.
\end{align*}
For every $\eta>0$, $x\in\mathbb R^d$, and
$\varepsilon\in\{-1,1\}^d$, we have $x_j+i\varepsilon_j\eta
    \in\mathbb C\setminus\mathbb R$,
for $j\in \{1,\ldots,d\}$.
Hence,
\[
    m_\rho
    \left(
        x_1+i\varepsilon_1\eta,\ldots,
        x_d+i\varepsilon_d\eta
    \right)
    =0.
\]
It follows that
\[
    \int_{\mathbb R^d}f\,d\rho=0
\]
for every $f\in C_0(\mathbb R^d;\mathbb C)$, and hence $\rho = 0$ by Riesz's representation theorem
\cite[Theorem~6.19]{Rudin1987}. 

It remains to prove the second assertion. Suppose that $\mu$ and $\nu$
are compactly supported and that
\[
    m_\mu(z_1,\ldots,z_d)
    =
    m_\nu(z_1,\ldots,z_d),
    \qquad
    z_1,\ldots,z_d\in\mathbb C^+.
\]
Set $\rho\coloneqq\mu-\nu$. Since $\mu$ and $\nu$ are compactly
supported, there exists $M>0$ such that
\[
    \operatorname{supp}(\rho)\subseteq[-M,M]^d.
\]
Let $\Omega_M\coloneqq\mathbb C\setminus[-M,M]$.
Then $\Omega_M$ is open and connected, and $m_\rho$ is holomorphic in
each coordinate on $\Omega_M^d$. 
If $d=1$, the identity theorem \cite[Corollary to Theorem~10.18]{Rudin1987} immediately gives
\[
    m_\rho(z)=0,
    \qquad z\in\Omega_M.
\]
If $d\geq2$, fix $z_2,\ldots,z_d\in\mathbb C^+$. The function
\[
    z_1\mapsto m_\rho(z_1,z_2,\ldots,z_d)
\]
is holomorphic on $\Omega_M$ and vanishes on $\mathbb C^+$. Hence, by
the identity theorem, it vanishes identically
on $\Omega_M$. Therefore,
\[
    m_\rho(z_1,\ldots,z_d)=0
\]
for every $z_1\in\Omega_M$ and $z_2,\ldots,z_d\in\mathbb C^+$.
Now fix $z_1\in\Omega_M$ and $z_3,\ldots,z_d\in\mathbb C^+$ and apply
the same argument to
\[
    z_2\mapsto m_\rho(z_1,z_2,z_3,\ldots,z_d).
\]
Repeating this argument successively in each of the remaining
coordinates, yields
\[
    m_\rho(z_1,\ldots,z_d)=0,
    \qquad
    \text{for every }z_1,\ldots,z_d\in\Omega_M.
\]
Since $\mathbb C\setminus\mathbb R\subseteq\Omega_M$, the first
assertion of the lemma applies and yields $\rho=0$, and hence, $\mu=\nu$.
\end{proof}

\subsection{Proof of Theorem~\ref{thm:vague-convergence}}
\begin{proof}[Proof of Theorem~\ref{thm:vague-convergence}]
Suppose first that $\mu_n$ converges vaguely to a finite signed Borel
measure $\mu$. Fix $z_1,\ldots,z_d\in\mathbb C\setminus\mathbb R$ and define
\[
    h(x_1,\ldots,x_d)
    \coloneqq
    \prod_{j=1}^d\frac{1}{x_j-z_j},
    \qquad
    (x_1,\ldots,x_d)\in\mathbb R^d.
\]
We first note that $h\in C_0(\mathbb R^d;\mathbb C)$. Indeed, if
$\|x\|_2\to\infty$, then
$\max_{1\leq j\leq d}|x_j|\to\infty$. Choosing $j$ such that
$|x_j|=\max_k|x_k|$, we have
\[
    \frac{1}{|x_j-z_j|}\longrightarrow0,
\]
whereas, for every $k$,
\[
    \frac{1}{|x_k-z_k|}
    \leq \frac{1}{|\operatorname{Im}(z_k)|}.
\]
Hence $h(x)\to0$ as $\|x\|_2\to\infty$, and therefore
$h\in C_0(\mathbb R^d;\mathbb C)$.

Since $C_c(\mathbb R^d;\mathbb C)$ is dense in
$C_0(\mathbb R^d;\mathbb C)$ with respect to the supremum norm, for
every $\varepsilon>0$ there exists
$h_\varepsilon\in C_c(\mathbb R^d;\mathbb C)$ such that $\|h-h_\varepsilon\|_\infty<\varepsilon$.
Therefore,
\begin{align*}
\left|
    \int_{\mathbb R^d}h\,d\mu_n
    -
    \int_{\mathbb R^d}h\,d\mu
\right|
&\leq
\left|
    \int_{\mathbb R^d}(h-h_\varepsilon)\,d\mu_n
\right|
+
\left|
    \int_{\mathbb R^d}h_\varepsilon\,d\mu_n
    -
    \int_{\mathbb R^d}h_\varepsilon\,d\mu
\right|
+
\left|
    \int_{\mathbb R^d}(h_\varepsilon-h)\,d\mu
\right|
\\
&\leq
\varepsilon
\left\{
    \|\mu_n\|_{\mathrm{TV}}
    +
    \|\mu\|_{\mathrm{TV}}
\right\}
+
\left|
    \int_{\mathbb R^d}h_\varepsilon\,d\mu_n
    -
    \int_{\mathbb R^d}h_\varepsilon\,d\mu
\right|.
\end{align*}
By vague convergence, the last term converges to zero as
$n\to\infty$, and hence,
\[
\limsup_{n\to\infty}
\left|
    \int_{\mathbb R^d}h\,d\mu_n
    -
    \int_{\mathbb R^d}h\,d\mu
\right|
\leq
\varepsilon
\left\{
    \sup_{n\geq1}\|\mu_n\|_{\mathrm{TV}}
    +
    \|\mu\|_{\mathrm{TV}}
\right\}.
\]
Since $\varepsilon>0$ was arbitrary,
\[
    \int_{\mathbb R^d}h\,d\mu_n
    \longrightarrow
    \int_{\mathbb R^d}h\,d\mu.
\]
By the definition of $h$, this is precisely
\[
    m_{\mu_n}(z_1,\ldots,z_d)
    \longrightarrow
    m_\mu(z_1,\ldots,z_d).
\]
Thus Condition~$2$ holds.

The same approximation argument shows more generally that
\begin{equation*}\label{eq:vague-C0-convergence}
    \int_{\mathbb R^d}f\,d\mu_n
    \longrightarrow
    \int_{\mathbb R^d}f\,d\mu
\end{equation*}
for every $f\in C_0(\mathbb R^d;\mathbb C)$.

Conversely, suppose that $m_{\mu_n}(z_1,\ldots,z_d)$ converges for every
$z_1,\ldots,z_d\in\mathbb C\setminus\mathbb R$. 
For every $n\geq1$, define the bounded linear functional $T_n:C_0(\mathbb R^d;\mathbb C)\longrightarrow\mathbb C$ by
\[
    T_n(f)
    \coloneqq
    \int_{\mathbb R^d}f\,d\mu_n.
\]
By the Riesz representation theorem, we have $\|T_n\|=\|\mu_n\|_{\mathrm{TV}}$,
and hence
\[
    \sup_{n\geq1}\|T_n\|\leq C,
\]
for some constant $C>0$.
The ball
\[
    B_C
    \coloneqq
    \left\{
        T\in C_0(\mathbb R^d;\mathbb C)^*:
        \|T\|\leq C
    \right\}.
\]
is weak* compact by the Banach--Alaoglu theorem, and since $C_0(\mathbb{R}^d; \mathbb{C})$ is separable, it is weak* sequentially compact.
Thus every subsequence $\mu_{n_k}$ of $\mu_n$ has a weak*-convergent subsequence $\mu_{n_{k_l}}$.
We denote the limit of this subsequence by $\mu$, and note that
\begin{equation*}
    m_\mu(z_1,\dots, z_d) = \lim_{l\to \infty} m_{\mu_{n_{k_l}}}(z_1,\dots, z_d),
\end{equation*}
for every $z_1,\dots, z_d\in \mathbb{C}\setminus \mathbb{R}$. 
By the same argument as in the proof of Theorem \ref{thm:weak-convergence}, we show that this limit is unique, and hence, every subsequence of $\mu_n$ contains a subsequence converging to $\mu$ with respect to the weak* topology, i.e., $\mu_n$ converges to $\mu$ vaguely.

\end{proof}

\subsection{Example~\ref{ex:noninjectivity-upper-poly-halfplane}}

\begin{example}\label{ex:noninjectivity-upper-poly-halfplane}
In this example we construct two probability measures $P$ and $Q$ on $\mathbb{R}^2$ with $m_P(z_1,z_2) = m_Q(z_1,z_2)$ for all $z_1,z_2\in \mathbb{C}^+$, and $P\neq Q$.
This shows that the behavior of the Stieltjes transform on $\mathbb{C}^+\times \mathbb{C}^+$ does not uniquely determine the underlying probability measures.
We consider probability measures $P$ and $Q$ with densities
\begin{equation*}
    \begin{split}
        p(x,y) &= \frac{2}{\pi^2} \frac{1}{1+ x^2} \frac{1}{1+y^2}\cos^2(\arctan(x) - \arctan(y)),\\
        q(x,y) &= \frac{2}{\pi^2} \frac{1}{1+ x^2} \frac{1}{1+y^2}\sin^2(\arctan(x) - \arctan(y))
    \end{split}
\end{equation*}
respectively. We first show that these are indeed densities. Clearly, $p,q\geq 0$. We show that 
\begin{equation*}
    I = \int_{\mathbb{R}^2} p(x,y) dx dy = 1. 
\end{equation*}
The function $q$ also integrates to $1$ via the same argument.
Substituting
\begin{equation*}
    u = \arctan(x), \qquad \text{and} \qquad v = \arctan(y)
\end{equation*}
we get that 
\begin{equation*}
    I = \frac{2}{\pi^2}\int_{-\frac{\pi}{2}}^{\frac{\pi}{2}}\int_{-\frac{\pi}{2}}^{\frac{\pi}{2}} \cos^2(u-v) du dv. 
\end{equation*}
Using the identity
\begin{equation*}
    \cos^2(x) = \frac{1 + \cos(2x)}{2},
\end{equation*}
we get that 
\begin{equation*}
    I = \frac{2}{\pi^2} \frac{\pi^2}{2} + \frac{1}{\pi^2}\int_{-\frac{\pi}{2}}^{\frac{\pi}{2}}\int_{-\frac{\pi}{2}}^{\frac{\pi}{2}} \cos(2u - 2v) du dv = 1. 
\end{equation*}

Next, we would like to show that
\begin{equation}\label{eq:ex:noninjectivity-upper-poly-halfplane:difference}
    p(x,y) - q(x,y) = \frac{1}{\pi^2} \biggl( \frac{1}{(x-i)^2 (y+i)^2} + \frac{1}{(x+i)^2(y-i)^2}\biggr). 
\end{equation}
First, recall that 
\begin{equation*}
    \cos^2(x) - \sin^2(x) = \cos(2x). 
\end{equation*}
Using this identity together with $\cos(x) = \frac{1}{2} \bigl(e^{ix} + e^{-ix} \bigr)$, we get that
\begin{equation*}
    p(x,y) - q(x,y) = \frac{2}{\pi^2} \frac{1}{1+x^2}\frac{1}{1+y^2} \frac{1}{2}\biggl( e^{2i\arctan(x)}e^{-2i\arctan(y)} + e^{-2i\arctan(x)}e^{2i\arctan(y)} \biggr)
\end{equation*}
The identity in Equation \eqref{eq:ex:noninjectivity-upper-poly-halfplane:difference} is now a simple consequence of the two identities
\begin{equation*}
    \exp(-2i\arctan(x)) = - \frac{1+x^2}{(x-i)^2} \quad \text{and} \quad\exp(2i\arctan(x)) = - \frac{1+x^2}{(x+i)^2}. 
\end{equation*}
In the final step, we are going to show that 
\begin{equation*}
    m_{P-Q}(z_1, z_2) = \int_{\mathbb{R}^2} \frac{p(x,y) - q(x,y)}{(x-z_1) (y-z_2)} dx dy = 0 
\end{equation*}
for all $z_1,z_2 \in \mathbb{C}^+$. 
Using Identity \eqref{eq:ex:noninjectivity-upper-poly-halfplane:difference}, we can write $m_{P-Q} =\frac{1}{\pi^2}( A_1 + A_2)$, where
\begin{equation*}
    \begin{split}
        A_1(z_1,z_2) & = \int_{\mathbb{R}^2} \frac{1}{(x-i)^2(y+i)^2(x-z_1)(y-z_2)} dx dy\\
        A_2(z_1,z_2) & = \int_{\mathbb{R}^2} \frac{1}{(x+i)^2(y-i)^2(x-z_1)(y-z_2)} dx dy.
    \end{split}
\end{equation*}
We show that $A_1(z_1,z_2) = 0$ for all $z_1,z_2\in \mathbb{C}^+$.
The second term $A_2$ is zero by the same argument.
Using Fubini's theorem we write $A_1(z_1,z_2) = B_1(z_1)C_1(z_2)$, where
\begin{equation*}
\begin{split}
    B_1(z_1) &= \int_\mathbb{R} \frac{1}{(x-i)^2(x-z_1)} dx, \\
    C_1(z_2) &= \int_\mathbb{R} \frac{1}{(y+i)^2(y-z_2)} dy. \\
\end{split}
\end{equation*}
We show that $B_1(z_1) = 0$ for all $z_1\in \mathbb{C}^+$. 
Fix $z\in \mathbb{C}^+$ and  define
\begin{equation*}
    f(x) = \frac{1}{(x-i)^2(x-z)}. 
\end{equation*}
Denote by $\Sigma_M$ a semicircle in the upper complex half-plane connecting $-M$ to $M$, and denote by $\Gamma_M$ the closed contour $[-M,M] \cup \Sigma_M$. 
Note that 
\begin{equation*}
    \int_{\Sigma_M} f(w) dw = iM\int_0^\pi e^{it}f(Me^{it})\,dt 
\end{equation*}
for $M\to \infty$, since $|f(w)| = O(|w|^{-3})$. 
This implies that 
\begin{equation*}
    \int_\mathbb{R} f(x) dx = \lim_{M\to \infty} \int_{\Gamma_M} f(w) dw. 
\end{equation*}
We fix an $M$ with 
$M>\max\{1,|z|\}$, and distinguish two cases: $z=i$ and $z\neq i$. 

First, if $z \neq i$,  the Residue Theorem implies 
\begin{equation*}
    \frac{1}{2\pi i}\int_{\Gamma_M} f(w) dw = \mathrm{Res}_f(i) + \mathrm{Res}_f(z). 
\end{equation*}
Since $z$ is a simple pole, we have
\begin{equation*}
    \mathrm{Res}_f(z)  = \frac{1}{(z-i)^2},
\end{equation*}
and since $i$ is a double pole, we get that
\begin{equation*}
    \mathrm{Res}_f(i) = \frac{d}{dw} \frac{1}{w-z} \Big|_{w = i} = - \frac{1}{(z-i)^2}. 
\end{equation*}
Hence, $\mathrm{Res}_f(i) + \mathrm{Res}_f(z) = 0$, which implies $B_1(z) = 0$ in this case. 
If $z = i$, then $f(w)=(w-i)^{-3}$ has a pole of order three with zero residue. Hence, 
\begin{equation*}
    \frac{1}{2\pi i}\int_{\Gamma_M} f(w) dw = \mathrm{Res}_f(i) = 0.
\end{equation*}
In total, we get that $B_1(z) = 0$ for all $z\in \mathbb{C}^+$, and hence, $m_{P-Q}(z_1,z_2) = 0$ for all $z_1,z_2\in \mathbb{C}^+$, but $P\neq Q$. 

\end{example}

\section{Technical details for Sections~\ref{sec:application}
and~\ref{sec:proof-counterexample}}\label{appendix_B}

\subsection{Proof of Proposition~\ref{prop:risk-decomposition}}
\begin{proof}[Proof of Proposition~\ref{prop:risk-decomposition}]
Since $y=X\beta_n^\star+\varepsilon$, we have from \eqref{eq:spectral-estimator}
\[ \hat\beta(f)=f(S_n)S_n\beta_n^\star + f(S_n)\frac{X^\top\varepsilon}{n}.
\] 
Because $S_n$ and $f(S_n)$ commute, 
\begin{align} \label{eq:formula_diff_hat_beta_beta_star}
\hat\beta(f)-\beta_n^\star = -\{I_p-S_nf(S_n)\}\beta_n^\star + f(S_n)\frac{X^\top\varepsilon}{n}. 
\end{align}
Moreover, $\mathbb E[\varepsilon\mid X]=0$, $\mathbb E[\varepsilon\varepsilon^\top\mid X] = \sigma^2 I_n$. 
Consequently,
\begin{align*} 
\mathrm R_{\Sigma_n}\bigl(\hat\beta(f)\bigr) &= (\beta_n^\star)^\top \{I_p-S_nf(S_n)\} \Sigma_n \{I_p-S_nf(S_n)\} \beta_n^\star  + \frac{\sigma^2}{n} \tr(f(S_n)\Sigma_nf(S_n)S_n). 
\end{align*}
For notational convenience we write $s_i\coloneqq \lambda_i(S_n)$. By the spectral decomposition of $S_n$ in \eqref{eq:sam-cov-decomposition} and the definition of $\widehat{\mathsf K}^{\beta_n^\star}_n$ in \eqref{eq:def_K}, we get
\begin{align*} 
&(\beta_n^\star)^\top \{I_p-S_nf(S_n)\} \Sigma_n \{I_p-S_nf(S_n)\} \beta_n^\star \\ 
&\quad= \sum_{i,j=1}^p \bigl(1-s_if(s_i)\bigr) \bigl(1-s_jf(s_j)\bigr) \langle\beta_n^\star,w_i\rangle \langle w_i,\Sigma_nw_j\rangle \langle w_j,\beta_n^\star\rangle \\ 
&\quad= \|\beta_n^\star\|_2^2 \int_{[0,\infty)^2}  \bigl(1-sf(s)\bigr) \bigl(1-tf(t)\bigr) \,d\widehat{\mathsf K}^{\beta_n^\star}_n(s,t). \end{align*} 
Similarly, 
\begin{align*} 
\frac{\sigma^2}{n} \tr(f(S_n)\Sigma_nf(S_n)S_n)&=\frac{\sigma^2}{n} \tr\{\Sigma_nS_nf(S_n)^2\} = \frac{\sigma^2}{n} \sum_{i=1}^p s_if(s_i)^2w_i^\top\Sigma_nw_i \\ &= \sigma^2\frac pn \int_{[0,\infty)} sf(s)^2\,d\mathsf Q_n(s) = \sigma^2\gamma_n \int_{[0,\infty)} sf(s)^2\,d\mathsf Q_n(s). \end{align*} 
Combining the two identities proves the claim.
\end{proof}

Next, we introduce some notation used throughout the rest of this section. Following \cite{LinPan2026}, let
$\mathfrak m$ denote the Stieltjes transform of the deterministic
self-consistent measure $\varrho$ associated with the population
covariance matrix $\Sigma_n$. Throughout this section, we suppress the
dependence on $n$ of the deterministic self-consistent measure and corresponding Stieltjes-transform. For  each $n\in \mathbb N$ and $z\in\mathbb C\setminus\mathbb R$, $\mathfrak m(z)$ is characterized
as the solution $m$ of $z=f(m)$,
satisfying $\Im(z)\Im(m)>0$,
where
\begin{equation}
\label{eq:self_consistent_equation}
    f(m)
    \coloneqq
    -\frac1m
    +
    \frac1n
    \operatorname{tr}
    \left\{
        \Sigma_n
        (I_p+m\Sigma_n)^{-1}
    \right\}.
\end{equation}
Moreover,
\[
    \mathfrak m(z)
    =
    \int_{[0,\infty)}
    \frac1{x-z}\,d\varrho(x),
    \qquad
    z\in\mathbb C\setminus\mathbb R.
\]
For further details on the self-consistent equation \eqref{eq:self_consistent_equation}, we refer to
\cite[Section~2]{LinPan2026} and the references therein.

To relate this notation to the asymptotic Marchenko--Pastur quantities
used later, recall that
\[
    \frac1n
    \operatorname{tr}
    \left\{
        \Sigma_n
        (I_p+m\Sigma_n)^{-1}
    \right\}
    =
    \gamma_n
    \int_{[0,\infty)}
    \frac{t}{1+tm}\,dH_n(t),
    \qquad
    \gamma_n=\frac pn.
\]
Hence,
\[
    f(m)
    =
    -\frac1m
    +
    \gamma_n
    \int_{[0,\infty)}
    \frac{t}{1+tm}\,dH_n(t),
\]
which is the finite-$n$ analogue of the inverse Marchenko--Pastur map
\[
    x_{\gamma,H}(v)
    =
    -\frac1v
    +
    \gamma
    \int_{[0,\infty)}
    \frac{t}{1+tv}\,dH(t)
\]
used later in Section~\ref{Out-of-sample prediction risk}; see also
\cite[Equation~(1.4)]{SilversteinChoi1995}.
For readability, we otherwise follow the notation of
\cite{LinPan2026} in the present section.

\subsection{Lemma~\ref{lem:verification_GR_assumptions} and its proof}

\begin{lemma}\label{lem:verification_GR_assumptions}
Consider the setting in Proposition~\ref{prop:GreenRomanov-TV-counterexample}. Then Assumptions~1--6 of \cite{GreenRomanov2025} are satisfied.   
\end{lemma}

\begin{proof}
We set $m_n$ in Assumption~1 in \cite{GreenRomanov2025} to be $m_n=p/2$, so that $p/n=m_n/p=1/2$.
Since $Z_n$ has i.i.d.\ standard Gaussian entries, $S_n$ has full rank almost surely by \cite[Theorem~1]{Okamoto1973}, and
\[
    H_n=G_n=H=G
    =
    \frac12\delta_1+\frac12\delta_{1+\kappa},
\]
Assumptions~1--5 of \cite{GreenRomanov2025}
are satisfied with $k_0=0$.

To verify Assumption~6 of \cite{GreenRomanov2025}, we note that for $x\notin\operatorname{supp}(H)$ and
$\tau_\star\in\{1,1+\kappa\}$, we have
\[
\int
\frac{\tau^2}{(x-\tau)^2}
\,dH(\tau)
\geq
\frac{\tau_\star^2}{2(x-\tau_\star)^2}.
\]
Since $\tau_\star>0$, the right-hand side tends to infinity as
$x\to\tau_\star$ with
$x\notin\operatorname{supp}(H)$. Consequently,
\[
\lim_{\substack{x\to\tau_\star\\
x\notin\operatorname{supp}(H)}}
\int
\frac{\tau^2}{(x-\tau)^2}
\,dH(\tau)
=
\infty
\neq
\frac1\gamma=2.
\]
Hence, the sufficient condition for Assumption~6 in \cite[Equation~(18)]{GreenRomanov2025} is satisfied.
\end{proof}

\subsection{Lemma~\ref{lem:verification-assump-LinPan} and its proof}

\begin{lemma}\label{lem:verification-assump-LinPan}
Consider the setting in Proposition~\ref{prop:GreenRomanov-TV-counterexample}. There exists
$\kappa_0\in(0,1/2)$ such that, for every fixed
$\kappa\in(0,\kappa_0)$, Assumptions~1--6 of
\cite{LinPan2026} are satisfied.

Moreover, the positive part of
the support of the deterministic self-consistent spectral measure
$\varrho=\varrho_\kappa$ consists of a single compact interval:
\[
\operatorname{supp}(\varrho_\kappa)\cap(0,\infty)
=
[a_-(\kappa),a_+(\kappa)],
\qquad
0<a_-(\kappa)<a_+(\kappa)<\infty.
\]
\end{lemma}

\begin{proof}[Proof of Lemma~\ref{lem:verification-assump-LinPan}]
We start by translating the present model into the notation of
\cite{LinPan2026}. Set
\[
M=p,
\qquad
N=n,
\qquad
Y_n
\coloneqq
X_n^\top
=
\Sigma_n^{1/2}Z_n^\top.
\]
The matrix denoted by $X$ in Assumption~2 of
\cite{LinPan2026} is therefore $n^{-1/2}Z_n^\top$.
We will specify the constant $\tau>0$ appearing in
Assumptions~1--5 later. 

Assumption~1  in \cite{LinPan2026} requires
\[
\tau
\leq
\frac MN
\leq
\tau^{-1}
\]
uniformly in $N$. Since
\[
\frac MN
=
\frac pn
=
\frac12,
\]
this condition is satisfied for every $\tau\in(0,1/2]$.
The entries of $n^{-1/2}Z_n^\top$ clearly satisfy Assumption 2 of \cite{LinPan2026}.
The eigenvalues of $\Sigma_n$ are $1$ and $1+\kappa$, each with
multiplicity $p/2$, and hence,
\[
\|\Sigma_n\|_2
=
1+\kappa.
\]
Consequently, Assumption~3 is satisfied whenever
\[
\tau
\leq
\frac1{1+\kappa}.
\]
Moreover, if $\tau<1$, then $\mathbf 1_{\{\lambda_i(\Sigma_n)\leq\tau\}}=0$, for $1\leq i\leq p$.
It follows that
\[
\frac1M
\sum_{i=1}^M
\mathbf 1_{\{\lambda_i(\Sigma_n)\leq\tau\}}
=
0
\leq
1-\tau,
\]
which verifies Assumption~4.

It remains to verify Assumptions~5 and~6. To verify Assumption~5, we
first determine all real critical points of the self-consistent function
$f_\kappa$ associated with $\Sigma_n=I_p+\kappa P_n$.
We show that, for $\kappa\in(0,1/2)$, the derivative $f_\kappa'$ has
exactly two real zeros, one on $(-\infty,-1)$ and one on
$(-(1+\kappa)^{-1},0)$, and that both are simple. For $\kappa=0$,
these critical points and the corresponding critical values can be
computed explicitly. Using the implicit function theorem and
continuity in $\kappa$, we then show that, for all sufficiently small
fixed $\kappa>0$, the lower critical value remains bounded away from
zero, the two critical values remain separated from each other, and the
critical points remain separated from the poles of $f_\kappa$. The
characterization of the spectral edges as the critical values of
$f_\kappa$ then yields Assumption~5. The self-consistent function from \eqref{eq:self_consistent_equation}
in our setting is
\begin{align*}
f_\kappa(m)
&=
-\frac1m
+
\frac1n
\operatorname{tr}
\left\{
\frac{\Sigma_n}{I_p+m\Sigma_n}
\right\} 
=
-\frac1m
+
\frac14
\left(
\frac1{1+m}
+
\frac{1+\kappa}{1+(1+\kappa)m}
\right).
\end{align*}
The derivative is given by
\begin{align}\label{eq:f_kappa_prime}
f_\kappa'(m)
=
\frac1{m^2}
-
\frac14
\left(
\frac1{(1+m)^2}
+
\frac{(1+\kappa)^2}
{\{1+(1+\kappa)m\}^2}
\right),
\end{align}
and the negative poles of $f_\kappa$ are $-1$ and
$-(1+\kappa)^{-1}$.
By \cite[Lemma~2.4]{KnowlesYin2017}, $f_\kappa'$
has exactly one zero on
$(-\infty,-1)\cup(0,\infty)$ and exactly one zero on
$(-(1+\kappa)^{-1},0)$. We show that no additional zeros
occur between the two poles or on $(0,\infty)$.

To this end, let $-1<m<-(1+\kappa)^{-1}$
and set $t=-m$. Then $(1+\kappa)^{-1}<t<1$.
If $\kappa<1/2$, then $t>2/3$, and therefore $t>2(1-t)$.
It follows that
\[
\frac1{(1+m)^2}
=
\frac1{(1-t)^2}
>
\frac4{t^2}
=
\frac4{m^2}.
\]
Furthermore, for $-1<m<-(1+\kappa)^{-1}$, the denominator
$1+(1+\kappa)m$ is nonzero. Hence,
\[
\frac{(1+\kappa)^2}
{\{1+(1+\kappa)m\}^2}
>
0,
\]
and therefore
\begin{align*}
f_\kappa'(m)
&<
\frac1{m^2}
-
\frac1{4(1+m)^2}<0.
\end{align*}
Thus $f_\kappa'$ has no zero between its two poles.

Finally, let $m>0$. Then
\[
0<
\frac{m}{1+m}
<
1
\qquad\text{and}\qquad
0<
\frac{(1+\kappa)m}{1+(1+\kappa)m}
<
1.
\]
Multiplying $f_\kappa'(m)$ by $m^2>0$, we therefore obtain
\begin{align*}
m^2f_\kappa'(m)
=
1-\frac14
\left[
\left(
\frac{m}{1+m}
\right)^2
+
\left(
\frac{(1+\kappa)m}{1+(1+\kappa)m}
\right)^2
\right]
>
1-\frac14(1+1)
=
\frac12.
\end{align*}
Thus $f_\kappa'(m)>0$, and hence $f_\kappa'$ has no zero
on $(0,\infty)$. Consequently, $f_\kappa'$ has exactly two real zeros.
Both are simple, since \cite[Lemma~2.4]{KnowlesYin2017}
counts degenerate critical points twice.

We next consider the case $\kappa=0$. In this case,
\[
f_0(m)
=
-\frac1m+\frac1{2(1+m)}
\]
and
\begin{align}\label{eq:critical-points}
f_0'(m)
=
\frac1{m^2}
-
\frac1{2(1+m)^2}
=
\frac{2(1+m)^2-m^2}
{2m^2(1+m)^2}
=
\frac{m^2+4m+2}
{2m^2(1+m)^2}.
\end{align}
Since for $m\notin\{-1,0\}$, the denominator in
\eqref{eq:critical-points} is strictly positive,
the critical points are the roots of
\[
m^2+4m+2=0.
\]
Consequently, the two critical points are
\[
m_-^{(0)}
=
-(2+\sqrt2)
\qquad\text{and}\qquad
m_+^{(0)}
=
-(2-\sqrt2).
\]
Moreover, both critical points are simple. Indeed,
the numerator in \eqref{eq:critical-points} is
\[
q(m)=m^2+4m+2,
\qquad
q'(m)=2m+4,
\]
and hence,
\[
q'(m_-^{(0)})
=
-2\sqrt2\neq0,
\qquad
q'(m_+^{(0)})
=
2\sqrt2\neq0.
\]
We now compute the corresponding critical values. Since
\[
a_-^{(0)}
\coloneqq
f_0(m_-^{(0)})
=
\frac{2-\sqrt2}{2}
-
\frac{\sqrt2-1}{2}
=
\frac32-\sqrt2
=
\left(
1-\frac1{\sqrt2}
\right)^2.
\]
Similarly,
\[
a_+^{(0)}
\coloneqq
f_0(m_+^{(0)})
=
\frac{2+\sqrt2}{2}
+
\frac{\sqrt2+1}{2}
=
\frac32+\sqrt2
=
\left(
1+\frac1{\sqrt2}
\right)^2.
\]
In particular,
$a_-^{(0)},a_+^{(0)}-a_-^{(0)}>0$.
Moreover, the critical points are separated from
the pole $-1$:
\[
\min\left\{
|1+m_-^{(0)}|,
|1+m_+^{(0)}|
\right\}
=
\sqrt2-1>0.
\]
Since the zeros $m_-^{(0)}$ and $m_+^{(0)}$ of $f_0'$
are simple, the implicit function theorem yields continuous
functions $\kappa\mapsto m_-(\kappa)$ and
$\kappa\mapsto m_+(\kappa)$ for all sufficiently small
$\kappa\geq0$, satisfying
\[
f_\kappa'(m_\pm(\kappa))=0
\qquad\text{and}\qquad
m_\pm(0)=m_\pm^{(0)}.
\]
The corresponding critical values
$a_-(\kappa)=f_\kappa(m_-(\kappa))$ and
$a_+(\kappa)=f_\kappa(m_+(\kappa))$
are therefore continuous at $\kappa=0$. Since
\[
a_-^{(0)}>0,
\qquad
a_+^{(0)}-a_-^{(0)}>0,
\qquad
\min\left\{
|1+m_-^{(0)}|,
|1+m_+^{(0)}|
\right\}>0,
\]
continuity implies that there exist $\kappa_0>0$ and $c>0$
such that, for every $\kappa\in[0,\kappa_0]$,
\[
a_-(\kappa)\geq c,
\qquad
a_+(\kappa)-a_-(\kappa)\geq c,
\]
and
\[
\min_{\sigma\in\{1,1+\kappa\}}
\min_{\nu\in\{-,+\}}
|1+\sigma m_\nu(\kappa)|
\geq c.
\]
By the characterization of the spectral edges as the critical
values of $f_\kappa$ (see \cite[Section~7.1]{LinPan2026}),
the positive part of the support of the deterministic
measure $\varrho_\kappa$ is
\[
\operatorname{supp}(\varrho_\kappa)\cap(0,\infty)
=
[a_-(\kappa),a_+(\kappa)],
\]
which establishes the second statement of the lemma.
Thus, in the notation of \cite{LinPan2026}, $K=1$,
$\mathfrak a_1=a_+(\kappa)$,
$\mathfrak a_2=a_-(\kappa)$, and
$\mathfrak m_\kappa(\mathfrak a_1)=m_+(\kappa)$,
$\mathfrak m_\kappa(\mathfrak a_2)=m_-(\kappa)$.
Choose $\tau>0$ sufficiently small such that
\[
\tau
<
\min
\left\{
\frac12,
\frac1{1+\kappa},c
\right\}.
\]
With this choice of $\tau$, the preceding bounds verify
Assumption~5 of \cite{LinPan2026}.

It remains to verify Assumption~6. Since $\Sigma_n$ has
two distinct positive eigenvalues with fixed positive spectral
weights, and the critical points of $f_\kappa$ are simple
with distinct critical values,
\cite[Example~2.8]{KnowlesYin2017} implies that the bulk
component $[a_-(\kappa),a_+(\kappa)]$ is regular.
In particular, for every $\tau'>0$ for which
$[a_-(\kappa)+\tau',a_+(\kappa)-\tau']\neq\varnothing$,
there exists $c_{\tau'}>0$ such that
\[
\min_{E\in[a_-(\kappa)+\tau',a_+(\kappa)-\tau']}
\varrho_\kappa(E)
\geq c_{\tau'}.
\]
Hence Assumption~6 of \cite{LinPan2026} is satisfied.
This concludes the proof.

\end{proof}

\subsection{Lemma~\ref{lem:distinct-eigenvalues} and its proof}

\begin{lemma}\label{lem:distinct-eigenvalues}
Under the assumptions of
Proposition~\ref{prop:GreenRomanov-TV-counterexample}, the eigenvalues of $S_n$ are strictly positive and pairwise
distinct almost surely. Furthermore, for any unit vector $v\in \mathbb R^p$, we have
\begin{align*}
\|\widehat{\mathsf K}_n^v\|_{\mathrm{TV}}
=
\sum_{i,j=1}^p
|\langle v,w_i\rangle|
\left|
\langle w_i,\Sigma_nw_j\rangle
\right|
|\langle w_j,v\rangle|
\end{align*}
\end{lemma}
\begin{proof}[Proof of Lemma~\ref{lem:distinct-eigenvalues}]
    Since $\Sigma_n$ is positive definite, we have $\operatorname{rank}(\Sigma_n)=p$ for each $n$.
Moreover, the joint distribution of the $p$ column vectors of $Z_n$
is absolutely continuous with respect to the $np$-dimensional
Lebesgue measure. Therefore, by \cite[Theorem~1]{Okamoto1973},
\[
\operatorname{rank}(Z_n\Sigma_n Z_n^\top)
=
\min\{n,p\}
=
p
\]
and the nonzero eigenvalues of $Z_n\Sigma_n Z_n^\top$ are pairwise
distinct almost surely, for every $n$.

Since $nS_n$ and $Z_n\Sigma_n Z_n^\top$ share the same nonzero eigenvalues, we conclude that $\operatorname{rank}(S_n)=p$ and the eigenvalues of $S_n$ are strictly positive and pairwise
distinct almost surely. We may therefore arrange them as $\lambda_1(S_n)>\cdots>\lambda_p(S_n)>0$.
It follows that the atoms
\[
\bigl(\lambda_i(S_n),\lambda_j(S_n)\bigr),
\qquad
1\leq i,j\leq p,
\]
are pairwise distinct almost surely.
Thus, 
\begin{align*}
\|\widehat{\mathsf K}_n^v\|_{\mathrm{TV}}
=
\sum_{i,j=1}^p
|\langle v,w_i\rangle|
\left|
\langle w_i,\Sigma_nw_j\rangle
\right|
|\langle w_j,v\rangle|
\end{align*}
almost surely.
\end{proof}

\subsection{Proof of Lemma~\ref{lem:bulk-overlap-control}}

\begin{proof}[Proof of Lemma~\ref{lem:bulk-overlap-control}]
Since $\Sigma_n=I_p+\kappa P_n$,
we have
\[
\langle w_i,\Sigma_nw_j\rangle
=
\delta_{ij}+\kappa c_{ij}.
\]
In particular, whenever \(i\neq j\), $\left|
\langle w_i,\Sigma_nw_j\rangle
\right|
=
\kappa|c_{ij}|$.
We first show that there exists a constant \(c_0>0\) such that
\begin{align} \label{eq:goal1}
\mathbb P\!\left(\sum_{i\ne j}c_{ij}^2\ge c_0p\right)\longrightarrow 1.
\end{align}

Write $Z_n=(Z_{n,1},Z_{n,2})$,
where $Z_{n,1},Z_{n,2}$ take values in $\mathbb R^{n\times p/2}$.
By the choice of $\Sigma_n$ in Proposition \ref{prop:GreenRomanov-TV-counterexample}, we have
\[
S_n
=
\begin{pmatrix}
n^{-1}Z_{n,1}^\top Z_{n,1}
&
\sqrt{1+\kappa}\,n^{-1}Z_{n,1}^\top Z_{n,2}
\\
\sqrt{1+\kappa}\,n^{-1}Z_{n,2}^\top Z_{n,1}
&
(1+\kappa)n^{-1}Z_{n,2}^\top Z_{n,2}
\end{pmatrix}.
\]
Consequently,
\[
S_nP_n-P_nS_n
=
\begin{pmatrix}
0
&
\sqrt{1+\kappa}\,n^{-1}Z_{n,1}^\top Z_{n,2}
\\
-\sqrt{1+\kappa}\,n^{-1}Z_{n,2}^\top Z_{n,1}
&
0
\end{pmatrix},
\]
and hence,
\begin{align}\label{eq:Frobenius-bound}
\|S_nP_n-P_nS_n\|_F^2
=
2(1+\kappa)
\left\|
\frac1nZ_{n,1}^\top Z_{n,2}
\right\|_F^2.
\end{align}
Conditional on $Z_{n,1}$, the $p/2$ columns of $Z_{n,2}$ are
independent standard Gaussian vectors in $\mathbb R^n$. Therefore,
\[
\mathbb E\left[
\left.
\left\|
\frac1nZ_{n,1}^\top Z_{n,2}
\right\|_F^2
\,\right|\,
Z_{n,1}
\right]
=
\frac{p}{2n^2}\|Z_{n,1}\|_F^2
\]
and
\[
\operatorname{Var}\left(
\left.
\left\|
\frac1nZ_{n,1}^\top Z_{n,2}
\right\|_F^2
\,\right|\,
Z_{n,1}
\right)
=
\frac{p}{n^4}
\operatorname{tr}
\left\{
(Z_{n,1}Z_{n,1}^\top)^2
\right\}.
\]
Moreover,
\[
\mathbb E\left(\|Z_{n,1}\|_F^2\right)
=
\frac{np}{2},
\qquad
\operatorname{Var}(\|Z_{n,1}\|_F^2)
=
np,
\]
and
\[
\mathbb E\left(
\operatorname{tr}\bigl\{
(Z_{n,1}Z_{n,1}^\top)^2\bigr\}\right)
=
\frac{np}{2}
\left(
n+\frac p2+1
\right).
\]
Hence, the law of total expectation gives
\[
\mathbb E
\left(
\left\|
\frac1nZ_{n,1}^\top Z_{n,2}
\right\|_F^2
\right)
=
\mathbb E\left[
\mathbb E\left[
\left.
\left\|
\frac1nZ_{n,1}^\top Z_{n,2}
\right\|_F^2
\,\right|\,
Z_{n,1}
\right]
\right]
=\frac{p}{2n^2}\frac{np}{2}
=
\frac{p^2}{4n}.
\] 
Furthermore, by the law of total variance we get
\begin{align*}
\operatorname{Var}
\left(
\left\|
\frac1nZ_{n,1}^\top Z_{n,2}
\right\|_F^2
\right)
&=
\mathbb E\left[
\operatorname{Var}\left(
\left.
\left\|
\frac1nZ_{n,1}^\top Z_{n,2}
\right\|_F^2
\,\right|\,
Z_{n,1}
\right)
\right]
+
\operatorname{Var}\left(
\mathbb E\left[
\left.
\left\|
\frac1nZ_{n,1}^\top Z_{n,2}
\right\|_F^2
\,\right|\,
Z_{n,1}
\right]
\right)
\\
&=
\frac{p^2(n+p+1)}{2n^3}.
\end{align*}
Since \(n=2p\), it follows that
\[
\mathbb E
\left[
\frac1p
\left\|
\frac1nZ_{n,1}^\top Z_{n,2}
\right\|_F^2
\right]
=
\frac18\quad\text{and}\quad
\operatorname{Var}
\left(
\frac1p
\left\|
\frac1nZ_{n,1}^\top Z_{n,2}
\right\|_F^2
\right)
=
O(p^{-2}).
\]
Hence, we conclude that,
\[
\frac1p
\left\|
\frac1nZ_{n,1}^\top Z_{n,2}
\right\|_F^2
\overset{\mathbb{P}}{\longrightarrow}
\frac18.
\]
Consequently, using \eqref{eq:Frobenius-bound} and the preceding display, we have
\[
\frac1p
\|S_nP_n-P_nS_n\|_F^2
\overset{\mathbb{P}}{\longrightarrow}
\frac{1+\kappa}{4}.
\]
This implies
\begin{align}\label{eq:step1} 
\mathbb P\left( \|S_nP_n-P_nS_n\|_F^2 \geq \frac{1+\kappa}{8}p \right) \longrightarrow1. 
\end{align} 
On the other hand, 
\[
\|S_nP_n-P_nS_n\|_F^2
=
\sum_{i,j=1}^p
\bigl(\lambda_i(S_n)-\lambda_j(S_n)\bigr)^2
c_{ij}^2=\sum_{\substack{1\leq i,j\leq p\\i\neq j}}\bigl(\lambda_i(S_n)-\lambda_j(S_n)\bigr)^2
c_{ij}^2.
\]
Since $S_n$ is positive definite by Lemma~\ref{lem:distinct-eigenvalues}, $|\lambda_i(S_n)-\lambda_j(S_n)|
\leq
\|S_n\|_2$,
and hence
\begin{align}\label{eq:Frobenius-estimate}
\|S_nP_n-P_nS_n\|_F^2
\leq
\|S_n\|_2^2
\sum_{\substack{1\leq i,j\leq p\\i\neq j}}
c_{ij}^2,\qquad\text{ almost surely}.
\end{align}
Moreover, using \cite[Theorem~3.1]{BaiYin1988} and since
\[
\|S_n\|_2
\leq
\|\Sigma_n\|_2
\left\|
\frac1nZ_n^\top Z_n
\right\|_2
=
(1+\kappa)
\left\|
\frac1nZ_n^\top Z_n
\right\|_2,
\]
we may choose a deterministic constant $0<M<\infty$ such that
\begin{align}\label{eq:step2}
\mathbb P\bigl(\|S_n\|_2\leq M\bigr)
\longrightarrow1.
\end{align}
On the intersection of the events considered in \eqref{eq:step1} and \eqref{eq:step2} the estimate in \eqref{eq:Frobenius-estimate}
implies that
\[ 
\sum_{\substack{1\leq i,j\leq p\\i\neq j}} c_{ij}^2 \geq \frac{1+\kappa}{8M^2}p. 
\] 
Therefore, setting $c_0\coloneqq (1+\kappa)(8M^2)^{-1}>0$, we obtain 
\begin{align*} 
\mathbb P\left( \sum_{\substack{1\leq i,j\leq p\\i\neq j}} c_{ij}^2 \geq c_0p \right) \geq 1- \mathbb P\left( \|S_nP_n-P_nS_n\|_F^2 
< \frac{1+\kappa}{8}p \right) - \mathbb P\bigl(\|S_n\|_2>M\bigr) \longrightarrow1,
\end{align*} 
which completes the proof of \eqref{eq:goal1}.

We next show \eqref{eq:definition_A_n}, that is, that the lower bound in \eqref{eq:goal1}
of order $p$ still holds with probability tending to one when the
summation is restricted to $i,j\in J_n$. Consider the decomposition in \eqref{eq:sam-cov-decomposition} and write $W_n=(w_1,\ldots,w_p)$. Since \(P_n\) is an orthogonal projection, the matrix $W_n^\top P_nW_n=(c_{ij})_{i,j=1}^p$ is also an orthogonal projection. Therefore,
\begin{align}\label{eq:orthoproj}
\sum_{j=1}^pc_{ij}^2
=
c_{ii}
\leq1
\end{align}
for every \(i\).

Since $|J_n^c|\leq2\eta p+2$ and using \eqref{eq:orthoproj},
we have
\begin{align*}
\sum_{\substack{1\leq i,j\leq p\\i\neq j}}
c_{ij}^2
-
\sum_{\substack{i,j\in J_n\\i\neq j}}
c_{ij}^2
\leq
\sum_{i\notin J_n}\sum_{j=1}^pc_{ij}^2
+
\sum_{j\notin J_n}\sum_{i=1}^pc_{ij}^2
\leq
2|J_n^c|
\leq
4\eta p+4.
\end{align*}
Choosing $\eta\in (0,1/2)$ such that $4\eta<c_0/2$, on the event \[\left\{
\sum_{\substack{1\leq i,j\leq p\\i\neq j}}
c_{ij}^2
\geq
c_0p\right\}
\]
we have
\[
\sum_{\substack{i,j\in J_n\\i\neq j}}
c_{ij}^2
\geq
\sum_{\substack{1\leq i,j\leq p\\i\neq j}}
c_{ij}^2
-
4\eta p-4
\geq
(c_0-4\eta)p-4
\geq
\frac{c_0}{2}p-4.
\]
For all sufficiently large $n$, $\frac{c_0}{2}p-4
\geq
\frac{c_0}{4}p$. 
Thus, setting $c_1\coloneqq\frac{c_0}{4}>0$
and
\begin{align*}
A_n
\coloneqq
\left\{
\sum_{\substack{i,j\in J_n\\i\neq j}}
c_{ij}^2
\geq
c_1p
\right\},
\end{align*}
we have 
\[
\left\{\sum_{\substack{1\leq i,j\leq p\\i\neq j}}
c_{ij}^2
\geq
c_0p\right\}\subseteq\left\{
\sum_{\substack{i,j\in J_n\\i\neq j}}
c_{ij}^2
\geq
c_1p
\right\},
\]
for $n$ large enough.
Therefore, $\mathbb P(A_n)\longrightarrow1$, which proves the first claim \eqref{eq:definition_A_n} of the lemma.

To show the second claim \eqref{eq:outline-step2}, we apply
\cite[Corollary~2.3]{LinPan2026} with the deterministic matrix $P_n$. Note that, by Lemma~\ref{lem:verification-assump-LinPan}, Assumptions~1--6
of \cite{LinPan2026} hold.
We first recall the notion of stochastic domination used in
\cite{LinPan2026}. Let $Y\equiv Y^{(N)}(t)$, $Z\equiv Z^{(N)}(t)$,
and 
$t\in T^{(N)}$,
be two families of nonnegative random variables. The notation
\[
Y\prec Z
\]
means that, for every fixed $\delta>0$ and $L>0$,
\[
\sup_{t\in T^{(N)}}
\mathbb P\left(
Y^{(N)}(t)>N^\delta Z^{(N)}(t)
\right)
\leq
N^{-L},\qquad\text{for all } N\geq N_0(\delta,L).
\]

We now translate the notation of \cite{LinPan2026} to the present
setting. We have
\[
M=p,
\qquad
N=n,
\qquad
Y_n=X_n^\top.
\]
Consequently, the left singular vectors of $Y_n$ are precisely the
eigenvectors $w_1,\ldots,w_p$ of $S_n$, with the same ordering as in
\eqref{eq:sam-cov-decomposition}.

Using the support analysis from Lemma~\ref{lem:verification-assump-LinPan}, the positive part of the support of the deterministic measure $\varrho_\kappa$ consists of a single bulk component. Hence, in the
notation of \cite{LinPan2026}, $K=1$.
Since all $M=p$ left singular vectors correspond to this unique bulk,
we have
\[
N_1=M=p,
\qquad
u_{1,i}=w_i,
\qquad
1\leq i\leq p.
\]
Therefore,
\[
\mathfrak n_{1,i}
=
i\wedge(N_1+1-i)
=
i\wedge(p+1-i).
\]
For notational simplicity, we set $\mathfrak n_i
\coloneqq
\mathfrak n_{1,i}
=
i\wedge(p+1-i)$.

Since $P_n$ is deterministic and $\|P_n\|_2=1$,
it is an admissible choice for the matrix $D_1$ in
\cite[Corollary~2.3]{LinPan2026}. The deterministic centering term in
\cite[Equation~(2.10a)]{LinPan2026} is multiplied by
$\delta_{k\ell}\delta_{ij}$. In the present one-bulk setting,
$k=\ell=1$  this factor reduces to $\delta_{ij}$. Consequently, the
centering term vanishes whenever $i\neq j$, and
\cite[Equation~(2.10a)]{LinPan2026} yields, uniformly for
$1\leq i,j\leq p$ with $i\neq j$,
\[
|c_{ij}|
=
|\langle w_i,P_nw_j\rangle|
\prec
\bigl(
n\mathfrak n_i\mathfrak n_j
\bigr)^{-1/6}.
\]
For every $i\in J_n$, we have $i
\geq
\lceil\eta p\rceil
\geq
\eta p$
and
\[
p+1-i
\geq
p+1-\lfloor(1-\eta)p\rfloor
\geq
\eta p.
\]
It follows that
\[
\mathfrak n_i
=
i\wedge(p+1-i)
\geq
\eta p,
\qquad
i\in J_n.
\]
Since $n=2p$, for every $i,j\in J_n$ we therefore have
\[
\bigl(
n\mathfrak n_i\mathfrak n_j
\bigr)^{-1/6}
\leq
\bigl\{
2p(\eta p)^2
\bigr\}^{-1/6}
=
2^{-1/6}\eta^{-1/3}p^{-1/2}.
\]
Thus, setting $C_\eta
\coloneqq
2^{-1/6}\eta^{-1/3}$,
we obtain, uniformly for $i,j\in J_n$ with $i\neq j$,
\[
|c_{ij}|
\prec
C_\eta p^{-1/2}.
\]
We now derive the asserted bound for the maximum. Fix
$\varepsilon\in(0,1/2)$ and set $\delta
\coloneqq
\frac{\varepsilon}{2}$.
By the definition of stochastic domination, choosing, for example,
$L=4$, we have, for all sufficiently large $n$,
\[
\max_{\substack{i,j\in J_n\\i\neq j}}
\mathbb P\left(
|c_{ij}|
>
C_\eta n^\delta p^{-1/2}
\right)
\leq
n^{-4}.
\]
Since $n=2p$, $C_\eta n^\delta p^{-1/2}
=
C_\eta2^\delta p^{-1/2+\delta}$.
Because $\delta=\varepsilon/2<\varepsilon$, for all sufficiently large
$p$, 
\[
C_\eta2^\delta p^{-1/2+\delta}
\leq
p^{-1/2+\varepsilon}.
\]
Consequently, for all sufficiently large $n$,
\[
\max_{\substack{i,j\in J_n\\i\neq j}}
\mathbb P\left(
|c_{ij}|
>
p^{-1/2+\varepsilon}
\right)
\leq
\frac{1}{n^4}.
\]
A union bound over the at most $p^2$ pairs $(i,j)$ now gives
\[
\mathbb P\left(
\max_{\substack{i,j\in J_n\\i\neq j}}
|c_{ij}|
>
p^{-1/2+\varepsilon}
\right)
\leq
\sum_{\substack{i,j\in J_n\\i\neq j}}
\mathbb P\left(
|c_{ij}|
>
p^{-1/2+\varepsilon}
\right)
\leq
\frac{p^2}{n^4}
=
\frac{p^2}{(2p)^4}
\longrightarrow
0.
\]
The second claim \eqref{eq:outline-step2} follows.
\end{proof}

\subsection{Proof of Lemma~\ref{lem:Gaussian-retention}}

\begin{proof}[Proof of Lemma~\ref{lem:Gaussian-retention}]
Fix an arbitrary $\varepsilon\in (0,1/2)$, and let
$
g_n'\sim N(0,I_p),
$
such that $g_n, g_n'$, and $Z_n$ are mutually independent of each other. Conditional
on $S_n$, the random variables
\[
\widetilde{g}_i
\coloneqq
\langle g_n,w_i\rangle,
\qquad
1\leq i\leq p,
\]
are independent Gaussian random variables.
For convenience write \[
\xi_i
=
\mathbf 1_{\{1\leq|\widetilde g_i|\leq2\}},
\qquad
1\leq i\leq p,
\]
and
\[
R_n
=
\sum_{\substack{i,j\in J_n\\i\neq j}}
\xi_i\xi_jc_{ij}^2.
\]
Conditional on $S_n$, the random variables
$\xi_1,\ldots,\xi_p$ are independent. 
Similarly, we define $\xi_1', \dots, \xi_p'$ by replacing $g_n$ with $g_n'$, and obtain $R_n^{(k)}$ from
$R_n$ by replacing $\xi_k$ with $\xi_k'$. 
If $k\notin J_n$,
then $R_n=R_n^{(k)}$, whereas, if $k\in J_n$, then, using
$c_{ij}=c_{ji}$,
\[
\begin{aligned}
R_n-R_n^{(k)}
&=
2(\xi_k-\xi_k')
\sum_{\substack{j\in J_n\\j\neq k}}
\xi_jc_{kj}^2.
\end{aligned}
\]
Consequently,
\[
\begin{aligned}
|R_n-R_n^{(k)}|
&\leq
2\sum_{j=1}^pc_{kj}^2
=
2c_{kk}
\leq2,
\end{aligned}
\]
where we used \eqref{eq:orthoproj}.
Note that $R_n^{(k)}$ and $R_n$ are independent given $\sigma(S_n, \xi_1,\dots, \xi_{k-1}, \xi_{k+1}, \dots, \xi_p)$.
Therefore, the conditional Efron--Stein inequality gives
\begin{align}\label{eq:conditional_ES}
\operatorname{Var}(R_n\mid S_n)
\leq
\frac12
\sum_{k=1}^p
\mathbb E\left[
\left(R_n-R_n^{(k)}\right)^2
\,\middle|\,
S_n
\right]
\leq
\frac12\sum_{k=1}^p4
=
2p.
\end{align}
Conditional on $S_n$, the random variables
$\xi_1,\ldots,\xi_p$ are independent and satisfy
\[
\mathbb E[\xi_i\mid S_n]
=
\mathbb P\bigl(1\leq|Y|\leq2\bigr)\eqqcolon q,\qquad Y\sim\mathcal{N}(0,1).
\]
Since the coefficients $c_{ij}^2$ are measurable with respect to the $\sigma$-field generated by
$S_n$, it follows that
\[
\mathbb E[R_n\mid S_n]
=
q^2
\sum_{\substack{i,j\in J_n\\i\neq j}}
c_{ij}^2.
\]
Set $c_2
\coloneqq
\frac{c_1q^2}{2}>0$ and consider the event in \eqref{eq:definition_A_n}, which we denote by $A_n$.  On the event $\{R_n<c_2p\}\cap A_n$, we have
\[
\left|
R_n-\mathbb E[R_n\mid S_n]
\right|
\geq
\mathbb E[R_n\mid S_n]-R_n
\\
>
q^2c_1p-\frac{q^2c_1p}{2}
=
c_2p.
\]
Consequently, conditional Chebyshev's inequality and
\eqref{eq:conditional_ES} give, almost surely,
\begin{align*}
\mathbf 1_{A_n}
\mathbb P\left(
R_n<c_2p
\,\middle|\,
S_n
\right)
&\leq
\mathbf 1_{A_n}
\mathbb P\left(
\left|
R_n-\mathbb E[R_n\mid S_n]
\right|
>
c_2p
\,\middle|\,
S_n
\right)
\leq
\mathbf 1_{A_n}
\frac{\operatorname{Var}(R_n\mid S_n)}
{c_2^2p^2}
\\
&\leq
\mathbf 1_{A_n}
\frac{2p}{c_2^2p^2}
=
\mathbf 1_{A_n}
\frac{8}{c_1^2q^4p}.
\end{align*}
Since $A_n$ is measurable with respect to $S_n$, the tower property
yields
\[
\mathbb P\left(
\{R_n<c_2p\}\cap A_n
\right)
=
\mathbb E\left[
\mathbf 1_{A_n}
\mathbb P\left(
R_n<c_2p
\,\middle|\,
S_n
\right)
\right]
\leq
\frac{8}{c_1^2q^4p}
\mathbb P(A_n)
\leq
\frac{8}{c_1^2q^4p}.
\]
Therefore,
\begin{align*}
\mathbb P(R_n<c_2p)
&=
\mathbb P\left(
\{R_n<c_2p\}\cap A_n
\right)
+
\mathbb P\left(
\{R_n<c_2p\}\cap A_n^c
\right)
\\
&\leq
\mathbb P(A_n^c)
+
\frac{8}{c_1^2q^4p}
\longrightarrow0,
\end{align*}
which yields
\begin{align}\label{eq:step4}
\mathbb P(R_n\geq c_2p)
\longrightarrow1.
\end{align}
Consequently, on the intersection of the event in \eqref{eq:outline-step2} with the event in \eqref{eq:step4}, we obtain
\[
\sum_{\substack{i,j\in J_n\\i\neq j}}
\xi_i\xi_j|c_{ij}|
\geq
\frac{R_n}{
\displaystyle
\max_{\substack{i,j\in J_n\\i\neq j}}
|c_{ij}|
}
\geq
c_2p^{3/2-\varepsilon},
\]
where the denominator in the second inequality is strictly positive, since $R_n\geq c_2p>0$.
Furthermore, since $\|g_n\|_2^2\sim \mathcal{X}_p^2$, it follows that $\mathbb P\left(\|g_n\|_2^2\leq2p\right)
\to1$.

Let \[
E_n
\coloneqq
\left\{
\|g_n\|_2^2\leq2p
\right\}
\cap
\left\{
R_n\geq c_2p
\right\}
\cap
\left\{
\max_{\substack{i,j\in J_n\\i\neq j}}
|c_{ij}|
\leq
p^{-1/2+\varepsilon}
\right\}.
\]
By the preceding estimates, it follows that
\[
\mathbb P(E_n)\longrightarrow1.
\]
Moreover, since $\langle w_i,\Sigma_nw_j\rangle
=
\delta_{ij}+\kappa c_{ij}$, and using Lemma~\ref{lem:distinct-eigenvalues}, we obtain
\[
\|\widehat{\mathsf K}_n^{u_n}\|_{\mathrm{TV}}
=
\frac1{\|g_n\|_2^2}
\sum_{i,j=1}^p
|\widetilde g_i|
\left|
\delta_{ij}+\kappa c_{ij}
\right|
|\widetilde g_j|.
\]
All summands are nonnegative. Therefore, restricting the sum to
$i,j\in J_n$ with $i\neq j$, for which $\delta_{ij}=0$, yields
\[
\|\widehat{\mathsf K}_n^{u_n}\|_{\mathrm{TV}}
\geq
\frac{\kappa}{\|g_n\|_2^2}
\sum_{\substack{i,j\in J_n\\i\neq j}}
|\widetilde g_i\widetilde g_j|\,|c_{ij}|.
\]
Since $
|\widetilde g_i\widetilde g_j|
\geq
\xi_i\xi_j$
for all $i,j$, we obtain on $E_n$
\begin{align*}
\|\widehat{\mathsf K}_n^{u_n}\|_{\mathrm{TV}}
&\geq
\frac{\kappa}{2p}
\sum_{\substack{i,j\in J_n\\i\neq j}}
\xi_i\xi_j|c_{ij}|
\geq
\frac{\kappa}{2p}
\frac{R_n}{
\displaystyle
\max_{\substack{i,j\in J_n\\i\neq j}}
|c_{ij}|
}
\\
&\geq
\frac{\kappa}{2p}
\frac{c_2p}{p^{-1/2+\varepsilon}}
=
\frac{\kappa c_2}{2}
p^{1/2-\varepsilon}.
\end{align*}
Thus, with $c_3\coloneqq2^{-1}\kappa c_2>0$,
we have
\[
\mathbb P(E_n)\leq \mathbb P\left(
\|\widehat{\mathsf K}_n^{u_n}\|_{\mathrm{TV}}
\geq
c_3p^{1/2-\varepsilon}
\right)
\longrightarrow1.
\] 
\end{proof}

\begin{proof}[Proof of Lemma~\ref{lem:Gaussian-block-spherical}]
Fix an arbitrary $\varepsilon\in (0,1/2)$. We first compare $u_n$ and $\overline u_n$. Set 
\[ 
a_n \coloneqq \frac{\|g_n^{(1)}\|_2^2}{\|g_n\|_2^2}. 
\] 
Since $\|g_n^{(1)}\|_2^2$ and $\|g_n^{(2)}\|_2^2$ are independent $\chi^2$-distributed random variables with \(p/2\) degrees of freedom, 
\[ 
a_n= \frac{\|g_n^{(1)}\|_2^2} {\|g_n^{(1)}\|_2^2+\|g_n^{(2)}\|_2^2} \sim \operatorname{Beta}\left(\frac p4,\frac p4\right). 
\] 
Therefore, 
\[ 
\mathbb E(a_n)=\frac12, \qquad \operatorname{Var}(a_n)=\frac1{2p+4}, 
\] 
and hence, 
\[ 
a_n-\frac12 = O_{\mathbb P}(p^{-1/2}). 
\]
Writing $u_n^{(1)} = g_n^{(1)}\|g_n^{(1)}\|_2^{-1}$, $u_n^{(2)} = g_n^{(2)}\|g_n^{(2)}\|_2^{-1}$, and observing that $\|g_n\|_2^2=\|g_n^{(1)}\|_2^2+\|g_n^{(2)}\|_2^2$, we have 
\[ 
u_n =\frac{g_n}{\|g_n\|_2}= \begin{pmatrix} \sqrt{a_n}\,u_n^{(1)} \\[1mm] 
\sqrt{1-a_n}\,u_n^{(2)} \end{pmatrix}, \qquad \overline u_n = \begin{pmatrix} 2^{-1/2}u_n^{(1)} \\[1mm] 2^{-1/2}u_n^{(2)} \end{pmatrix}. 
\] 
Since $\|u_n^{(1)}\|_2=\|u_n^{(2)}\|_2=1$, it follows that 
\[ 
\|u_n-\overline u_n\|_2^2 = \left(\sqrt{a_n}-\frac1{\sqrt2}\right)^2 + \left(\sqrt{1-a_n}-\frac1{\sqrt2}\right)^2. 
\] 
Moreover, 
\begin{align*}
\left| \sqrt{a_n}-\frac1{\sqrt2} \right| = \frac{|a_n-\frac12|} {\sqrt{a_n}+\frac1{\sqrt2}} \leq \sqrt2\left|a_n-\frac12\right| 
\end{align*}
and, similarly,
\begin{align*}
\left| \sqrt{1-a_n}-\frac1{\sqrt2} \right| \leq \sqrt2\left|a_n-\frac12\right|. 
\end{align*}
Consequently, 
\[ 
\|u_n-\overline u_n\|_2^2 \leq 4\left|a_n-\frac12\right|^2 = O_{\mathbb P}(p^{-1}), 
\] 
and therefore 
\begin{align} \label{eq:u-unbar-distance} 
\|u_n-\overline u_n\|_2 = O_{\mathbb P}(p^{-1/2}). \end{align}

Next, we show that the total variation is Lipschitz continuous in the
signal direction. Let $u,v\in\mathbb R^p$ be unit vectors. Since, by
Lemma~\ref{lem:distinct-eigenvalues}, the atoms
\[
\bigl(\lambda_i(S_n),\lambda_j(S_n)\bigr),
\qquad
1\leq i,j\leq p,
\]
are pairwise distinct almost surely, we have
\[
\left|
\|\widehat{\mathsf K}_n^{u}\|_{\mathrm{TV}}
-
\|\widehat{\mathsf K}_n^{v}\|_{\mathrm{TV}}
\right|
\leq
\sum_{i,j=1}^p
\left|
\langle w_i,\Sigma_nw_j\rangle
\right|
\left|
|\langle u,w_i\rangle|
|\langle u,w_j\rangle|
-
|\langle v,w_i\rangle|
|\langle v,w_j\rangle|
\right|.
\]
Using $|ab-cd|\leq|a-c||b|+|c||b-d|$,
we obtain
\begin{align}
\label{eq:Lipschitz}
\begin{split}
\left|
\|\widehat{\mathsf K}_n^{u}\|_{\mathrm{TV}}
-
\|\widehat{\mathsf K}_n^{v}\|_{\mathrm{TV}}
\right|
&\leq
\sum_{i,j=1}^p
\left|
\langle w_i,\Sigma_nw_j\rangle
\right|
\left|
|\langle u,w_i\rangle|
-
|\langle v,w_i\rangle|
\right|
|\langle u,w_j\rangle|
\\
&\qquad+
\sum_{i,j=1}^p
\left|
\langle w_i,\Sigma_nw_j\rangle
\right|
|\langle v,w_i\rangle|
\left|
|\langle u,w_j\rangle|
-
|\langle v,w_j\rangle|
\right|.
\end{split}
\end{align}

By the Cauchy--Schwarz inequality,
\begin{align*}
&\sum_{i,j=1}^p
\left|
\langle w_i,\Sigma_nw_j\rangle
\right|
\bigl|
|\langle u,w_i\rangle|
-
|\langle v,w_i\rangle|
\bigr|
|\langle u,w_j\rangle|\\
&\qquad\leq
\left(
\sum_{i,j=1}^p
\langle w_i,\Sigma_nw_j\rangle^2
\right)^{1/2}
\left(
\sum_{i,j=1}^p
\left(
|\langle u,w_i\rangle|
-
|\langle v,w_i\rangle|
\right)^2
\langle u,w_j\rangle^2
\right)^{1/2}
\\
&\qquad=
\|\Sigma_n\|_F
\left(
\sum_{i=1}^p
\left(
|\langle u,w_i\rangle|
-
|\langle v,w_i\rangle|
\right)^2
\right)^{1/2}
\left(
\sum_{j=1}^p
\langle u,w_j\rangle^2
\right)^{1/2}
\\
&\qquad\leq
\|\Sigma_n\|_F
\left(
\sum_{i=1}^p
|\langle u-v,w_i\rangle|^2
\right)^{1/2}
\\
&\qquad=
\|\Sigma_n\|_F\|u-v\|_2,
\end{align*}
where we used $|
|\langle u,w_i\rangle|
-
|\langle v,w_i\rangle||
\leq
|\langle u-v,w_i\rangle|$.
The second term on the right-hand side of
\eqref{eq:Lipschitz} is bounded in the same way, and hence,
\begin{align}
\label{eq:TV-Lipschitz}
\left|
\|\widehat{\mathsf K}_n^{u}\|_{\mathrm{TV}}
-
\|\widehat{\mathsf K}_n^{v}\|_{\mathrm{TV}}
\right|
\leq
2\|\Sigma_n\|_F\|u-v\|_2.
\end{align}
Since $\|\Sigma_n\|_F^2
=
(p/2)+(p/2)(1+\kappa)^2
\leq
p(1+\kappa)^2$,
combining \eqref{eq:u-unbar-distance} and
\eqref{eq:TV-Lipschitz} yields
\[
\left|
\|\widehat{\mathsf K}_n^{u_n}\|_{\mathrm{TV}}
-
\|\widehat{\mathsf K}_n^{\overline u_n}\|_{\mathrm{TV}}
\right|
=
O_{\mathbb P}(1).
\]                                     
Since $\varepsilon\in(0,1/2)$, we have $p^{1/2-\varepsilon}\longrightarrow\infty$.
Therefore,
\begin{align*}
\mathbb P\left(
\|\widehat{\mathsf K}_n^{\overline u_n}\|_{\mathrm{TV}}
<
\frac{c_3}{2}p^{1/2-\varepsilon}
\right)
&\leq
\mathbb P\left(
\|\widehat{\mathsf K}_n^{u_n}\|_{\mathrm{TV}}
<
c_3p^{1/2-\varepsilon}
\right)
\\
&\qquad+
\mathbb P\left(
\left|
\|\widehat{\mathsf K}_n^{u_n}\|_{\mathrm{TV}}
-
\|\widehat{\mathsf K}_n^{\overline u_n}\|_{\mathrm{TV}}
\right|
>
\frac{c_3}{2}p^{1/2-\varepsilon}
\right)
\longrightarrow0.
\end{align*}
Indeed, the first probability converges to zero by Lemma~\ref{lem:Gaussian-retention}, while the second converges to
zero because its random quantity is $O_{\mathbb P}(1)$ and its
threshold tends to infinity.
Hence, setting $c_4\coloneqq c_3/2>0$,
we conclude that
\begin{align*}
\mathbb P\left(
\|\widehat{\mathsf K}_n^{\overline u_n}\|_{\mathrm{TV}}
\geq
c_4p^{1/2-\varepsilon}
\right)
\longrightarrow1.
\end{align*}
This completes the proof.
\end{proof}

\subsection{Proof of Lemma~\ref{lem:block-spherical-deterministic-transfer}}

\begin{proof}[Proof of Lemma~\ref{lem:block-spherical-deterministic-transfer}]
Let $Q_{1,n}$ and $Q_{2,n}$ be independent Haar-distributed orthogonal
matrices in $\mathbb R^{(p/2)\times(p/2)}$, independent of $Z_n$, and
define
\begin{align} \label{eq:def_Q_n_block}
    Q_n
    \coloneqq
    \begin{pmatrix}
        Q_{1,n} & 0
        \\
        0 & Q_{2,n}
    \end{pmatrix}.
\end{align}
By the block-diagonal structure of $\Sigma_n$, the matrices $Q_n$ and
$\Sigma_n$ commute, and the same is true for $Q_n$ and
$\Sigma_n^{1/2}$. Moreover, conditional on $Q_n$, orthogonal invariance
of the standard Gaussian distribution gives
\[
    Z_nQ_n
    \overset{d}=
    Z_n.
\]
Since $Q_n$ is independent of $Z_n$, we therefore have
\[
    (Z_nQ_n,Q_n)
    \overset{d}=
    (Z_n,Q_n).
\]

Next, consider the transformed data matrix $X_nQ_n=Z_n\Sigma_n^{1/2}Q_n$, whose
sample covariance matrix is given by
\begin{align*}
Q_n^\top
\left(
\frac1n
\Sigma_n^{1/2}Z_n^\top Z_n\Sigma_n^{1/2}
\right)
Q_n
=
Q_n^\top S_nQ_n.
\end{align*}
Using the representation in \eqref{eq:sam-cov-decomposition}, we have
\[
    Q_n^\top S_nQ_n
    =
    \sum_{i=1}^p
    \lambda_i(S_n)
    (Q_n^\top w_i)(Q_n^\top w_i)^\top.
\]
Moreover, since $Q_n$ commutes with $\Sigma_n$,
\[
    \langle Q_n^\top w_i,\Sigma_nQ_n^\top w_j\rangle
    =
    \langle w_i,Q_n\Sigma_nQ_n^\top w_j\rangle
    =
    \langle w_i,\Sigma_nw_j\rangle.
\] 
Consequently, the signed measure obtained from the transformed design
$Z_nQ_n$ and the signal $\beta_n^\star$ is
\begin{align}
\frac{1}{\|\beta_n^\star\|_2^2}
&\sum_{i,j=1}^p
\langle\beta_n^\star,Q_n^\top w_i\rangle
\langle Q_n^\top w_i,\Sigma_nQ_n^\top w_j\rangle
\langle Q_n^\top w_j,\beta_n^\star\rangle
\delta_{(\lambda_i(S_n),\lambda_j(S_n))}
\label{eq:signed_measure1}\\
&=
\frac{1}{\|Q_n\beta_n^\star\|_2^2}
\sum_{i,j=1}^p
\langle Q_n\beta_n^\star,w_i\rangle
\langle w_i,\Sigma_nw_j\rangle
\langle w_j,Q_n\beta_n^\star\rangle
\delta_{(\lambda_i(S_n),\lambda_j(S_n))}
=
\widehat{\mathsf K}_n^{Q_n\beta_n^\star}. \nonumber 
\end{align}
Since $(Z_nQ_n,Q_n)\overset{d}=
(Z_n,Q_n)$, the measure in \eqref{eq:signed_measure1}
has the same distribution as $\widehat{\mathsf K}_n^{\beta_n^\star}$.
Taking total variation norms in \eqref{eq:signed_measure1} therefore gives
\begin{align}
\label{eq:TV-rotation-invariance}
    \|\widehat{\mathsf K}_n^{\beta_n^\star}\|_{\mathrm{TV}}
    \overset{d}=
    \|\widehat{\mathsf K}_n^{Q_n\beta_n^\star}\|_{\mathrm{TV}}.
\end{align}
It remains to identify the distribution of $Q_n\beta_n^\star$. Write
\[
    \beta_n^\star
    =
    \begin{pmatrix}
        \beta_{n,1}^\star
        \\
        \beta_{n,2}^\star
    \end{pmatrix},
    \qquad
    \beta_{n,1}^\star,\beta_{n,2}^\star\in\mathbb R^{p/2}.
\]
By the block-diagonal structure of $Q_n$ in \eqref{eq:def_Q_n_block}, 
\[
    Q_n\beta_n^\star
    =
    \begin{pmatrix}
        Q_{1,n}\beta_{n,1}^\star
        \\[2mm]
        Q_{2,n}\beta_{n,2}^\star
    \end{pmatrix}.
\]
Note that, by construction of $\beta_n^\star$, $\|\beta_{n,1}^\star\|_2^2
    =
    \|\beta_{n,2}^\star\|_2^2
    =
    \frac12$.
Since $Q_{1,n}$ and $Q_{2,n}$ are independent Haar-distributed
orthogonal matrices, the two blocks are independent and uniformly
distributed on the sphere of radius $2^{-1/2}$ in
$\mathbb R^{p/2}$.
On the other hand, since $g_n^{(1)}$ and $g_n^{(2)}$ are independent
standard Gaussian vectors in $\mathbb R^{p/2}$, the vectors
\[
    u_n^{(r)}
    =
    \frac{g_n^{(r)}}{\|g_n^{(r)}\|_2},
    \qquad r=1,2,
\]
are independent and uniformly distributed on the unit sphere in
$\mathbb R^{p/2}$. Hence,
\[
    Q_n\beta_n^\star
    \overset{d}=
    \begin{pmatrix}
        2^{-1/2}u_n^{(1)}
        \\[1mm]
        2^{-1/2}u_n^{(2)}
    \end{pmatrix}
    =
    \overline u_n.
\] 
Since $Q_n$ is independent of $Z_n$ and $\beta_n^\star$ is
deterministic, $Q_n\beta_n^\star$ is independent of $S_n$. Similarly,
$\overline u_n$ is independent of $S_n$, since $g_n$ is independent of
$Z_n$. Therefore,
\[
    \bigl(S_n,Q_n\beta_n^\star\bigr)
    \overset{d}=
    \bigl(S_n,\overline u_n\bigr).
\]
Combining this identity with \eqref{eq:TV-rotation-invariance} yields
\[
    \|\widehat{\mathsf K}_n^{\beta_n^\star}\|_{\mathrm{TV}}
    \overset{d}=
    \|\widehat{\mathsf K}_n^{\overline u_n}\|_{\mathrm{TV}}.
\]
\end{proof}

\section{Technical details for Sections~\ref{sec:proof-bilinear}
and~\ref{sec:proof-overlap-measure}}\label{appendix_C}

\subsection{Lemma~\ref{lem:eventual-spectral-confinement} and its proof}

\begin{lemma}
\label{lem:eventual-spectral-confinement}
Suppose that Assumptions~\ref{ass:GR-basic-1}--\ref{ass:GR-basic-3} hold.
Then there exists a deterministic constant $M<\infty$ such that,
almost surely,
\[
    \operatorname{spec}(S_n)\subseteq[0,M]
\]
for all sufficiently large $n$.
\end{lemma}

\begin{proof}
By Assumption~\ref{ass:GR-basic-3},
\begin{equation}
\label{eq:uniform-Sigma-bound}
    C_\Sigma
    \coloneqq
    \sup_{n\geq1}\|\Sigma_n\|_2
    <\infty.
\end{equation}
It therefore remains to obtain an almost-sure deterministic
upper bound on $n^{-1/2}\|Z_n\|_2$.
We define
\[
    C_\gamma
    \coloneqq
    \sup_{n\geq 1}\gamma_n,
    \qquad
    C_8
    \coloneqq
    \sup_{n\geq1}
    \max_{1\leq i\leq n,\,1\leq j\leq p}
    \mathbb E|z_{ij}|^8
    <\infty.
\]
Set $L_n\coloneqq n^{2/5}$ and define
\[
    z_{ij}^{L_n}
    \coloneqq
    z_{ij}\mathbf1\{|z_{ij}|\leq L_n\}
    -
    \mathbb E\!\left[
        z_{ij}\mathbf1\{|z_{ij}|\leq L_n\}
    \right]
\]
and
\[
    z_{ij}^{L_n+}
    \coloneqq
    z_{ij}\mathbf1\{|z_{ij}|>L_n\}
    -
    \mathbb E\!\left[
        z_{ij}\mathbf1\{|z_{ij}|>L_n\}
    \right].
\]
Since $\mathbb E(z_{ij})=0$, we have
\[
    z_{ij}=z_{ij}^{L_n}+z_{ij}^{L_n+}.
\]
Consequently, writing $Z_n^{L_n}$ and $Z_n^{L_n+}$
for the corresponding $n\times p$ matrices,
\begin{equation}
\label{eq:Z-truncation-decomposition}
    Z_n=Z_n^{L_n}+Z_n^{L_n+}.
\end{equation}

We first control the truncated part $Z_n^{L_n}$.
Its entries are independent and centered. Moreover,
\[
    |z_{ij}^{L_n}|
    \leq
    |z_{ij}|\mathbf1\{|z_{ij}|\leq L_n\}
    +
    \mathbb E\!\left[
        |z_{ij}|\mathbf1\{|z_{ij}|\leq L_n\}
    \right]
    \leq2L_n.
\]
Furthermore,
\[
    \mathbb E|z_{ij}^{L_n}|^2
    =
    \operatorname{Var}\!\left(
        z_{ij}\mathbf1\{|z_{ij}|\leq L_n\}
    \right)
    \leq
    \mathbb E\!\left[
        z_{ij}^2\mathbf1\{|z_{ij}|\leq L_n\}
    \right]
    \leq1.
\]
Using $|a-b|^4\leq8(|a|^4+|b|^4)$ and Jensen's
inequality, we also obtain
\begin{align*}
    \mathbb E|z_{ij}^{L_n}|^4
    &\leq
    8\mathbb E\!\left[
        |z_{ij}|^4\mathbf1\{|z_{ij}|\leq L_n\}
    \right]
    +
    8\left|
        \mathbb E\!\left[
            z_{ij}\mathbf1\{|z_{ij}|\leq L_n\}
        \right]
    \right|^4
    \\
    &\leq
    16\mathbb E|z_{ij}|^4
    \leq
    16\bigl(\mathbb E|z_{ij}|^8\bigr)^{1/2}
    \leq
    16C_8^{1/2}.
\end{align*}
In particular, this bound is independent of
the truncation level $L_n$. Using Lata{\l}a's operator-norm estimate
\cite[Theorem~2]{Latala2005} therefore gives
\begin{align*}
    \mathbb E\|Z_n^{L_n}\|_2
    \leq C\Bigg\{
        \max_{1\leq i\leq n}
        \left(
            \sum_{j=1}^{p}
            \mathbb E|z_{ij}^{L_n}|^2
        \right)^{1/2}
        +
        \max_{1\leq j\leq p}
        \left(
            \sum_{i=1}^{n}
            \mathbb E|z_{ij}^{L_n}|^2
        \right)^{1/2}
        +
        \left(
            \sum_{i=1}^{n}\sum_{j=1}^{p}
            \mathbb E|z_{ij}^{L_n}|^4
        \right)^{1/4}
    \Bigg\},
\end{align*}
where $C$ is a universal constant.
Using the preceding moment bounds,
\[
    \mathbb E\|Z_n^{L_n}\|_2
    \leq
    C\left\{
        \sqrt p+\sqrt n+
        \bigl(16C_8^{1/2}np\bigr)^{1/4}
    \right\}.
\]
Since $p/n\leq C_\gamma$, there exists a deterministic
constant $C_0<\infty$, independent of $n$, such that
\begin{equation}
\label{eq:truncated-expectation-bound}
    \mathbb E\|Z_n^{L_n}\|_2
    \leq C_0\sqrt n,
    \qquad n\geq1.
\end{equation}

We next upgrade this expectation bound to an almost-sure
bound. Since the entries of $Z_n^{L_n}$ are independent
and satisfy $|z_{ij}^{L_n}|\leq2L_n$ almost surely,
\cite[Lemma~28]{GreenRomanov2025} yields, for every $t>0$,
\[
    \mathbb P\left(
        \|Z_n^{L_n}\|_2
        \geq
        \mathbb E\|Z_n^{L_n}\|_2+t
    \right)
    \leq
    \exp\left\{
        -\frac{t^2}{8L_n^2}
    \right\}.
\]
Taking $t=\sqrt n$, using
\eqref{eq:truncated-expectation-bound} and the preceding display, we obtain
\[
    \mathbb P\left(
        \|Z_n^{L_n}\|_2
        \geq
        (C_0+1)\sqrt n
    \right)
    \leq
    \exp\left\{
        -\frac{n}{8L_n^2}
    \right\}
    =
    \exp\left\{
        -\frac18n^{1/5}
    \right\}.
\]
Since
\[
    \sum_{n=1}^{\infty}
    \exp\left\{
        -\frac18n^{1/5}
    \right\}
    <\infty,
\]
the Borel--Cantelli lemma gives, almost surely,
\[
    \|Z_n^{L_n}\|_2
    \leq
    (C_0+1)\sqrt n
\]
for all sufficiently large $n$.
In particular,
\begin{equation}
\label{eq:truncated-as-bound}
    \limsup_{n\to\infty}
    \frac{\|Z_n^{L_n}\|_2}{\sqrt n}
    \leq C_0+1
    \qquad\text{almost surely}.
\end{equation}

We now control the tail matrix $Z_n^{L_n+}$.
By the union bound and Markov's inequality,
\begin{align*}
    \mathbb P\left(
        \max_{\substack{
            1\leq i\leq n\\
            1\leq j\leq p
        }}
        |z_{ij}|>L_n
    \right)
    \leq
    \sum_{i=1}^{n}\sum_{j=1}^{p}
    \mathbb P(|z_{ij}|>L_n)
    \leq
    np\,C_8L_n^{-8}
    \leq
    C_\gamma C_8n^{-6/5}.
\end{align*}
The right-hand side is summable in $n$.
Hence, another application of the Borel--Cantelli lemma
gives an event $\Omega_1$ with $\mathbb P(\Omega_1)=1$
such that, for every $\omega\in\Omega_1$, there exists
$N_1(\omega)\in\mathbb N$ for which
\begin{equation}
\label{eq:Z-eventual-entry-truncation}
    \max_{\substack{
        1\leq i\leq n\\
        1\leq j\leq p
    }}
    |z_{ij}(\omega)|
    \leq L_n,
    \qquad n\geq N_1(\omega).
\end{equation}
Let $\Omega_0$ be the intersection of $\Omega_1$
with the probability-one event on which
\eqref{eq:truncated-as-bound} holds.
Then $\mathbb P(\Omega_0)=1$.
Fix $\omega\in\Omega_0$. For every $n\geq N_1(\omega)$,
\eqref{eq:Z-eventual-entry-truncation} and the definition
of $z_{ij}^{L_n+}$ give
\[
    z_{ij}^{L_n+}(\omega)
    =
    -\mathbb E\!\left[
        z_{ij}\mathbf1\{|z_{ij}|>L_n\}
    \right],
    \qquad
    1\leq i\leq n,\quad1\leq j\leq p.
\]
For every $n\in\mathbb N$ and every index pair $(i,j)$, we have
\begin{align*}
    \left|
        \mathbb E\!\left[
            z_{ij}\mathbf1\{|z_{ij}|>L_n\}
        \right]
    \right|
    \leq
    \mathbb E\!\left[
        |z_{ij}|\mathbf1\{|z_{ij}|>L_n\}
    \right]
    \leq
    L_n^{-7}\mathbb E|z_{ij}|^8
    \leq
    C_8L_n^{-7}.
\end{align*}
Moreover, using Assumption~\ref{ass:GR-basic-1},
there exists a deterministic $n_0\in\mathbb N$ such that
\[
    p\leq2\gamma n,
    \qquad n\geq n_0.
\]
We therefore obtain, for every
$n\geq\max\{N_1(\omega),n_0\}$,
\begin{align*}
    \frac{\|Z_n^{L_n+}(\omega)\|_2}{\sqrt n}
    &\leq
    \frac{\|Z_n^{L_n+}(\omega)\|_F}{\sqrt n}
    \leq
    \frac{\sqrt{np}}{\sqrt n}\,
    C_8L_n^{-7}
    \leq
    C_8(2\gamma)^{1/2}
    n^{1/2-14/5}
    =C_8(2\gamma)^{1/2}n^{-23/10}.
\end{align*}
Consequently,
\begin{equation}
\label{eq:tail-as-bound}
    \frac{\|Z_n^{L_n+}(\omega)\|_2}{\sqrt n}
    \longrightarrow0.
\end{equation}

Combining \eqref{eq:Z-truncation-decomposition},
\eqref{eq:truncated-as-bound}, and
\eqref{eq:tail-as-bound}, we conclude that
\begin{align*}
    \limsup_{n\to\infty}
    \frac{\|Z_n(\omega)\|_2}{\sqrt n}
    \leq
    \limsup_{n\to\infty}
    \frac{\|Z_n^{L_n}(\omega)\|_2}{\sqrt n}
    +
    \limsup_{n\to\infty}
    \frac{\|Z_n^{L_n+}(\omega)\|_2}{\sqrt n}
    \leq C_0+1.
\end{align*}
Thus, choosing the deterministic constant
$C_Z\coloneqq C_0+2$, there exists
$N_2(\omega)\in\mathbb N$ such that
\begin{equation}
\label{eq:Z-eventual-bound}
    \frac{\|Z_n(\omega)\|_2}{\sqrt n}
    \leq C_Z,
    \qquad n\geq N_2(\omega).
\end{equation}

Finally, recall that
\[
    S_n(\omega)
    =
    \frac1nX_n(\omega)^\top X_n(\omega)
    =
    \frac1n
    \Sigma_n^{1/2}
    Z_n(\omega)^\top Z_n(\omega)
    \Sigma_n^{1/2}.
\]
Hence,
\begin{align*}
    \|S_n(\omega)\|_2
    \leq
    \frac1n
    \|\Sigma_n^{1/2}\|_2^2
    \|Z_n(\omega)^\top Z_n(\omega)\|_2
    =
    \|\Sigma_n\|_2
    \frac{\|Z_n(\omega)\|_2^2}{n}.
\end{align*}
Using \eqref{eq:uniform-Sigma-bound} and
\eqref{eq:Z-eventual-bound}, it follows that
\[
    \|S_n(\omega)\|_2
    \leq C_\Sigma C_Z^2,
    \qquad n\geq N_2(\omega).
\]
Since $S_n(\omega)$ is positive semidefinite,
\[
    \operatorname{spec}(S_n(\omega))
    \subseteq[0,\|S_n(\omega)\|_2].
\]
Therefore, choosing the deterministic constant
$M\coloneqq C_\Sigma C_Z^2+1$, we obtain
\[
    \operatorname{spec}(S_n(\omega))
    \subseteq[0,M],
    \qquad n\geq N_2(\omega).
\]
Since $\omega\in\Omega_0$ was arbitrary and
$\mathbb P(\Omega_0)=1$, this proves that, almost surely,
\[
    \operatorname{spec}(S_n)\subseteq[0,M]
\]
for all sufficiently large $n$.
This completes the proof.
\end{proof}

\subsection{Marchenko--Pastur notation and the Silverstein--Choi support characterization}

We briefly recall some notation and results used throughout this section.
Recall the generalized Marchenko--Pastur probability measure $F=F_{\gamma,H}$ and its
companion measure $\overline F=\overline F_{\gamma,H}$ defined in \eqref{eq:GR-MP-weak-convergence}.
Since $\overline F
    =
    (1-\gamma)\delta_0
    +
    \gamma F$,
it follows that, for every Borel set
$A\subseteq\mathbb R\setminus\{0\}$,
\begin{align}
\label{eq:MP-companion-relation-away-zero}
    \overline F(A) = \gamma\,F(A).
\end{align}
Moreover, recall that
$\mathcal S=
    \operatorname{supp}(F)$ and $\mathcal S^+
    =\mathcal S\cap(0,\infty)$.
By \eqref{eq:MP-companion-relation-away-zero},
\begin{align}
\label{eq:nonzero-support-identification}
    \operatorname{supp}(\overline F)
    \cap(0,\infty)
    =\operatorname{supp}( F)
    \cap(0,\infty)=
    \mathcal S^+.
\end{align}
Recall that we denote by $v:\mathbb C^+\to\mathbb C^+$ the Stieltjes transform of
$\overline F$ defined in Lemma~\ref{lem:GR-Stieltjes-convergence}. By
\cite[Equation~(1.3)]{SilversteinChoi1995}, for every
$z\in\mathbb C^+$, $v(z)$ is the unique $w\in\mathbb C^+$ satisfying
\begin{align}
\label{eq:MP-Silverstein-equation}
    z
    =
    -\frac1w
    +
    \gamma
    \int_{[0,\infty)}
    \frac{t}{1+tw}
    \,dH(t).
\end{align}
Following \cite[Equation~(1.4)]{SilversteinChoi1995}, define
\begin{align}
\label{eq:inverse_mp_map}
    z_{\overline F}(w)
    \coloneqq
    -\frac1w
    +
    \gamma
    \int_{[0,\infty)}
    \frac{t}{1+tw}
    \,dH(t),
    \qquad
    w\in
    \left\{
        w\in\mathbb C\setminus\{0\}:
        -\frac1w\notin\operatorname{supp}(H)
    \right\}.
\end{align}
Combining \eqref{eq:MP-Silverstein-equation} and
\eqref{eq:inverse_mp_map} yields
\begin{align}
\label{eq:fixed_point_id}
    z_{\overline F}(v(z))
    =
    z,
    \qquad
    z\in\mathbb C^+.
\end{align}
By \cite[Theorem~1.1]{SilversteinChoi1995}, the limit
\[
    v_0(x)
    \coloneqq
    \lim_{\substack{z\to x\\ z\in\mathbb C^+}}
    v(z)
\]
exists for every $x\in\mathbb R\setminus\{0\}$, and $v_0$ is
continuous on $\mathbb R\setminus\{0\}$. Moreover, by \cite[Theorem~1.1]{SilversteinChoi1995} and \eqref{eq:MP-companion-relation-away-zero},
$F$ has a continuous density $f_{\gamma,H}$ on
$(0,\infty)$ satisfying
\begin{equation}
\label{eq:GR-density-boundary-value}
    \frac1\pi\Im\{v_0(x)\}
    =
    \gamma f_{\gamma,H}(x),
    \qquad
    x>0.
\end{equation}
For the analysis on the real axis, define
\[
    \mathcal B_H
    \coloneqq
    \left\{
        w\in\mathbb R\setminus\{0\}:
        -\frac1w\notin\operatorname{supp}(H)
    \right\},
\]
and
\begin{equation}
\label{eq:MP-real-inverse-map}
    x_{\overline F}(w)
    \coloneqq
    -\frac1w
    +
    \gamma
    \int_{[0,\infty)}
    \frac{t}{1+tw}
    \,dH(t),
    \qquad
    w\in\mathcal B_H.
\end{equation}

\subsection{Lemma~\ref{lem:positive-spectrum-separated-zero} and its proof}

\begin{lemma}
\label{lem:positive-spectrum-separated-zero}
Suppose that Assumptions~\ref{ass:GR-basic-1}, \ref{ass:GR-basic-4} and~\ref{eq:GR-main-positive-population-gap} hold.
Then there exists $\delta_0>0$ such that
\[
    (0,\delta_0)
    \cap
    \mathcal S
    =
    \varnothing.
\]
In particular, the positive part $\mathcal S^+$ of the
Marchenko--Pastur support is separated from zero.
\end{lemma}

\begin{proof}
By \cite[Theorem~4.2]{SilversteinChoi1995}, if
$v\in\mathcal B_H$ and
$\frac{d}{dv}x_{\overline F}(v)>0$, then
\[
    x_{\overline F}(v)
    \notin
    \operatorname{supp}(\overline F).
\]
Thus, whenever in addition $x_{\overline F}(v)>0$,
\eqref{eq:nonzero-support-identification} implies
\[
    x_{\overline F}(v)
    \notin
    \mathcal S.
\]
It therefore suffices to find an interval
$J\subseteq\mathcal B_H$ and some $\delta>0$ such that
$\frac{d}{dv}x_{\overline F}(v)>0$ for $v\in J$ and
\[
    x_{\overline F}(J)=(0,\delta).
\]
We construct such an interval in the two cases
$\gamma[1-H(\{0\})]<1$ and $\gamma[1-H(\{0\})]>1$.

Suppose first that $\gamma[1-H(\{0\})]<1$.
Then, for every $v<-\tau_-^{-1}$, we have
$0<-1/v<\tau_-$ and therefore
\[
    -\frac1v
    \notin
    \operatorname{supp}(H).
\]
Thus,
\[
    \left(
        -\infty,-\frac1{\tau_-}
    \right)
    \subseteq
    \mathcal B_H.
\]
Multiplying the definition of $x_{\overline F}$ in
\eqref{eq:MP-real-inverse-map} by $v<-\tau_-^{-1}$ gives
\begin{equation}
\label{eq:positive-gap-vx}
    v x_{\overline F}(v)
    =
    -1
    +
    \gamma
    \int_{(0,\infty)}
    \frac{tv}{1+tv}
    \,dH(t),
\end{equation}
where the contribution of the atom of $H$ at zero vanishes. For every $v\leq-2\tau_-^{-1}$ and
$t\in\operatorname{supp}(H)\setminus\{0\}$, we have
$tv\leq-2$.
Consequently,
\[
    |1+tv|
    =
    |tv|-1
    \geq
    \frac12|tv|,
\]
and hence
\[
    \left|
        \frac{tv}{1+tv}
    \right|
    \leq2.
\]
For every fixed $t>0$,
\[
    \frac{tv}{1+tv}
    =
    \frac1{1+\frac1{tv}}
    \longrightarrow1
    \qquad
    \text{as }v\to-\infty.
\]
The dominated convergence theorem applied to
\eqref{eq:positive-gap-vx} therefore yields
\[
    v x_{\overline F}(v)
    \longrightarrow
    -1+\gamma[1-H(\{0\})]
    <0
    \qquad
    \text{as }v\to-\infty.
\]
It follows that $x_{\overline F}(v)\to0$ as $v\to-\infty$.
Moreover, since $v<0$, there exists
$v_0<-\tau_-^{-1}$ such that $x_{\overline F}(v)>0$
for $v<v_0$.

We next differentiate $x_{\overline F}$.
To justify differentiation under the integral sign,
fix $v_0<-\tau_-^{-1}$ and choose $\varepsilon>0$ such that
\[
    b
    \coloneqq
    v_0+\varepsilon
    <
    -\tau_-^{-1}.
\]
For $v\in(v_0-\varepsilon,v_0+\varepsilon)$ and
$t\in\operatorname{supp}(H)\setminus\{0\}$,
we have $v\leq b<0$ and $t\geq\tau_-$, and hence
\[
    |1+tv|
    =
    t|v|-1
    \geq
    t|b|-1
    \geq
    t\bigl(|b|-\tau_-^{-1}\bigr).
\]
Consequently,
\[
    \left|
        \frac{t}{1+tv}
    \right|
    \leq
    \frac1{|b|-\tau_-^{-1}}
    \qquad\text{and}\qquad
    \frac{t^2}{|1+tv|^2}
    \leq
    \frac1{(|b|-\tau_-^{-1})^2}.
\]
These constant majorants are $H$-integrable.
Hence, \cite[Theorem~6.28]{Klenke}, applied on the open interval
$(v_0-\varepsilon,v_0+\varepsilon)$, permits differentiation
under the integral sign.
Since $v_0<-\tau_-^{-1}$ was arbitrary,
\begin{equation}
\label{eq:diff_under_int}
    \frac{d}{dv}x_{\overline F}(v)
    =
    \frac1{v^2}
    -
    \gamma
    \int_{(0,\infty)}
    \frac{t^2}{(1+tv)^2}
    \,dH(t),
    \qquad
    v<-\tau_-^{-1}.
\end{equation}
Multiplying \eqref{eq:diff_under_int} by $v^2$ gives
\[
    v^2
    \frac{d}{dv}x_{\overline F}(v)
    =
    1
    -
    \gamma
    \int_{(0,\infty)}
    \left(
        \frac{tv}{1+tv}
    \right)^2
    \,dH(t).
\]
For $v\leq-2/\tau_-$, the integrand is bounded by $4$,
and for every fixed $t>0$ it converges to $1$
as $v\to-\infty$.
Hence another application of dominated convergence yields
\[
    v^2
    \frac{d}{dv}x_{\overline F}(v)
    \longrightarrow
    1-\gamma[1-H(\{0\})]
    >0
    \qquad
    \text{as }v\to-\infty.
\]
Thus there exists $v_-<-2/\tau_-$ such that
\[
    x_{\overline F}(v)>0
    \qquad\text{and}\qquad
    \frac{d}{dv}x_{\overline F}(v)>0
\]
for every $v\leq v_-$.
Set $\delta_-
    \coloneqq
    x_{\overline F}(v_-)
    >0$.
Then $x_{\overline F}$ is continuous and strictly increasing
on $(-\infty,v_-]$, with limit zero as $v\to-\infty$.
It therefore maps $(-\infty,v_-)$ bijectively onto
$(0,\delta_-)$. For every $x\in(0,\delta_-)$, there is thus
$v_x\in\mathcal B_H$ satisfying
\[
    x=x_{\overline F}(v_x)
    \qquad\text{and}\qquad
    \frac{d}{dv}x_{\overline F}(v_x)>0.
\]
Theorem~4.2 of \cite{SilversteinChoi1995} gives
$x\notin\operatorname{supp}(\overline F_{\gamma,H})$.
Since $x>0$, \eqref{eq:nonzero-support-identification}
then gives $x\notin\mathcal S$.
Consequently,
\[
    (0,\delta_-)
    \cap
    \mathcal S
    =
    \varnothing.
\]

Suppose now that $\gamma[1-H(\{0\})]>1$.
For every $v>0$, $-v^{-1}<0$, and therefore $(0,\infty)\subseteq\mathcal B_H$.
Define
\[
    h(v)
    \coloneqq
    v x_{\overline F}(v)
    =
    -1
    +
    \gamma
    \int_{(0,\infty)}
    \frac{tv}{1+tv}
    \,dH(t),
    \qquad
    v>0.
\]
Since $0<(tv)(1+tv)^{-1}<1$ for $t,v>0$,
dominated convergence gives
\[
    \lim_{v\downarrow0}h(v)
    =
    -1
    \qquad\text{and}\qquad
    \lim_{v\to\infty}h(v)
    =
    -1+\gamma[1-H(\{0\})]
    >0.
\]
We again justify differentiation under the integral sign.
Fix $v_0>0$ and choose $0<a<v_0<b$.
For $v\in(a,b)$ and $t>0$, the integrand $tv/(1+tv)$
is bounded by $1$, while, writing $y=tv$,
\[
    \frac{t}{(1+tv)^2}
    =
    \frac1v
    \frac{y}{(1+y)^2}
    \leq
    \frac1{4a},
\]
because $\sup_{y\geq0}y(1+y)^{-2}=1/4$.
Thus \cite[Theorem~6.28]{Klenke} applies on the open
interval $(a,b)$ and gives
\begin{equation*}\label{eq:h_diff}
    \frac{d}{dv}h(v)
    =
    \gamma
    \int_{(0,\infty)}
    \frac{t}{(1+tv)^2}
    \,dH(t)
    >0,
    \qquad
    v>0.
\end{equation*}
The strict inequality follows from
$H((0,\infty))=1-H(\{0\})>0$.
The same local bound and dominated convergence also show that
$\frac{d}{dv}h$ is continuous on $(0,\infty)$.

Thus $h$ is continuous and strictly increasing on $(0,\infty)$.
Since its limits at zero and infinity have opposite signs,
there exists a unique $v_\star>0$ such that $h(v_\star)=0$.
Since $h(v)=v x_{\overline F}(v)$, we have
$x_{\overline F}(v_\star)=0$.
Moreover,
\[
    \frac{d}{dv}x_{\overline F}(v)
    =
    \frac{
        v\frac{d}{dv}h(v)-h(v)
    }{
        v^2
    },
\]
and therefore
\[
    \frac{d}{dv}x_{\overline F}(v_\star)
    =
    \frac{
        \frac{d}{dv}h(v_\star)
    }{
        v_\star
    }
    >0.
\]
Since $h$ and $\frac{d}{dv}h$ are continuous,
$\frac{d}{dv}x_{\overline F}$ is continuous on $(0,\infty)$.
Consequently, there exists $\eta>0$ such that
\[
    \frac{d}{dv}x_{\overline F}(v)>0,
    \qquad
    v\in[v_\star,v_\star+\eta].
\]
Set $\delta_+
    \coloneqq
    x_{\overline F}(v_\star+\eta)
    >0.$ Then $x_{\overline F}$ maps $(v_\star,v_\star+\eta)$
bijectively onto $(0,\delta_+)$.
For every $v\in(v_\star,v_\star+\eta)$,
Theorem~4.2 of \cite{SilversteinChoi1995} gives $x_{\overline F}(v)
    \notin
    \operatorname{supp}(\overline F)$.
Since
\[
    x_{\overline F}\bigl((v_\star,v_\star+\eta)\bigr)
    =
    (0,\delta_+)
    \subset(0,\infty),
\]
\eqref{eq:nonzero-support-identification} therefore implies
\[
    (0,\delta_+)
    \cap
    \mathcal S
    =
    \varnothing.
\]
Since $\gamma[1-H(\{0\})]\neq1$ by
Assumption~\ref{eq:GR-main-positive-population-gap},
     the choice
\[
    \delta_0
    \coloneqq
    \begin{cases}
        \delta_-,
        &\gamma[1-H(\{0\})]<1,\\
        \delta_+,
        &\gamma[1-H(\{0\})]>1,
    \end{cases}
\]
proves the claim.
\end{proof}

\subsection{Lemma~\ref{lem:GR-auxiliary-positive-measure} and its proof}

\begin{lemma}
\label{lem:GR-auxiliary-positive-measure}
Suppose that Assumption~\ref{ass:GR-basic} holds, and let
$M>0$ be the deterministic constant in
Lemma~\ref{lem:eventual-spectral-confinement}. 
Then there exists a unique deterministic finite positive
Borel measure $\mathsf C$ supported on $[0,M]$ such that,
almost surely,
\[
    \mathsf C_n\Rightarrow\mathsf C,
\]
where $\mathsf{C}_n$ is given in \eqref{eq:def_measure_Cn}. 
The Stieltjes transform of $\mathsf{C}$ is given by
\begin{align}
\label{eq:GR-C-transform}
    m_{\mathsf C}(z)
    =
    -\frac1z
    \int_{[0,\infty)}
    \frac{t}{1+t v(z)}
    \,dG(t),
    \qquad z\in\mathbb C^+,
\end{align}
where $G$ is the weak limit of $G_n$ from
Assumption~\ref{ass:GR-basic-5}.
\end{lemma}

\begin{proof}
By the spectral decomposition of $S_n$, the Stieltjes transform of $\mathsf C_n$ can be written as
\[
    m_{\mathsf C_n}(z)
    =
    \frac1{\|\beta_n^\star\|_2^2}
    (\beta_n^\star)^\top
    \Sigma_n^{1/2}
    R(z)
    \Sigma_n^{1/2}
    \beta_n^\star,
    \qquad
    R(z)=(S_n-zI_p)^{-1}.
\]
Moreover, by Parseval's identity and the definition of $G_n$,
\begin{align}
\label{eq:GR-Cn-mass}
    \mathsf C_n([0,\infty))
    =
    \frac{
        (\beta_n^\star)^\top
        \Sigma_n
        \beta_n^\star
    }{
        \|\beta_n^\star\|_2^2
    }
    =
    \int_{[0,\infty)}t\,dG_n(t)
    \eqqcolon a_n.
\end{align}
By Assumption~\ref{ass:GR-basic-3},
\[
    C_\Sigma
    \coloneqq
    \sup_{n\geq1}\|\Sigma_n\|_2
    <\infty.
\]
Since $\Sigma_n$ is nonnegative definite and
$\operatorname{supp}(G_n)\subseteq\operatorname{spec}(\Sigma_n)$,
we have
\[
    \operatorname{supp}(G_n)
    \subseteq[0,C_\Sigma],
    \qquad n\geq1.
\]
Assumption~\ref{ass:GR-basic-5} gives $G_n\Rightarrow G$,
and hence $G$ is also supported on $[0,C_\Sigma]$.
Since $t\mapsto t$ is continuous and bounded on this interval,
\[
    a_n
    =
    \int_{[0,\infty)}t\,dG_n(t)
    \longrightarrow
    \int_{[0,\infty)}t\,dG(t)
    \eqqcolon a.
\]

Suppose first that $a>0$.
Then there exists a deterministic $n_1\in\mathbb N$ such that
$a_n>0$ for every $n\geq n_1$.
Set $\widetilde\beta_n
    \coloneqq
    \Sigma_n^{1/2}\beta_n^\star$.
For $n\geq n_1$, we have
\[
    \|\widetilde\beta_n\|_2^2
    =
    a_n\|\beta_n^\star\|_2^2
    >0.
\]
Using the representation from Assumption~\ref{ass:GR-basic-5}, the signal-weighted
population spectral measure corresponding to
$\widetilde\beta_n$ is given by
\begin{align*}\label{eq:GR-Gtilde}
    \widetilde G_n
    =
    \frac1{\|\widetilde\beta_n\|_2^2}
    \sum_{j=1}^p
    |\langle\widetilde\beta_n,q_j\rangle|^2
    \delta_{\lambda_j(\Sigma_n)}
    =
    \frac1{a_n}
    \sum_{j=1}^p
    \lambda_j(\Sigma_n)
    \frac{
        |\langle\beta_n^\star,q_j\rangle|^2
    }{
        \|\beta_n^\star\|_2^2
    }
    \delta_{\lambda_j(\Sigma_n)}.
\end{align*}
Hence, $\widetilde G_n\ll G_n$, with
\[
    d\widetilde G_n
    =
    \frac{t}{a_n}dG_n(t),
    \qquad G_n\text{-almost everywhere}.
\]
Therefore, for every bounded continuous function
$\varphi:[0,\infty)\to\mathbb C$,
\begin{align*}
    \int_{[0,\infty)}
    \varphi(t)\,d\widetilde G_n(t)
    =
    \frac1{a_n}
    \int_{[0,\infty)}
    t\varphi(t)\,dG_n(t)
    \longrightarrow
    \frac1a
    \int_{[0,\infty)}
    t\varphi(t)\,dG(t),
\end{align*}
where we used that $t\mapsto t\varphi(t)$ is continuous
and bounded on the common compact support $[0,C_\Sigma]$.
Consequently,
\[
    \widetilde G_n\Rightarrow\widetilde G,
    \qquad
    d\widetilde G(t)=\frac{t}{a}\,dG(t).
\]
For $n\geq n_1$, define 
\[
    \Theta_n
    \coloneqq
    \frac{
        \widetilde\beta_n\widetilde\beta_n^\top
    }{
        \|\widetilde\beta_n\|_2^2
    },
\]
and set $\Theta_n=0$ for $n<n_1$.
For $n\geq n_1$, $\|\Theta_n\|_{\mathrm{tr}}=1$ and in particular,
$\sup_{n\geq1}\|\Theta_n\|_{\mathrm{tr}}\leq1$.
By \cite[Proposition~2.1 and Remark~2.2]{DoernemannStark},
applied with $\Theta_n$ under
Assumptions~\ref{ass:GR-basic-1}--\ref{ass:GR-basic-4}, we obtain, almost surely,
for every $z\in\mathbb C^+$,
\[
    \operatorname{tr}\!\left[
        \Theta_n
        \left\{
            R(z)
            -
            \bigl(-zv(z)\Sigma_n-zI_p\bigr)^{-1}
        \right\}
    \right]
    \longrightarrow0.
\]
Since
\[
    \bigl(-zv(z)\Sigma_n-zI_p\bigr)^{-1}
    =
    -\frac1z
    \bigl(I_p+v(z)\Sigma_n\bigr)^{-1},
\]
the spectral decomposition of $\Sigma_n$ gives,
for $n\geq n_1$,
\begin{align*}
    &
    \operatorname{tr}\!\left[
        \Theta_n
        \left\{
            R(z)
            -
            \bigl(-zv(z)\Sigma_n-zI_p\bigr)^{-1}
        \right\}
    \right]
    \\
    &\qquad=
    \frac{
        \widetilde\beta_n^\top
        R(z)
        \widetilde\beta_n
    }{
        \|\widetilde\beta_n\|_2^2
    }
    +
    \frac1z
    \int_{[0,\infty)}
    \frac1{1+t v(z)}
    \,d\widetilde G_n(t)
    \\
    &\qquad=
    \frac1{a_n}m_{\mathsf C_n}(z)
    +
    \frac1z
    \int_{[0,\infty)}
    \frac1{1+t v(z)}
    \,d\widetilde G_n(t).
\end{align*}
For fixed $z\in\mathbb C^+$, we have $\Im\{v(z)\}>0$.
Thus $1+t v(z)\neq0$ for every $t\geq0$, and
$t\mapsto(1+t v(z))^{-1}$ is continuous and bounded
on $[0,C_\Sigma]$.
The weak convergence $\widetilde G_n\Rightarrow\widetilde G$
therefore yields
\[
    \int_{[0,\infty)}
    \frac1{1+t v(z)}
    \,d\widetilde G_n(t)
    \longrightarrow
    \int_{[0,\infty)}
    \frac1{1+t v(z)}
    \,d\widetilde G(t).
\]
Combining the preceding displays, we obtain, almost surely, for every $z\in\mathbb C^+$,
\begin{align*}
    \frac1{a_n}m_{\mathsf C_n}(z)
    \longrightarrow
    -\frac1z
    \int_{[0,\infty)}
    \frac1{1+t v(z)}
    \,d\widetilde G(t)
    =
    -\frac1{az}
    \int_{[0,\infty)}
    \frac{t}{1+t v(z)}
    \,dG(t).
\end{align*}
It follows that, almost surely, for every $z\in\mathbb C^+$,
\begin{equation}
\label{eq:GR-Cn-transform-convergence}
    m_{\mathsf C_n}(z)
    \longrightarrow
    -\frac1z
    \int_{[0,\infty)}
    \frac{t}{1+t v(z)}
    \,dG(t)\eqqcolon h_{\mathsf C}(z).
\end{equation}

If $a=0$, then $a_n\to0$.
Since $\mathsf C_n$ is a positive measure, almost surely, for every $z\in \mathbb  C^+$,
\[
    |m_{\mathsf C_n}(z)|
    \leq
    \frac{\mathsf C_n([0,\infty))}{\Im(z)}
    =
    \frac{a_n}{\Im(z)}
    \longrightarrow0.
\]
Since $G$ is supported on $[0,\infty)$ and
$\int t\,dG(t)=0$, we have $G=\delta_0$.
Hence, the right-hand side of
\eqref{eq:GR-C-transform} is zero as well.
Thus \eqref{eq:GR-Cn-transform-convergence}
also holds when $a=0$.

By Lemma~\ref{lem:eventual-spectral-confinement},
there exists an event $\Omega_{\mathrm{0}}$
with $\mathbb P(\Omega_{\mathrm{0}})=1$
such that, for every $\omega\in\Omega_{\mathrm{0}}$,
there exists $N(\omega)\in\mathbb N$ satisfying
\[
    \operatorname{spec}(S_n(\omega))
    \subseteq[0,M],
    \qquad n\geq N(\omega).
\]
Since $\mathsf C_n$ is supported on
$\operatorname{spec}(S_n)$, it follows that
\[
    \operatorname{supp}(\mathsf C_n(\omega))
    \subseteq[0,M],
    \qquad n\geq N(\omega).
\]
Moreover, positivity of $\mathsf C_n$ and
\eqref{eq:GR-Cn-mass} give
\[
    \|\mathsf C_n\|_{\mathrm{TV}}
    =
    \mathsf C_n([0,\infty))
    =
    a_n
    \leq
    \|\Sigma_n\|_2
    \leq C_\Sigma
\]
for every $n$ and every realization.

Let $\Omega_1$ denote the probability-one event on which
\eqref{eq:GR-Cn-transform-convergence} holds simultaneously
for every $z\in\mathbb C^+$, and define $\Omega_2\coloneqq
    \Omega_1\cap\Omega_{\mathrm{0}}$.
Then $\mathbb P(\Omega_2)=1$.
Fix $\omega\in\Omega_2$.
The sequence $\{\mathsf C_n(\omega)\}_{n\geq1}$
satisfies
\[
    \sup_{n\geq1}
    \|\mathsf C_n(\omega)\|_{\mathrm{TV}}
    \leq C_\Sigma,
\]
and
\[
    \operatorname{supp}(\mathsf C_n(\omega))
    \subseteq[0,M],
    \qquad n\geq N(\omega).
\]
Moreover, its Stieltjes transforms converge to
$h_{\mathsf C}(z)$ for every $z\in\mathbb C^+$.
Proposition~\ref{cor:weak-convergence-compact-support}
therefore yields a unique finite signed Borel measure
$\mathsf C_\omega$ supported on $[0,M]$ such that
\[
    \mathsf C_n(\omega)\Rightarrow\mathsf C_\omega,
    \qquad
   m_{\mathsf C_\omega}(z)=h_{\mathsf C}(z),
    \quad z\in\mathbb C^+.
\]
Since $h_{\mathsf C}$ is deterministic and the
Stieltjes transform uniquely determines a finite
signed Borel measure supported on $[0,M]$,
$\mathsf C_\omega$ does not depend on $\omega$.
Denote this common measure by $\mathsf C$.
Since every $\mathsf C_n(\omega)$ is positive,
its weak limit $\mathsf C$ is positive as well.

Finally, $\omega\in\Omega_2$ was arbitrary and
$\mathbb P(\Omega_2)=1$.
Thus, almost surely,
\[
    \mathsf C_n\Rightarrow\mathsf C,
\]
and $m_{\mathsf C}=h_{\mathsf C}$ gives
\eqref{eq:GR-C-transform}.
This completes the proof.
\end{proof}

\subsection{Lemma~\ref{lem:GR-boundary-transform-nonzero} and its proof}

\begin{lemma}
\label{lem:GR-boundary-transform-nonzero}
Suppose that Assumptions~\ref{ass:GR-basic-1},
\ref{ass:GR-basic-3}, and \ref{ass:GR-basic-4} hold. Then
\[
    v_0(s)\neq0,
    \qquad s>0.
\]
Moreover, $\mathcal S$ is compact and
\[
    \mathcal S\cap(0,\infty)\neq\varnothing.
\]
\end{lemma}

\begin{proof}
By Assumptions~\ref{ass:GR-basic-3} and
\ref{ass:GR-basic-4},
\[
    \operatorname{supp}(H)\subseteq[0,C_\Sigma],
    \qquad
    C_\Sigma
    \coloneqq
    \sup_{n\geq1}\|\Sigma_n\|_2
    <\infty.
\]
Indeed, every $H_n$ is supported on $[0,C_\Sigma]$,
and $H_n\Rightarrow H$. Thus $\operatorname{supp}(H)$
is compact. Since Assumption~\ref{ass:GR-basic-4} gives
$H((0,\infty))>0$, we may define
\[
    \tau_+
    \coloneqq
    \max\operatorname{supp}(H)
    \in(0,\infty).
\]

We start by proving the first assertion of the lemma. Suppose, to the contrary, that $v_0(s)=0$ for some
$s\in(0,\infty)$.
Choose a sequence $z_k\in\mathbb C^+$ such that $z_k\to s$.
By the definition of $v_0(s)$,
\[
    v(z_k)\longrightarrow0.
\]
Thus, for all sufficiently large $k$,
$|v(z_k)|\leq(2\tau_+)^{-1}$.
For such $k$ and every $t\in\operatorname{supp}(H)$,
\[
    |1+t v(z_k)|
    \geq
    1-t|v(z_k)|
    \geq
    1-\tau_+|v(z_k)|
    \geq
    \frac12.
\]
Since $H$ is a probability measure,
\[
    \left|
        \int_{[0,\infty)}
        \frac{t}{1+t v(z_k)}
        \,dH(t)
    \right|
    \leq
    \int_{[0,\infty)}
    \frac{t}{|1+t v(z_k)|}
    \,dH(t)
    \leq
    2\int_{[0,\infty)}t\,dH(t)
    \leq
    2\tau_+.
\]
On the other hand, using \eqref{eq:inverse_mp_map},
the reverse triangle inequality, and the preceding display,
we obtain
\[
    |z_k|
    \geq
    \frac1{|v(z_k)|}
    -
    \gamma
    \left|
        \int_{[0,\infty)}
        \frac{t}{1+t v(z_k)}
        \,dH(t)
    \right|
    \geq
    \frac1{|v(z_k)|}
    -
    2\gamma\tau_+
    \longrightarrow\infty,
\]
as $k\to\infty$, contradicting $z_k\to s\in\mathbb R$.
Therefore, $v_0(s)\neq0$ for every $s>0$.

Next, we prove that $\mathcal S$ is compact. Define
\[
    \Phi(r)
    \coloneqq
    x_{\overline F}\left(-\frac1r\right)
    =
    r
    +
    \gamma
    \int_{[0,\infty)}
    \frac{rt}{r-t}
    \,dH(t),
    \qquad r>\tau_+.
\]
Fix $r_0>\tau_+$ and choose
$0<\eta<r_0-\tau_+$.
For $r\in(r_0-\eta,r_0+\eta)$ and
$t\in\operatorname{supp}(H)$,
\[
    r-t>r_0-\eta-\tau_+>0.
\]
Moreover,
\[
    \left|
        \frac{rt}{r-t}
    \right|
    \leq
    \frac{(r_0+\eta)\tau_+}
         {r_0-\eta-\tau_+}
\]
and
\[
    \left|
        \frac{d}{dr}\frac{rt}{r-t}
    \right|
    =
    \frac{t^2}{(r-t)^2}
    \leq
    \frac{\tau_+^2}
         {(r_0-\eta-\tau_+)^2}.
\]
The preceding bounds provide $H$-integrable majorants for the
integrand in $\Phi$ and its derivative, uniformly for
$r\in(r_0-\eta,r_0+\eta)$.
Hence, by \cite[Theorem~6.28]{Klenke}, differentiation under the
integral sign is permitted and yields
\[
    \frac{d}{dr}\Phi(r)
    =
    1
    -
    \gamma
    \int_{[0,\infty)}
    \frac{t^2}{(r-t)^2}
    \,dH(t),
    \qquad
    r\in(r_0-\eta,r_0+\eta).
\]
Since $r_0>\tau_+$ was arbitrary, this identity holds
for every $r>\tau_+$. Choose $R>\tau_+$ such that
\[
    \gamma
    \frac{\tau_+^2}{(R-\tau_+)^2}
    <1.
\]
For $r\geq R$ and $t\in\operatorname{supp}(H)$,
we have $0\leq t\leq\tau_+$ and
$r-t\geq r-\tau_+$. Since $H$ is a probability measure,
\[
    \gamma
    \int_{[0,\infty)}
    \frac{t^2}{(r-t)^2}
    \,dH(t)
    \leq
    \gamma
    \frac{\tau_+^2}{(r-\tau_+)^2}
    \leq
    \gamma
    \frac{\tau_+^2}{(R-\tau_+)^2}
    <1.
\]
Together with the preceding expression for
$\frac{d}{dr}\Phi(r)$, this yields
\[
    \frac{d}{dr}\Phi(r)>0,
    \qquad r\geq R.
\]
Furthermore, since $r>\tau_+=\max\operatorname{supp}(H)$,
we have $r\notin\operatorname{supp}(H)$,
and hence $-1/r\in\mathcal B_H$.
By the chain rule,
\[
    \frac{d}{dr}\Phi(r)
    =
    \frac1{r^2}
    \left.
    \frac{d}{dw}x_{\overline F}(w)
    \right|_{w=-1/r},
\]
and therefore
\[
    \left.
    \frac{d}{dw}x_{\overline F}(w)
    \right|_{w=-1/r}
    =
    r^2\frac{d}{dr}\Phi(r)
    >0,
    \qquad r\geq R.
\]
Thus, Theorem~4.2 of \cite{SilversteinChoi1995} gives
\[
    \Phi(r)
    =
    x_{\overline F}\left(-\frac1r\right)
    \notin
    \operatorname{supp}(\overline F),
    \qquad r\geq R.
\]
Moreover, for every $r>\tau_+$ and
$t\in\operatorname{supp}(H)$, we have
\[
    \frac{rt}{r-t}\geq0.
\]
Consequently,
\[
    \Phi(r)
    =
    r
    +
    \gamma
    \int_{[0,\infty)}
    \frac{rt}{r-t}
    \,dH(t)
    \geq r>0.
\]
Hence, by
\eqref{eq:nonzero-support-identification}, $\Phi(r)\notin\mathcal S$, for $r\geq R$.
Since $\Phi$ is differentiable, it is continuous.
Moreover, it is strictly increasing on $[R,\infty)$
and $\Phi(r)\longrightarrow\infty$ as $r\to \infty$.
It follows that $\Phi\bigl([R,\infty)\bigr)
    =
    [\Phi(R),\infty)$.
Together with $\Phi(r)\notin\mathcal S$ for $r\geq R$, this gives $[\Phi(R),\infty)\cap \mathcal S = \varnothing$. 
Since $F$ is supported on $[0,\infty)$, we have $\mathcal S\subseteq[0,\infty)$, and therefore, $\mathcal S \subseteq [0,\Phi(R))$.
Since the support of a Borel measure is closed,
$\mathcal S$ is closed and bounded, and hence compact.

It remains to show that $\mathcal S \cap (0,\infty)\neq \varnothing$.
The discussion following
\cite[Equation~(1.4)]{SilversteinChoi1995} gives
\[
    \overline F(\{0\})
    =
    \max\left\{
        1-\gamma\bigl(1-H(\{0\})\bigr),
        0
    \right\}.
\]
Since $\overline F
    =
    (1-\gamma)\delta_0+\gamma F$,
we obtain
\[
    F(\{0\})
    =
    \max\left\{
        H(\{0\}),
        1-\frac1\gamma
    \right\}.
\]
Using Assumption~\ref{ass:GR-basic-4}, we have $H(\{0\})<1$.
Moreover, since $\gamma\in(0,\infty)$, $1-\gamma^{-1}<1$.
Consequently, $F(\{0\})<1$ and thus, $F((0,\infty))>0$,
which implies $
    \mathcal S\cap(0,\infty)
    \neq
    \varnothing$, completing the proof.
\end{proof}

\subsection{Lemma~\ref{lem:GR-population-gap-holder} and its proof}

\begin{lemma}
\label{lem:GR-population-gap-holder}
Suppose that
Assumptions~\ref{ass:GR-basic-1},
\ref{ass:GR-basic-3}, and
\ref{ass:GR-basic-4} hold, together with
Assumptions~\ref{eq:GR-main-positive-population-gap} and
\ref{eq:GR-main-zero-density-regularity}.
Set $\mathcal S^+
    \coloneqq
    \mathcal S\cap(0,\infty)$.
Then there exists $C_v<\infty$ such that
\begin{equation}
\label{eq:GR-global-one-third-holder}
    |v_0(s)-v_0(t)|
    \leq
    C_v|s-t|^{1/3},
    \qquad
    s,t\in\mathcal S^+.
\end{equation}
Moreover, if $s\in\mathcal S^+$ and
$f_{\gamma,H}(s)=0$, then
\begin{align}
\label{eq:GR-zero-density-critical-properties}
    v_0(s)\in\mathcal B_H,
    \qquad
    \left.
    \frac{d}{dw}x_{\overline F}(w)
    \right|_{w=v_0(s)}
    =
    0.
\end{align}
\end{lemma}

\paragraph*{Proof outline of Lemma~\ref{lem:GR-population-gap-holder}}

The following proof is divided into six steps.

\medskip
\noindent
\textbf{Step 1:}\label{step1:lem16}
We fix $s_0\in\mathcal S^+$ with
$f_{\gamma,H}(s_0)>0$. We then show that, for all
$s,t\in\mathcal S^+$ sufficiently close to $s_0$,
\begin{align}
\label{eq:Lip_continuous1}
    |v_0(s)-v_0(t)|
    =
    \left|
        \int_t^s
        \frac{d}{dx}v_0(x)\,dx
    \right|
    \leq
    C_{s_0}|s-t|.
\end{align}

\medskip
\noindent
\textbf{Step 2:}\label{step2:lem16}
Next, we fix $s_0\in\mathcal S^+$ with $f_{\gamma,H}(s_0)=0$. Then $s_0$ is either a boundary point
of $\mathcal S^+$ or an interior zero-density point,
that is, a cusp.
Using Assumption~\ref{eq:GR-main-zero-density-regularity},
we first establish that
\begin{align}
\label{eq:step1_in_sec_claim}
    v_0(s_0)\in\mathcal B_H,
    \qquad
    x_{\overline F}(v_0(s_0))
    =
    s_0,
\end{align}
and
\begin{align}
\label{eq:step2_in_sec_claim}
    \left.
    \frac{d}{dw}x_{\overline F}(w)
    \right|_{w=v_0(s_0)}
    =
    0.
\end{align}
Combining \eqref{eq:step1_in_sec_claim} and
\eqref{eq:step2_in_sec_claim} yields the second assertion
of the lemma displayed in
\eqref{eq:GR-zero-density-critical-properties}.

\medskip
\noindent
\textbf{Step 3:}\label{step3:lem16}
Fix $s_0\in\mathcal S^+$ with $f_{\gamma,H}(s_0)=0$.
We show that $r\mapsto z_{\overline F}(-1/r)-s_0$ has a zero
of order $\nu=2$ at $-v_0(s_0)^{-1}$ if its second derivative
at $-v_0(s_0)^{-1}$ is nonzero, and of order $\nu=3$
if this second derivative vanishes.
Consequently,
\[
    z_{\overline F}(w)-s_0
    =
    \{w-v_0(s_0)\}^{\nu}g(w),
    \qquad
    \nu\in\{2,3\},
\]
where $g$ is holomorphic near $v_0(s_0)$ and
$g(v_0(s_0))\neq0$.

\medskip
\noindent
\textbf{Step 4:}\label{step4:lem16}
Using the preceding step for $s_0\in\mathcal S^+$ with
$f_{\gamma,H}(s_0)=0$, we show that, for 
$\nu\in\{2,3\}$ determined there, $s_0$ is the unique point
in a neighborhood of $s_0$ within $\mathcal S^+$ satisfying
$f_{\gamma,H}(s)=0$, and
\begin{align}
\label{eq:GR-holder-outline-pointwise}
    |v_0(s)-v_0(s_0)|
    \leq
    C|s-s_0|^{1/\nu}
\end{align}
for $s\in\mathcal S^+$ sufficiently close to $s_0$.
At every such $s\neq s_0$, we have $f_{\gamma,H}(s)>0$, and
\begin{align}
\label{eq:GR-holder-outline-derivative}
    \left|
        \frac{d}{ds}v_0(s)
    \right|
    \leq
    C|s-s_0|^{-1+1/\nu}.
\end{align}

\medskip
\noindent
\textbf{Step 5:}\label{step5:lem16}
Consider again $s_0\in\mathcal S^+$ with
$f_{\gamma,H}(s_0)=0$. We show that, on either side of $s_0$,
the intersection of $\mathcal S^+$ with a sufficiently small
neighborhood of $s_0$ is connected whenever it is nonempty.
Using this fact, together with
\eqref{eq:GR-holder-outline-pointwise} and
\eqref{eq:GR-holder-outline-derivative}, we then obtain
\begin{align}
\label{eq:local_1_3_hölder}
    |v_0(s)-v_0(t)|
    \leq
    C|s-t|^{1/3}
\end{align}
for all $s,t\in\mathcal S^+$ sufficiently close to $s_0$.
Combining this with \eqref{eq:Lip_continuous1}, it follows
that $v_0$ is $1/3$-H\"older in a neighborhood of every
$s\in\mathcal S^+$.

\medskip
\noindent
\textbf{Step 6:}\label{step6:lem16}
Since $\mathcal S^+$ is compact by
Lemmas~\ref{lem:positive-spectrum-separated-zero} and
\ref{lem:GR-boundary-transform-nonzero}, we use
\eqref{eq:local_1_3_hölder} to establish
\eqref{eq:GR-global-one-third-holder}.

\begin{proof}
\textbf{\hyperref[step1:lem16]{Proof of Step~1}}
We first consider $s_0\in\mathcal S^+$ satisfying
$f_{\gamma,H}(s_0)>0$. By
\eqref{eq:GR-density-boundary-value},
\[
    \Im\{v_0(s_0)\}
    =
    \pi\gamma f_{\gamma,H}(s_0)
    >0.
\]
Hence, by the argument in
\cite[Proof of Theorem~1.1]{SilversteinChoi1995},
given in the final paragraph of Section~3,
\begin{equation}
\label{eq:GR-positive-density-inverse-identity}
    z_{\overline F}(v_0(s_0))
    =
    s_0,
\end{equation}
and
\begin{equation}
\label{eq:GR-positive-density-nonzero-derivative}
    \left.
    \frac{d}{dw}
    z_{\overline F}(w)
    \right|_{w=v_0(s_0)}
    =
    \frac1{v_0(s_0)^2}
    -
    \gamma
    \int_{[0,\infty)}
    \frac{t^2}{\{1+t v_0(s_0)\}^2}
    \,dH(t)
    \neq0.
\end{equation}
Moreover, by the holomorphic inverse function theorem, the companion Stieltjes transform $v$ admits
a holomorphic extension to a complex neighborhood $U$
of $s_0$, which we may choose so that
$U\cap\mathbb R\subset(0,\infty)$.
Let $\widetilde v$ denote this holomorphic extension.
Since $\widetilde v$ agrees with $v$ on $U\cap\mathbb C^+$,
the continuity of $\widetilde v$ and the definition of $v_0$
imply that the restriction of $\widetilde v$ to
$U\cap\mathbb R$ agrees with $v_0$.
Choose $\varepsilon>0$ such that
$[s_0-\varepsilon,s_0+\varepsilon]\subset U$.
Since $\widetilde v$ is holomorphic and agrees with $v_0$
on $U\cap\mathbb R$, the function $v_0$ is continuously
differentiable on $U\cap\mathbb R$.
Hence,
\[
    C_{s_0}
    \coloneqq
    \sup_{x\in[s_0-\varepsilon,s_0+\varepsilon]}
    \left|
        \frac{d}{dx}v_0(x)
    \right|
    <\infty.
\]
Therefore, for
$s,t\in\mathcal S^+\cap(s_0-\varepsilon,s_0+\varepsilon)$,
\begin{align}
\label{eq:GR-positive-density-local-Lipschitz}
    |v_0(s)-v_0(t)|
    =
    \left|
        \int_t^s
        \frac{d}{dx}v_0(x)\,dx
    \right|
    \leq
    C_{s_0}|s-t|.
\end{align}

\textbf{\hyperref[step2:lem16]{Proof of Step~2}}
We next consider a point $s_0\in\mathcal S^+$ satisfying
$f_{\gamma,H}(s_0)=0$.
By \eqref{eq:GR-density-boundary-value},
$\Im\{v_0(s_0)\}=0$.
Moreover, Lemma~\ref{lem:GR-boundary-transform-nonzero}
gives $v_0(s_0)\neq0$, and together with
Assumption~\ref{eq:GR-main-zero-density-regularity},
we obtain
\[
    v_0(s_0)\in\mathcal B_H.
\]
Since $\operatorname{supp}(H)$ is compact by
Assumptions~\ref{ass:GR-basic-3} and \ref{ass:GR-basic-4},
there exist a neighborhood $U\subseteq \mathbb C$ of $v_0(s_0)$
and a constant $c>0$ such that
\begin{align}
\label{eq:local_den_bound}
    |w|\geq c
    \qquad\text{and}\qquad
    |1+tw|\geq c
\end{align}
for every $w\in U$ and every $t\in\operatorname{supp}(H)$. 
Moreover, since $H$ has compact support, \eqref{eq:local_den_bound}
justifies differentiation under the integral sign and
shows that $z_{\overline F}$ is holomorphic on $U$.
Since $v(z)\longrightarrow v_0(s_0)$ as $z\to s_0$,
with $z\in\mathbb C^+$, by
\cite[Theorem~1.1]{SilversteinChoi1995},
we have $v(z)\in U$ whenever $z\in\mathbb C^+$
is sufficiently close to $s_0$.
Letting $z\to s_0$ in \eqref{eq:fixed_point_id}
and using the continuity of $z_{\overline F}$ therefore gives
\[
    z_{\overline F}(v_0(s_0))
    =
    s_0.
\]
Since $v_0(s_0)\in\mathcal B_H$, it follows from
\eqref{eq:MP-real-inverse-map} that
\begin{equation}
\label{eq:GR-holder-zero-density-image}
    x_{\overline F}(v_0(s_0))
    =
    s_0.
\end{equation}

We next prove that
\begin{equation}
\label{eq:GR-holder-zero-density-critical}
    \left.
    \frac{d}{dw}x_{\overline F}(w)
    \right|_{w=v_0(s_0)}
    =
    0.
\end{equation}
By \eqref{eq:local_den_bound} and
\cite[Theorem~6.28]{Klenke}, differentiation under
the integral sign in \eqref{eq:MP-real-inverse-map} gives
\[
    \left.
    \frac{d}{dw}x_{\overline F}(w)
    \right|_{w=v_0(s_0)}
    =
    \frac1{v_0(s_0)^2}
    -
    \gamma
    \int_{[0,\infty)}
    \frac{t^2}
         {\{1+t v_0(s_0)\}^2}
    \,dH(t).
\]
For $z\in\mathbb C^+$, combining
\eqref{eq:fixed_point_id} with \eqref{eq:inverse_mp_map}
and taking imaginary parts, using the identity
$\Im\{w^{-1}\}=-\Im\{w\}|w|^{-2}$ for
$w\in\mathbb C\setminus\{0\}$, yields
\[
    \Im\{z\}
    =
    \Im\{v(z)\}
    \left\{
        \frac1{|v(z)|^2}
        -
        \gamma
        \int_{[0,\infty)}
        \frac{t^2}{|1+t v(z)|^2}
        \,dH(t)
    \right\}.
\]
Since $\Im\{z\}>0$ and $\Im\{v(z)\}>0$,
the expression in braces is strictly positive.
Letting $z\to s_0$, with $z\in\mathbb C^+$,
using \eqref{eq:local_den_bound} and dominated convergence,
we conclude that
\[
    \left.
    \frac{d}{dw}x_{\overline F}(w)
    \right|_{w=v_0(s_0)}
    \geq0.
\]
If this derivative were strictly positive, then,
since $v_0(s_0)\in\mathcal B_H$,
\cite[Theorem~4.2]{SilversteinChoi1995} together with
\eqref{eq:GR-holder-zero-density-image} would imply that
\[
    s_0\notin\operatorname{supp}(\overline F).
\]
Since $s_0>0$, the support identification in \eqref{eq:nonzero-support-identification} would then give
$s_0\notin\mathcal S^+$, contradicting $s_0\in\mathcal S^+$.
Hence \eqref{eq:GR-holder-zero-density-critical} holds.
Since $s_0$ was an arbitrary zero-density point,
this proves \eqref{eq:GR-zero-density-critical-properties}.

\textbf{\hyperref[step3:lem16]{Proof of Step~3}}
In this step, we establish the local factorization at each
point $s_0\in\mathcal S^+$ where $f_{\gamma,H}(s_0)=0$.
Recall that for such points $v_0(s_0)\in\mathcal B_H$.
Note that $r_0\coloneqq-v_0(s_0)^{-1}$ satisfies $r_0\neq 0$, and $r_0\notin\operatorname{supp}(H)$.
Since $\operatorname{supp}(H)$ is compact, there exist
a bounded complex neighborhood $U$ of $r_0$,
with $0\notin U$, and a constant $c>0$ such that
$|r-t|\geq c$ uniformly over $r\in U$ and
$t\in\operatorname{supp}(H)$.
Hence, the function
\[
    r
    \mapsto
    z_{\overline F}\left(-\frac1r\right)
    =
    r
    +
    \gamma
    \int_{[0,\infty)}
    \frac{rt}{r-t}
    \,dH(t)
\]
is holomorphic on $U$, since the integrand is holomorphic
in $r$ and uniformly bounded for $r\in U$ and
$t\in\operatorname{supp}(H)$. 
Differentiating under the integral sign, we obtain,
for every $r\in U$,
\begin{align*}
    \frac{d}{dr}
    z_{\overline F}\left(-\frac1r\right)
    &=
    1
    -
    \gamma
    \int_{[0,\infty)}
    \frac{t^2}{(r-t)^2}
    \,dH(t),
    \\
    \frac{d^2}{dr^2}
    z_{\overline F}\left(-\frac1r\right)
    &=
    2\gamma
    \int_{[0,\infty)}
    \frac{t^2}{(r-t)^3}
    \,dH(t),
    \\
    \frac{d^3}{dr^3}
    z_{\overline F}\left(-\frac1r\right)
    &=
    -6\gamma
    \int_{[0,\infty)}
    \frac{t^2}{(r-t)^4}
    \,dH(t).
\end{align*}
At $r=r_0$,
\[
    \left.
    \frac{d^3}{dr^3}
    z_{\overline F}\left(-\frac1r\right)
    \right|_{r=r_0}
    <0,
\]
because the integrand is nonnegative and strictly positive
for $t>0$, and $H((0,\infty))>0$ by
Assumption~\ref{ass:GR-basic-4}. Since $r_0$ is real and does not belong to
$\operatorname{supp}(H)$, the functions $z_{\overline F}$
and $x_{\overline F}$ agree on a real neighborhood
of $v_0(s_0)$.
Therefore, evaluating at $r=r_0$, using the chain rule
and \eqref{eq:GR-holder-zero-density-critical}, yields
\[
    \left.
    \frac{d}{dr}
    z_{\overline F}\left(-\frac1r\right)
    \right|_{r=r_0}
    =
    \frac1{r_0^2}
    \left.
    \frac{d}{dw}x_{\overline F}(w)
    \right|_{w=v_0(s_0)}
    =
    0.
\]
Using \eqref{eq:GR-holder-zero-density-image}, we also have
\[
    z_{\overline F}\left(-\frac1{r_0}\right)-s_0
    =
    x_{\overline F}(v_0(s_0))-s_0
    =
    0.
\]
Together with the preceding display, a Taylor expansion
around $r_0$ shows that, if
\[
    \left.
    \frac{d^2}{dr^2}
    z_{\overline F}\left(-\frac1r\right)
    \right|_{r=r_0}
    \neq0,
\]
then
\[
    z_{\overline F}\left(-\frac1r\right)-s_0
    =
    \frac12
    \left.
    \frac{d^2}{dr^2}
    z_{\overline F}\left(-\frac1r\right)
    \right|_{r=r_0}
    (r-r_0)^2
    +
    o\bigl(|r-r_0|^2\bigr),
\]
so $r_0$ is a zero of order two of
$r\mapsto z_{\overline F}(-1/r)-s_0$.
If instead the second derivative vanishes at $r_0$,
then the strictly negative third derivative gives
\[
    z_{\overline F}\left(-\frac1r\right)-s_0
    =
    \frac16
    \left.
    \frac{d^3}{dr^3}
    z_{\overline F}\left(-\frac1r\right)
    \right|_{r=r_0}
    (r-r_0)^3
    +
    o\bigl(|r-r_0|^3\bigr),
\]
so $r_0$ is a zero of order three of
$r\mapsto z_{\overline F}(-1/r)-s_0$. Since $w\mapsto-1/w$ is biholomorphic in a neighborhood
of $v_0(s_0)$ and maps $v_0(s_0)$ to $r_0$,
composition with this map preserves the order of the zero.
Consequently, $w\mapsto z_{\overline F}(w)-s_0$ has a zero
of order $\nu\in\{2,3\}$ at $v_0(s_0)$.
Hence, there exists a function $g$, holomorphic in a complex
neighborhood of $v_0(s_0)$, such that
\begin{align}
\label{eq:GR-holder-local-factorization}
    z_{\overline F}(w)-s_0
    =
    -\frac1w
    +
    \gamma
    \int_{[0,\infty)}
    \frac{t}{1+tw}
    \,dH(t)
    -
    s_0
    =
    \{w-v_0(s_0)\}^{\nu}g(w),
\end{align}
and $g(v_0(s_0))\neq0$.

\textbf{\hyperref[step4:lem16]{Proof of Step~4}}
Since $g$ is holomorphic, it is continuous, and since
$g(v_0(s_0))\neq0$, after decreasing the complex neighborhood
of $v_0(s_0)$ of the preceding step, if necessary,
there exist constants $c_1,c_2>0$ such that
\[
    c_1
    \leq
    |g(w)|
    \leq
    c_2,
\]
for all points $w$ in this neighborhood.
Taking absolute values in
\eqref{eq:GR-holder-local-factorization} therefore gives
\begin{align}
\label{eq:GR-holder-local-power-bounds}
    c_1|w-v_0(s_0)|^\nu
    \leq
    \left|
        -\frac1w
        +
        \gamma
        \int_{[0,\infty)}
        \frac{t}{1+tw}
        \,dH(t)
        -
        s_0
    \right|
    \leq
    c_2|w-v_0(s_0)|^\nu.
\end{align}
Differentiating
\eqref{eq:GR-holder-local-factorization} gives
\begin{align}
\label{eq:diff_fix_point}
    \frac1{w^2}
    -
    \gamma
    \int_{[0,\infty)}
    \frac{t^2}{(1+tw)^2}
    \,dH(t)
    &=
    \{w-v_0(s_0)\}^{\nu-1}
    \left[
        \nu g(w)
        +
        \{w-v_0(s_0)\}\frac{d}{dw}g(w)
    \right].
\end{align}
At $w=v_0(s_0)$, the expression in square brackets equals
\[
    \nu g(v_0(s_0))\neq0.
\]
Since $g$ is holomorphic, both $g$ and its derivative
are continuous. Hence, after decreasing the neighborhood
if necessary,
\[
    \left|
        \nu g(w)
        +
        \{w-v_0(s_0)\}\frac{d}{dw}g(w)
        -
        \nu g(v_0(s_0))
    \right|
    \leq
    \frac12|\nu g(v_0(s_0))|,
\]
and, by the reverse triangle inequality,
\[
    \left|
        \nu g(w)
        +
        \{w-v_0(s_0)\}\frac{d}{dw}g(w)
    \right|
    \geq
    \frac12|\nu g(v_0(s_0))|
    \eqqcolon c_3
    >0
\]
throughout this neighborhood.
Consequently, taking absolute values in
\eqref{eq:diff_fix_point} and using the preceding bounds,
we obtain
\begin{equation}
\label{eq:GR-holder-local-derivative-bound}
    \left|
        \frac1{w^2}
        -
        \gamma
        \int_{[0,\infty)}
        \frac{t^2}{(1+tw)^2}
        \,dH(t)
    \right|
    \geq
    c_3|w-v_0(s_0)|^{\nu-1}.
\end{equation}
By continuity of $v_0$, we may choose $\varepsilon>0$
so small that $v_0(s)$ lies in the preceding neighborhood
whenever $|s-s_0|<\varepsilon$.
We claim that, after decreasing $\varepsilon$ if necessary,
we have
\begin{align}
\label{eq:uniqueness_s_0}
    \left\{
        s\in
        \mathcal S^+
        \cap
        (s_0-\varepsilon,s_0+\varepsilon)
        :
        f_{\gamma,H}(s)=0
    \right\}
    =
    \{s_0\}.
\end{align}
Indeed, suppose the claim is false.
Then, for every $k\in\mathbb N$, there exists
\[
    s_k\in
    \mathcal S^+
    \cap
    (s_0-k^{-1},s_0+k^{-1}),
    \qquad
    s_k\neq s_0,
\]
such that $f_{\gamma,H}(s_k)=0$.
In particular, $s_k\to s_0$.
Since \eqref{eq:GR-zero-density-critical-properties}
has already been proved for every zero-density point,
we have
\[
    \left.
    \frac{d}{dw}x_{\overline F}(w)
    \right|_{w=v_0(s_k)}
    =
    0.
\]
By continuity of $v_0$, $v_0(s_k)\to v_0(s_0)$.
Hence, for all sufficiently large $k$, $v_0(s_k)$
belongs to the neighborhood on which
\eqref{eq:GR-holder-local-derivative-bound} holds.
Since the derivative can vanish there only at $v_0(s_0)$
by \eqref{eq:GR-holder-local-derivative-bound},
it follows that $v_0(s_k)=v_0(s_0)$.
Moreover, \eqref{eq:GR-holder-zero-density-image},
applied at $s_k$, gives
$x_{\overline F}(v_0(s_k))=s_k$.
Therefore,
\[
    s_k
    =
    x_{\overline F}(v_0(s_k))
    =
    x_{\overline F}(v_0(s_0))
    =
    s_0,
\]
a contradiction to $s_k\neq s_0$.

Next, we derive a pointwise estimate for $v_0$ near $s_0$.
For every $s\in\mathcal S^+$ with
$|s-s_0|<\varepsilon$, either $s=s_0$, in which case
\eqref{eq:GR-holder-zero-density-image} applies,
or $f_{\gamma,H}(s)>0$, in which case
\eqref{eq:GR-positive-density-inverse-identity} applies.
Hence,
\begin{align}
\label{eq:GR-holder-local-boundary-identity}
    z_{\overline F}(v_0(s))
    =
    -\frac1{v_0(s)}
    +
    \gamma
    \int_{[0,\infty)}
    \frac{t}{1+t v_0(s)}
    \,dH(t)
    =
    s.
\end{align}
Substituting $w=v_0(s)$ into
\eqref{eq:GR-holder-local-power-bounds} and using
\eqref{eq:GR-holder-local-boundary-identity} gives
\[
    c_1
    |v_0(s)-v_0(s_0)|^\nu
    \leq
    |s-s_0|.
\]
Since $\nu\in\{2,3\}$ and $c_1>0$, it follows that
\begin{align}
\label{eq:GR-holder-pointwise-critical}
    |v_0(s)-v_0(s_0)|
    \leq
    C|s-s_0|^{1/\nu},
\end{align}
where $C=c_1^{-1/\nu}$.

Next, we derive a pointwise estimate for
$\frac{d}{ds}v_0(s)$ near $s_0$, with $s\neq s_0$.
By \eqref{eq:uniqueness_s_0},
$f_{\gamma,H}(s)>0$ for every
\[
    s\in
    \mathcal S^+
    \cap
    (s_0-\varepsilon,s_0+\varepsilon)
    \setminus\{s_0\}.
\]
Recall from Step~1 that $v_0$ is continuously
differentiable at every such $s$.
Differentiating \eqref{eq:GR-holder-local-boundary-identity}
with respect to $s$ gives
\[
    \left.
    \frac{d}{dw}z_{\overline F}(w)
    \right|_{w=v_0(s)}
    \frac{d}{ds}v_0(s)
    =
    1.
\]
Since the first factor is nonzero by
\eqref{eq:GR-positive-density-nonzero-derivative}, we obtain
\[
    \frac{d}{ds}v_0(s)
    =
    \left[
        \frac1{v_0(s)^2}
        -
        \gamma
        \int_{[0,\infty)}
        \frac{t^2}
             {\{1+t v_0(s)\}^2}
        \,dH(t)
    \right]^{-1}.
\]
Therefore, \eqref{eq:GR-holder-local-derivative-bound} gives
\[
    \left|
        \frac{d}{ds}v_0(s)
    \right|
    \leq
    c_3^{-1}
    |v_0(s)-v_0(s_0)|^{1-\nu}.
\]
On the other hand, the upper bound in
\eqref{eq:GR-holder-local-power-bounds}, together with
\eqref{eq:GR-holder-local-boundary-identity}, gives
\[
    |s-s_0|
    \leq
    c_2
    |v_0(s)-v_0(s_0)|^\nu,
\]
and hence
\[
    |v_0(s)-v_0(s_0)|
    \geq
    c_2^{-1/\nu}|s-s_0|^{1/\nu}.
\]
Since $1-\nu<0$, it follows, after increasing $C$
if necessary, that
\begin{equation}
\label{eq:GR-holder-derivative-singularity}
    \left|
        \frac{d}{ds}v_0(s)
    \right|
    \leq
    C|s-s_0|^{-1+1/\nu}.
\end{equation}

\textbf{\hyperref[step5:lem16]{Proof of Step~5}}
After decreasing $\varepsilon$ if necessary,
we assume that $\varepsilon<s_0$, so that
\[
    (s_0-\varepsilon,s_0+\varepsilon)
    \subset
    (0,\infty).
\]
Next, we show that each of the sets
$\mathcal S^+\cap(s_0-\varepsilon,s_0)$ and
$\mathcal S^+\cap(s_0,s_0+\varepsilon)$ is connected
whenever it is nonempty.
Consider first $I_+\coloneqq(s_0,s_0+\varepsilon)$.
Since $I_+\subset(0,\infty)$,
$
    \mathcal S^+\cap I_+
    =
    \mathcal S\cap I_+
$.
As $\mathcal S$ is closed in $\mathbb R$,
$\mathcal S^+\cap I_+$ is relatively closed in $I_+$.
On the other hand, by \eqref{eq:uniqueness_s_0},
$f_{\gamma,H}(s)>0$ for every $s\in\mathcal S^+\cap I_+$.
Since $f_{\gamma,H}$ is continuous on $(0,\infty)$ by
\cite[Theorem~1.1]{SilversteinChoi1995},
every such $s$ has a neighborhood in $I_+$ on which
$f_{\gamma,H}$ remains strictly positive.
Hence, this neighborhood is contained in
$\mathcal S^+\cap I_+$, and therefore
$\mathcal S^+\cap I_+$ is also relatively open in $I_+$.
Since $I_+$ is connected and $\mathcal S^+\cap I_+$
is both relatively open and relatively closed in $I_+$,
the latter set must be either empty or equal to $I_+$.
Hence, whenever it is nonempty, it equals $I_+$
and is therefore connected.
The same argument applies to
$I_-\coloneqq(s_0-\varepsilon,s_0)$.
Thus, each of the sets
$\mathcal S^+\cap I_-$ and
$\mathcal S^+\cap I_+$ is either empty
or equal to the corresponding interval, and is therefore
connected whenever it is nonempty.

Next, let
$
    s,t\in
    \mathcal S^+
    \cap
    (s_0-\varepsilon,s_0+\varepsilon)
$.
Suppose first that $s,t\neq s_0$ and that $s$ and $t$
lie on the same side of $s_0$.
Then, for the corresponding interval $I_\pm$, we have
$\mathcal S^+\cap I_\pm\neq\varnothing$, and hence
$\mathcal S^+\cap I_\pm=I_\pm$.
Moreover, by \eqref{eq:uniqueness_s_0},
$f_{\gamma,H}$ is strictly positive throughout this interval.
Consequently, $v_0$ is continuously differentiable there,
and \eqref{eq:GR-holder-derivative-singularity} gives,
for $\nu\in\{2,3\}$,
\begin{align*}
    |v_0(s)-v_0(t)|
    &\leq
    \int_{\min\{s,t\}}^{\max\{s,t\}}
    \left|
        \frac{d}{du}v_0(u)
    \right|
    \,du
    \\
    &\leq
    C
    \int_{\min\{s,t\}}^{\max\{s,t\}}
    |u-s_0|^{-1+1/\nu}
    \,du
    \\
    &=
    C\nu
    \left|
        |s-s_0|^{1/\nu}
        -
        |t-s_0|^{1/\nu}
    \right|
    \\
    &\leq
    C\nu|s-t|^{1/\nu}.
\end{align*}
Here, in the last step we used
$|a^{1/\nu}-b^{1/\nu}|\leq|a-b|^{1/\nu}$,
$a,b\geq0$, together with
$\bigl||s-s_0|-|t-s_0|\bigr|=|s-t|$,
which holds because $s$ and $t$ lie on the same side of $s_0$.
If one of the two points equals $s_0$, the same estimate
follows directly from \eqref{eq:GR-holder-pointwise-critical}.

Suppose next that $s$ and $t$ lie on opposite sides of $s_0$.
Then, by the triangle inequality and
\eqref{eq:GR-holder-pointwise-critical},
\begin{align*}
    |v_0(s)-v_0(t)|
    &\leq
    |v_0(s)-v_0(s_0)|
    +
    |v_0(t)-v_0(s_0)|
    \\
    &\leq
    C
    \left(
        |s-s_0|^{1/\nu}
        +
        |t-s_0|^{1/\nu}
    \right)
    \\
    &\leq
    C2^{1-1/\nu}
    \left(
        |s-s_0|
        +
        |t-s_0|
    \right)^{1/\nu}
    \\
    &=
    C2^{1-1/\nu}|s-t|^{1/\nu}.
\end{align*}
Thus, in all cases, after increasing $C$ if necessary,
\[
    |v_0(s)-v_0(t)|
    \leq
    C|s-t|^{1/\nu}.
\]
Since $\nu\in\{2,3\}$, after decreasing $\varepsilon$
if necessary so that $\varepsilon\leq1/2$,
we have $|s-t|\leq1$ and hence,
\[
    |s-t|^{1/\nu}
    \leq
    |s-t|^{1/3}.
\]
Consequently,
\[
    |v_0(s)-v_0(t)|
    \leq
    C_{s_0}|s-t|^{1/3},
    \qquad
    s,t\in
    \mathcal S^+
    \cap
    (s_0-\varepsilon,s_0+\varepsilon).
\]
Hence, $v_0$ is $1/3$-H\"older in a neighborhood of every
$s_0\in\mathcal S^+$ satisfying $f_{\gamma,H}(s_0)=0$.
If instead $s_0\in\mathcal S^+$ satisfies
$f_{\gamma,H}(s_0)>0$, then
\eqref{eq:GR-positive-density-local-Lipschitz}
gives Lipschitz continuity in a neighborhood of $s_0$.
After decreasing the corresponding neighborhood
so that its diameter is at most one,
\[
    |v_0(s)-v_0(t)|
    \leq
    C_{s_0}|s-t|
    \leq
    C_{s_0}|s-t|^{1/3}.
\]
Thus, $v_0$ is $1/3$-H\"older in a neighborhood
of every $s\in\mathcal S^+$.

\textbf{\hyperref[step6:lem16]{Proof of Step~6}}
By Lemma~\ref{lem:positive-spectrum-separated-zero},
there exists $\delta_0>0$ such that
$
    (0,\delta_0)\cap\mathcal S
    =
    \varnothing
$.
Moreover, Lemma~\ref{lem:GR-boundary-transform-nonzero}
shows that $\mathcal S$ is compact.
Since $\mathcal S^+=\mathcal S\cap[\delta_0,\infty)$,
it follows that $\mathcal S^+$ is compact.
For every $s_0\in\mathcal S^+$,
choose $\varepsilon_{s_0}>0$ and $C_{s_0}<\infty$ such that
\[
    |v_0(s)-v_0(t)|
    \leq
    C_{s_0}|s-t|^{1/3},
\]
for all
$s,t\in\mathcal S^+\cap
(s_0-2\varepsilon_{s_0},s_0+2\varepsilon_{s_0})$.
The intervals
$
    (s_0-\varepsilon_{s_0},s_0+\varepsilon_{s_0})$,
    $s_0\in\mathcal S^+,
$
form an open cover of $\mathcal S^+$.
By compactness, there exist $s_1,\ldots,s_m\in\mathcal S^+$
such that
\[
    \mathcal S^+
    \subset
    \bigcup_{j=1}^m
    (s_j-\varepsilon_{s_j},s_j+\varepsilon_{s_j}).
\]
Set
\[
    \delta
    \coloneqq
    \min_{1\leq j\leq m}\varepsilon_{s_j}
    >0,
    \qquad
    C
    \coloneqq
    \max_{1\leq j\leq m}C_{s_j}
    <\infty.
\]
Let $s,t\in\mathcal S^+$ with $|s-t|<\delta$,
and choose $j$ such that $|s-s_j|<\varepsilon_{s_j}$.
Then
\[
    |t-s_j|
    \leq
    |t-s|+|s-s_j|
    <
    \delta+\varepsilon_{s_j}
    \leq
    2\varepsilon_{s_j}.
\]
Thus $s$ and $t$ both belong to
$\mathcal S^+\cap
(s_j-2\varepsilon_{s_j},s_j+2\varepsilon_{s_j})$,
and hence
\[
    |v_0(s)-v_0(t)|
    \leq
    C|s-t|^{1/3}.
\]

Finally, continuity of $v_0$ and compactness of
$\mathcal S^+$ imply that
\[
    M_v
    \coloneqq
    \sup_{s\in\mathcal S^+}|v_0(s)|
    <\infty.
\]
If $|s-t|\geq\delta$, then
$|s-t|^{1/3}\delta^{-1/3}\geq1$ and
\[
    |v_0(s)-v_0(t)|
    \leq
    2M_v
    \leq
    2M_v\delta^{-1/3}|s-t|^{1/3}.
\]
Therefore, with
$C_v\coloneqq\max\{C,2M_v\delta^{-1/3}\}$,
we obtain
\[
    |v_0(s)-v_0(t)|
    \leq
    C_v|s-t|^{1/3},
    \qquad
    s,t\in\mathcal S^+,
\]
which proves \eqref{eq:GR-global-one-third-holder}.
\end{proof}

\begin{remark}
The conclusions of Lemma~\ref{lem:GR-population-gap-holder} are closely
related to the general regularity theory for self-consistent spectral
distributions. In particular,
\cite[Proposition~2.4]{AltErdosKruger2020} establishes
$1/3$-Hölder regularity for solutions of a general Dyson equation under
suitable boundedness and flatness assumptions, while
\cite[Theorem~2.5]{AltErdosKruger2020} shows that a zero-density point
in the support is either an edge with square-root behavior or an
interior cusp with cubic-root behavior. This corresponds to the two
possible orders $\nu=2$ and $\nu=3$ appearing in the proof above.
Their setting and assumptions are different from the present scalar
Marchenko--Pastur setting, so we give a direct argument tailored to the
assumptions used here.

Related results for correlated Wishart matrices are given in
\cite{Hachem2016}. They use the opposite sign convention for the
Stieltjes transform, so that their transform $m$ corresponds to $-v$
in our notation.
For an interior zero-density point, Proposition~1 in
\cite{Hachem2016} identifies the corresponding cubic-root cusp
behavior. In particular, in their notation the inverse map has a
critical point at the preimage of the cusp, with vanishing first and
second derivatives and nonvanishing third derivative. These results characterize the behavior at individual cusps, whereas
Lemma~\ref{lem:GR-population-gap-holder} establishes a uniform
$1/3$-Hölder estimate on the whole positive Marchenko--Pastur support.
\end{remark}

\subsection{Lemma~\ref{lem:GR-radon-nikodym-energy} and its proof}

\begin{lemma}
\label{lem:GR-radon-nikodym-energy}
Suppose that Assumptions~\ref{ass:GR-basic},\ref{eq:GR-main-positive-population-gap},
\ref{eq:GR-main-zero-density-regularity}, and
\ref{eq:GR-main-no-critical-signal} hold. Then
\begin{equation}
\label{eq:GR-energy-finite}
    \mathcal E
    \coloneqq
    \int_{[0,M]^2}
    |\kappa(s,t)|
    \,d\mu(s)
    \,d\mu(t)
    <\infty.
\end{equation}
\end{lemma}

\begin{proof}
Recall that $\mathcal S^+=\mathcal S\cap(0,\infty)$.
Under Assumptions~\ref{ass:GR-basic-1},
\ref{ass:GR-basic-4}, and
\ref{eq:GR-main-positive-population-gap},
Lemma~\ref{lem:positive-spectrum-separated-zero}
yields a $\delta_0>0$ such that
\begin{equation}
\label{eq:GR-representation-support-separated}
    (0,\delta_0)\cap\mathcal S=\varnothing.
\end{equation}
Moreover, Lemma~\ref{lem:GR-boundary-transform-nonzero}
shows that $\mathcal S$ is compact.
Since $\mathcal S^+=\mathcal S\cap[\delta_0,\infty)$,
the set $\mathcal S^+$ is compact as well.
The assumptions of
Lemma~\ref{lem:GR-population-gap-holder} are satisfied,
and hence, there exists a $C_v<\infty$ such that
\begin{align}
\label{eq:GR-representation-v-holder}
    |v_0(s)-v_0(u)|
    \leq
    C_v|s-u|^{1/3},
    \qquad
    s,u\in\mathcal S^+.
\end{align}
We first establish that
\begin{align}
\label{eq:GR-representation-denominator-separated}
    \inf_{\substack{
        s\in\mathcal S^+\\
        t\in\operatorname{supp}(G)
    }}
    |1+t v_0(s)|
    >0.
\end{align}
By Assumption~\ref{ass:GR-basic-3}, $C_\Sigma
    \coloneqq
    \sup_{n\geq1}\|\Sigma_n\|_2
    <\infty$.
Since
$\operatorname{supp}(G_n)\subseteq\operatorname{spec}(\Sigma_n)$,
all $G_n$ are supported on $[0,C_\Sigma]$.
Assumption~\ref{ass:GR-basic-5} gives $G_n\Rightarrow G$,
so $G$ is also supported on $[0,C_\Sigma]$.
In particular, $\operatorname{supp}(G)$ is compact.
Fix $s\in\mathcal S^+$ and $t\in\operatorname{supp}(G)$.
Suppose first that $f_{\gamma,H}(s)>0$.
By \eqref{eq:GR-density-boundary-value},
\[
    \Im\{v_0(s)\}
    =
    \pi\gamma f_{\gamma,H}(s)
    >0.
\]
If $t=0$, then $1+t v_0(s)=1$.
If $t>0$, then
\[
    \Im\{1+t v_0(s)\}
    =
    t\Im\{v_0(s)\}
    >0,
\]
and thus $1+t v_0(s)\neq0$.

Suppose next that $f_{\gamma,H}(s)=0$.
By Lemma~\ref{lem:GR-population-gap-holder},
\begin{equation*}\label{eq:GR-representation-zero-density-critical}
    v_0(s)\in\mathbb R\setminus\{0\},
    \qquad
    -\frac1{v_0(s)}
    \notin
    \operatorname{supp}(H),
    \qquad
    \left.
    \frac{d}{dw}x_{\overline F}(w)
    \right|_{w=v_0(s)}
    =
    0.
\end{equation*}
Suppose, to the contrary, that $1+t v_0(s)=0$.
Then $t=-1/v_0(s)$.
Since $t\in\operatorname{supp}(G)$ and
$f_{\gamma,H}(s)=0$, this contradicts
Assumption~\ref{eq:GR-main-no-critical-signal}.
Thus $1+t v_0(s)\neq0$ for every $s\in\mathcal S^+$
and every $t\in\operatorname{supp}(G)$. Since $\mathcal S^+$ and $\operatorname{supp}(G)$ are compact,
and $v_0$ is continuous on $\mathcal S^+$, the map
\[
    (s,t)\mapsto |1+t v_0(s)|
\]
is continuous on the compact set
$\mathcal S^+\times\operatorname{supp}(G)$.
It therefore attains its minimum there.
Since this minimum cannot be zero by the preceding argument,
we obtain
\eqref{eq:GR-representation-denominator-separated}.

Next, we construct a well-behaved representative of the Radon--Nikodym derivative of
$\mathsf C$ with respect to $\mu$.
Recall from Lemma~\ref{lem:GR-auxiliary-positive-measure} that
\begin{equation}
\label{eq:GR-representation-C-transform}
    m_{\mathsf C}(z)
    =
    -\frac1z
    \int_{[0,\infty)}
    \frac{t}{1+t v(z)}
    \,dG(t).
\end{equation}
Moreover, by Lemma~\ref{lem:GR-Stieltjes-convergence},
\begin{equation}
\label{eq:GR-representation-Q-transform}
    m_{\mathsf Q}(z)
    =
    -\frac1\gamma
    \left(
        1+\frac1{zv(z)}
    \right).
\end{equation}
Let
\[
    \mathsf C
    =
    \mathsf C^{\mathrm{ac}}
    +
    \mathsf C^{\mathrm{s}},
    \qquad
    \mathsf Q
    =
    \mathsf Q^{\mathrm{ac}}
    +
    \mathsf Q^{\mathrm{s}}
\]
denote the Lebesgue decompositions of
$\mathsf C$ and $\mathsf Q$ with respect to
Lebesgue measure on $\mathbb R$.
Thus $\mathsf C^{\mathrm{ac}}$ and $\mathsf Q^{\mathrm{ac}}$
are absolutely continuous with respect to Lebesgue measure,
whereas $\mathsf C^{\mathrm{s}}$ and $\mathsf Q^{\mathrm{s}}$
are singular with respect to Lebesgue measure.
Since $\mu=\mathsf C+\mathsf Q$,
its Lebesgue decomposition satisfies
\[
    \mu^{\mathrm{ac}}
    =
    \mathsf C^{\mathrm{ac}}
    +
    \mathsf Q^{\mathrm{ac}}.
\]
Since $s>0$ and $v_0(s)\neq0$ for every $s\in\mathcal S^+$
by Lemma~\ref{lem:GR-boundary-transform-nonzero},
Stieltjes inversion in
\cite[Theorem~B.8]{BaiSilverstein2010}, applied to
\eqref{eq:GR-representation-Q-transform}, gives
\begin{align*}\label{eq:GR-representation-Q-density}
    \frac{d\mathsf Q^{\mathrm{ac}}}{ds}(s)
    =
    \frac1\pi
    \lim_{\eta\downarrow0}
    \Im\bigl\{
        m_{\mathsf Q}(s+i\eta)
    \bigr\}
    =
    \frac{f_{\gamma,H}(s)}
         {s|v_0(s)|^2}
\end{align*}
for Lebesgue-almost every $s\in\mathcal S^+$. By \eqref{eq:GR-representation-denominator-separated}
and dominated convergence, the limit as $z\to s$,
$z\in\mathbb C^+$, may be taken under the integral
with respect to $G$ in \eqref{eq:GR-representation-C-transform}
for every $s\in\mathcal S^+$.
Using $\Im\{t/(1+t v_0(s))\}=
-(t^2\Im\{v_0(s)\})/|1+t v_0(s)|^2$
and $\Im\{v_0(s)\}=\pi\gamma f_{\gamma,H}(s)$, the same Stieltjes inversion applied to
\eqref{eq:GR-representation-C-transform} therefore gives
\begin{equation}
\label{eq:GR-representation-C-density}
    \frac{d\mathsf C^{\mathrm{ac}}}{ds}(s)
    =
    \frac1\pi
    \lim_{\eta\downarrow0}
    \Im\bigl\{
        m_{\mathsf C}(s+i\eta)
    \bigr\}
    =
    \frac{\gamma f_{\gamma,H}(s)}s
    \int_0^\infty
    \frac{t^2}{|1+t v_0(s)|^2}
    \,dG(t)
\end{equation}
for Lebesgue-almost every $s\in\mathcal S^+$. Therefore,
\begin{align}
\label{eq:GR-representation-mu-density}
    \frac{d\mu^{\mathrm{ac}}}{ds}(s)
    =
    \frac{f_{\gamma,H}(s)}s
    \left[
        \frac1{|v_0(s)|^2}
        +
        \gamma
        \int_0^\infty
        \frac{t^2}{|1+t v_0(s)|^2}
        \,dG(t)
    \right]
\end{align}
for Lebesgue-almost every $s\in\mathcal S^+$.

We next show that this density is uniformly bounded.
By \eqref{eq:GR-representation-support-separated},
$s\geq\delta_0$ on $\mathcal S^+$.
Moreover, by \cite[Theorem~1.1]{SilversteinChoi1995},
$v_0$ is continuous on $\mathcal S^+$, while
Lemma~\ref{lem:GR-boundary-transform-nonzero} shows that
$v_0(s)\neq0$ for every $s\in\mathcal S^+$ and that
$\mathcal S$ is compact.
Since $\mathcal S^+=\mathcal S\cap[\delta_0,\infty)$,
it follows that $\mathcal S^+$ is compact and
\begin{align}
\label{eq:bound_step0}
    \inf_{s\in\mathcal S^+}|v_0(s)|>0.
\end{align}
Using \eqref{eq:GR-representation-denominator-separated}
and the compactness of $\operatorname{supp}(G)$,
it follows that
\begin{align*}\sup_{s\in\mathcal S^+}
    \int_0^\infty
    \frac{t^2}{|1+t v_0(s)|^2}
    \,dG(t)
    <\infty.
\end{align*}
Finally,
\begin{align}
\label{eq:bound_step2}
    f_{\gamma,H}(s)
    =
    \frac1{\pi\gamma}\Im\{v_0(s)\}
\end{align}
is bounded on $\mathcal S^+$ because $v_0$ is continuous there.
Combining
\eqref{eq:GR-representation-mu-density}--\eqref{eq:bound_step2},
there exists $C_\mu<\infty$ such that
\begin{equation}
\label{eq:GR-representation-mu-density-bounded}
    \frac{d\mu^{\mathrm{ac}}}{ds}(s)
    \leq C_\mu
\end{equation}
for Lebesgue-almost every $s\in\mathcal S^+$.
We now define $\rho$ on $\mathcal S^+$ by
\begin{equation}
\label{eq:GR-representation-rho-explicit}
    \rho(s)
    \coloneqq
    \frac{
        \displaystyle
        \gamma|v_0(s)|^2
        \int_0^\infty
        \frac{t^2}{|1+t v_0(s)|^2}
        \,dG(t)
    }{
        \displaystyle
        1+
        \gamma|v_0(s)|^2
        \int_0^\infty
        \frac{t^2}{|1+t v_0(s)|^2}
        \,dG(t)
    }.
\end{equation}
For Lebesgue-almost every $s\in\mathcal S^+$ such that
$f_{\gamma,H}(s)>0$, Equations
\eqref{eq:GR-representation-C-density} and
\eqref{eq:GR-representation-mu-density} give
\begin{align*}
    \frac{
        \displaystyle
        \frac{d\mathsf C^{\mathrm{ac}}}{ds}(s)
    }{
        \displaystyle
        \frac{d\mu^{\mathrm{ac}}}{ds}(s)
    }
    &=
    \frac{
        \displaystyle
        \gamma|v_0(s)|^2
        \int_0^\infty
        \frac{t^2}{|1+t v_0(s)|^2}
        \,dG(t)
    }{
        \displaystyle
        1+
        \gamma|v_0(s)|^2
        \int_0^\infty
        \frac{t^2}{|1+t v_0(s)|^2}
        \,dG(t)
    }
    =
    \rho(s).
\end{align*}
On the other hand,
\eqref{eq:GR-representation-mu-density} implies that
\[
    \frac{d\mu^{\mathrm{ac}}}{ds}(s)=0
\]
for Lebesgue-almost every $s\in\mathcal S^+$
where $f_{\gamma,H}(s)=0$.
Hence, $\{s\in\mathcal S^+:f_{\gamma,H}(s)=0\}$
is a $\mu^{\mathrm{ac}}$-null set.
Consequently,
\begin{equation}
\label{eq:GR-representation-rho-ac-RN}
    \rho
    =
    \frac{
        d\mathsf C^{\mathrm{ac}}
    }{
        d\mu^{\mathrm{ac}}
    }
    \qquad
    \mu^{\mathrm{ac}}\text{-a.e. on }\mathcal S^+.
\end{equation}

We next establish H\"older continuity of $\rho$
on $\mathcal S^+$.
By \eqref{eq:GR-representation-denominator-separated},
there exists $c>0$ such that
\[
    |1+t v_0(s)|\geq c
\]
for every $s\in\mathcal S^+$ and every
$t\in\operatorname{supp}(G)$.
Moreover, since $\operatorname{supp}(G)$ and $\mathcal S^+$
are compact and $v_0$ is continuous on $\mathcal S^+$,
there exist constants $M_G,M_v<\infty$ such that
\[
    t\leq M_G,
    \qquad
    |v_0(s)|\leq M_v
\]
for every $t\in\operatorname{supp}(G)$ and every
$s\in\mathcal S^+$. Fix $s,u\in\mathcal S^+$.
For every $t\in\operatorname{supp}(G)$, we have
\begin{align*}
    \left|
        \frac{t v_0(s)}{1+t v_0(s)}
        -
        \frac{t v_0(u)}{1+t v_0(u)}
    \right|
    =
    \frac{
        t|v_0(s)-v_0(u)|
    }{
        |1+t v_0(s)|\,|1+t v_0(u)|
    }
    \leq
    \frac{M_G}{c^2}|v_0(s)-v_0(u)|.
\end{align*}
Furthermore,
\[
    \left|
        \frac{t v_0(s)}{1+t v_0(s)}
    \right|
    +
    \left|
        \frac{t v_0(u)}{1+t v_0(u)}
    \right|
    \leq
    \frac{2M_GM_v}{c}.
\]
Hence, using
$\bigl||a|^2-|b|^2\bigr|\leq|a-b|(|a|+|b|)$
for $a,b\in\mathbb C$, we obtain
\begin{align*}
    &
    \left|
        |v_0(s)|^2
        \int_0^\infty
        \frac{t^2}{|1+t v_0(s)|^2}
        \,dG(t)
        -
        |v_0(u)|^2
        \int_0^\infty
        \frac{t^2}{|1+t v_0(u)|^2}
        \,dG(t)
    \right|
    \\
    &\quad\leq
    \int_0^\infty
    \left|
        \left|
            \frac{t v_0(s)}{1+t v_0(s)}
        \right|^2
        -
        \left|
            \frac{t v_0(u)}{1+t v_0(u)}
        \right|^2
    \right|
    \,dG(t)
    \\
    &\quad\leq
    \frac{2M_G^2M_v}{c^3}|v_0(s)-v_0(u)|.
\end{align*}
Since, for every $a,b,\lambda\geq0$,
\[
    \left|
        \frac{\lambda a}{1+\lambda a}
        -
        \frac{\lambda b}{1+\lambda b}
    \right|
    =
    \frac{
        \lambda|a-b|
    }{
        (1+\lambda a)(1+\lambda b)
    }
    \leq
    \lambda|a-b|,
\]
the definition of $\rho$ now gives
\[
    |\rho(s)-\rho(u)|
    \leq
    \frac{2\gamma M_G^2M_v}{c^3}
    |v_0(s)-v_0(u)|
\]
for all $s,u\in\mathcal S^+$.
Set $C_\rho=2\gamma M_G^2M_v/c^3$.
Combining this with
\eqref{eq:GR-representation-v-holder} yields
\begin{equation}
\label{eq:GR-representation-rho-holder}
    |\rho(s)-\rho(u)|
    \leq
    C_\rho C_v|s-u|^{1/3},
    \qquad
    s,u\in\mathcal S^+.
\end{equation}
Hence, the $1/3$-H\"older continuity of $v_0$ on
$\mathcal S^+$ implies the $1/3$-H\"older continuity
of $\rho$ on $\mathcal S^+$.

We can now control the contribution of the absolutely
continuous part of $\mu$ to \eqref{eq:GR-energy-finite}
on $\mathcal S^+$.
By \eqref{eq:GR-representation-mu-density-bounded},
\eqref{eq:GR-representation-rho-holder},
and the definition of $\kappa$,
\begin{align}
\label{eq:GR-representation-energy-ac}
\begin{split}
    \int_{(\mathcal S^+)^2}
    |\kappa(s,u)|
    \,d\mu^{\mathrm{ac}}(s)
    \,d\mu^{\mathrm{ac}}(u)
    &\leq
    C_\rho C_v C_\mu^2
    \int_{\delta_0}^M
    \int_{\delta_0}^M
    |s-u|^{-2/3}
    \,ds\,du
    \\
    &=
    \frac92
    C_\rho C_v C_\mu^2
    (M-\delta_0)^{4/3}
    <\infty.
\end{split}
\end{align}
Here $\kappa(s,s)=0$, and the diagonal has zero measure
with respect to
$\mu^{\mathrm{ac}}\otimes\mu^{\mathrm{ac}}$.

It remains to control the remaining contributions to
\eqref{eq:GR-energy-finite}.
We first exclude singular mass on $\mathcal S^+$
and identify a neighborhood on which $\mu$ has no mass
outside $\mathcal S^+$.
By \eqref{eq:GR-representation-support-separated}, \eqref{eq:GR-representation-denominator-separated}, and
\eqref{eq:bound_step0},
\[
    a
    \coloneqq
    \inf_{s\in\mathcal S^+}s>0,
    \qquad
    b
    \coloneqq
    \inf_{s\in\mathcal S^+}|v_0(s)|>0,
    \qquad
    d
    \coloneqq
    \inf_{\substack{
        s\in\mathcal S^+\\
        t\in\operatorname{supp}(G)
    }}
    |1+t v_0(s)|>0.
\]
Moreover, since $\operatorname{supp}(G)$ is compact, $T\coloneqq\sup\operatorname{supp}(G)<\infty$.
Set $\varepsilon
    \coloneqq
    \min\left\{
        \frac b2,
        \frac{d}{2(1+T)}
    \right\}$.
By \cite[Theorem~1.1]{SilversteinChoi1995},
for every $s\in\mathcal S^+$,
\[
    \lim_{\substack{z\to s\\z\in\mathbb C^+}}v(z)
    =
    v_0(s).
\]
Hence, for every $s\in\mathcal S^+$, there exists
$r_s>0$, which we may choose so that $r_s<a/2$, such that
\[
    |v(z)-v_0(s)|<\varepsilon
\]
for $z\in\mathbb C^+$ and $|z-s|<r_s$.
For such $z$ and every $t\in\operatorname{supp}(G)$,
\[
    |z|
    \geq
    s-|z-s|
    \geq
    \frac a2,
    \qquad
    |v(z)|
    \geq
    |v_0(s)|-|v(z)-v_0(s)|
    \geq
    \frac b2,
\]
and
\[
    |1+t v(z)|
    \geq
    |1+t v_0(s)|
    -
    t|v(z)-v_0(s)|
    \geq
    d-T\varepsilon
    \geq
    \frac d2.
\]
By the compactness of $\mathcal S^+$,
finitely many intervals $(
        s_j-\frac{r_{s_j}}2,
        s_j+\frac{r_{s_j}}2
    )$, with $j=1,\ldots,m$,
cover $\mathcal S^+$.
Set
\[
    V
    \coloneqq
    \bigcup_{j=1}^m
    \left(
        s_j-\frac{r_{s_j}}2,
        s_j+\frac{r_{s_j}}2
    \right),
    \qquad
    \eta_0
    \coloneqq
    \frac12\min_{1\leq j\leq m}r_{s_j},
    \qquad
    c
    \coloneqq
    \frac12\min\{a,b,d\}.
\]
Then $V\subset(0,\infty)$ is an open neighborhood
of $\mathcal S^+$.
If $x\in V$ and $0<\eta\leq\eta_0$, then for some
$j\in\{1,\ldots,m\}$,
\[
    |x-s_j|<\frac{r_{s_j}}2
\]
and therefore
\[
    |x+i\eta-s_j|
    \leq
    |x-s_j|+\eta
    <
    r_{s_j}.
\]
The preceding estimates consequently give
\[
    |x+i\eta|\geq c,
    \qquad
    |v(x+i\eta)|\geq c,
    \qquad
    |1+t v(x+i\eta)|\geq c
\]
for every $x\in V$, every $0<\eta\leq\eta_0$,
and every $t\in\operatorname{supp}(G)$. Using the representations of $m_{\mathsf C}$ and
$m_{\mathsf Q}$ in \eqref{eq:GR-representation-C-transform} and
\eqref{eq:GR-representation-Q-transform}, respectively, we obtain
\[
    |m_{\mathsf Q}(x+i\eta)|
    \leq
    \frac1\gamma
    \left(
        1+\frac1{c^2}
    \right)
    \qquad\text{and}\qquad
    |m_{\mathsf C}(x+i\eta)|
    \leq
    \frac1{c^2}
    \int_0^\infty t\,dG(t)
\]
for every $x\in V$ and every $0<\eta\leq\eta_0$. Since $\operatorname{supp}(G)$ is compact, the integral
on the right-hand side is finite.
Hence, there exists $C<\infty$ such that
\[
    |m_{\mathsf Q}(x+i\eta)|
    +
    |m_{\mathsf C}(x+i\eta)|
    \leq C
\]
for every $x\in V$ and every $0<\eta\leq\eta_0$. For $\eta\geq\eta_0$, we use the following Stieltjes bound
for finite positive measures supported on $[0,\infty)$:
\[
    |m_\nu(x+i\eta)|
    =
    \left|
        \int_{[0,\infty)}
        \frac1{s-(x+i\eta)}
        \,d\nu(s)
    \right|
    \leq
    \int_{[0,\infty)}
    \frac1{|s-(x+i\eta)|}
    \,d\nu(s)
    \leq
    \frac{\nu([0,\infty))}{\eta}
    \leq
    \frac{\nu([0,\infty))}{\eta_0},
\]
where $\nu$ denotes either $\mathsf Q$ or $\mathsf C$.
Thus, combining the preceding two displays, we obtain
\begin{equation}
\label{eq:GR-representation-transforms-locally-bounded}
    \sup_{\substack{x\in V\\\eta>0}}
    \left\{
        |m_{\mathsf Q}(x+i\eta)|
        +
        |m_{\mathsf C}(x+i\eta)|
    \right\}
    \leq
    \max\left\{
        C,
        \frac{
            \mathsf Q([0,\infty))
            +
            \mathsf C([0,\infty))
        }{\eta_0}
    \right\}
    <\infty.
\end{equation}
Let $\varphi\in C_c(V)$ and extend $\varphi$ by zero
to $\mathbb R$.
Since $\operatorname{supp}(\varphi)$ is a compact subset
of $V$, it is also a compact subset of $\mathbb R$,
and hence this extension belongs to
$C_c(\mathbb R)\subseteq C_0(\mathbb R)$.
Hence, Proposition~\ref{prop:multivariate-stieltjes-inversion},
applied in dimension one, gives
\[
    \int_{[0,\infty)}
    \varphi(s)\,d\mathsf Q(s)
    =
    \lim_{\eta\downarrow0}
    \frac1\pi
    \int_{[0,\infty)}
    \varphi(s)
    \Im\bigl\{
        m_{\mathsf Q}(s+i\eta)
    \bigr\}
    \,ds
\]
and
\[
    \int_{[0,\infty)}
    \varphi(s)\,d\mathsf C(s)
    =
    \lim_{\eta\downarrow0}
    \frac1\pi
    \int_{[0,\infty)}
    \varphi(s)
    \Im\bigl\{
        m_{\mathsf C}(s+i\eta)
    \bigr\}
    \,ds.
\]
By \cite[Theorem~1.1]{SilversteinChoi1995},
$v(s+i\eta)\to v_0(s)$ as $\eta\downarrow0$
for every $s\in V$.
Using the representations \eqref{eq:Q:stielt:def} and
\eqref{eq:GR-C-transform}, together with the preceding
uniform lower bounds on $|v(s+i\eta)|$ and
$|1+t v(s+i\eta)|$, and dominated convergence in
the integral with respect to $G$, we therefore obtain
\[
    \frac1\pi
    \Im\bigl\{
        m_{\mathsf Q}(s+i\eta)
    \bigr\}
    \longrightarrow
    \mathbf1\{s\in\mathcal S^+\}
    \frac{f_{\gamma,H}(s)}
         {s|v_0(s)|^2}
\]
and
\[
    \frac1\pi
    \Im\bigl\{
        m_{\mathsf C}(s+i\eta)
    \bigr\}
    \longrightarrow
    \mathbf1\{s\in\mathcal S^+\}
    \frac{\gamma f_{\gamma,H}(s)}s
    \int_{[0,\infty)}
    \frac{t^2}{|1+t v_0(s)|^2}
    \,dG(t)
\]
for every $s\in V$. Moreover, after increasing $C$ if necessary,
\eqref{eq:GR-representation-transforms-locally-bounded} gives
\[
    \left|
        \frac1\pi
        \varphi(s)
        \Im\bigl\{
            m_{\mathsf Q}(s+i\eta)
        \bigr\}
    \right|
    \leq
    \frac C\pi|\varphi(s)|
\]
and
\[
    \left|
        \frac1\pi
        \varphi(s)
        \Im\bigl\{
            m_{\mathsf C}(s+i\eta)
        \bigr\}
    \right|
    \leq
    \frac C\pi|\varphi(s)|
\]
for every $s\in V$ and every $\eta>0$.
Since $\varphi$ vanishes outside $V$, these bounds hold
on $\mathbb R$, and $C\pi^{-1}|\varphi|$
is Lebesgue integrable. Dominated convergence therefore yields
\[
    \int_{[0,\infty)}
    \varphi(s)\,d\mathsf Q(s)
    =
    \int_V
    \varphi(s)
    \mathbf1\{s\in\mathcal S^+\}
    \frac{f_{\gamma,H}(s)}
         {s|v_0(s)|^2}
    \,ds
\]
and
\[
    \int_{[0,\infty)}
    \varphi(s)\,d\mathsf C(s)
    =
    \int_V
    \varphi(s)
    \mathbf1\{s\in\mathcal S^+\}
    \frac{\gamma f_{\gamma,H}(s)}s
    \left(
        \int_0^\infty
        \frac{t^2}{|1+t v_0(s)|^2}
        \,dG(t)
    \right)
    \,ds.
\]
Since these identities hold for every $\varphi\in C_c(V)$,
it follows that
\begin{align}
\label{eq:no_singular_Q}
    \left.\mathsf Q\right|_V
    =
    \mathbf1\{s\in\mathcal S^+\}
    \frac{f_{\gamma,H}(s)}
         {s|v_0(s)|^2}
    \,ds
\end{align}
and
\begin{align}
\label{eq:no_singular_C}
    \left.\mathsf C\right|_V
    =
    \mathbf1\{s\in\mathcal S^+\}
    \frac{\gamma f_{\gamma,H}(s)}s
    \left(
        \int_0^\infty
        \frac{t^2}{|1+t v_0(s)|^2}
        \,dG(t)
    \right)
    \,ds.
\end{align}
Thus neither $\mathsf Q$ nor $\mathsf C$
has a singular component on $\mathcal S^+$.
Moreover,
\begin{equation}
\label{eq:GR-representation-neighborhood-gap}
    \mu(V\setminus\mathcal S^+)=0.
\end{equation}

We next locate the possible mass of $\mathsf Q$
outside $\mathcal S^+\cup\{0\}$.
Since $\mathsf Q$ is supported on $[0,M]$,
it suffices to consider
$x\in(0,\infty)\setminus\mathcal S^+$.
Since $\operatorname{supp}(\overline F)\subseteq[0,\infty)$
and by \eqref{eq:nonzero-support-identification}, we have
\[
    x\notin\operatorname{supp}(\overline F)\cup\{0\}.
\]
Since $\operatorname{supp}(\overline F)$ is closed,
there exist an open interval $I_x$ containing $x$
and an open complex neighborhood $U$ of $\overline{I_x}$
such that $0\notin U$ and  $U\cap\operatorname{supp}(\overline F)
    =
    \varnothing$.
Since $U\cap\operatorname{supp}(\overline F)=\varnothing$,
the Stieltjes transform $v$ extends holomorphically to $U$.
Moreover, $v$ is real-valued on $U\cap\mathbb R$.
We denote this extension again by $v$.
Its value at $x$ equals $v_0(x)$, and
\cite[Theorem~4.1]{SilversteinChoi1995} gives $v_0(x)\in\mathcal B_H$. In particular, $v_0(x)\neq0$.
After decreasing $I_x$ and $U$ if necessary,
we may therefore assume that $v(z)\neq0$
for all $z\in U$. Define
\[
    \widetilde m_{\mathsf Q}(z)
    \coloneqq
    -\frac1\gamma
    \left(
        1+\frac1{zv(z)}
    \right),
    \qquad z\in U.
\]
Then $\widetilde m_{\mathsf Q}$ is holomorphic on $U$
and, by \eqref{eq:GR-representation-Q-transform},
\[
    \widetilde m_{\mathsf Q}(z)
    =
    m_{\mathsf Q}(z),
    \qquad z\in U\cap\mathbb C^+.
\]
Moreover, $\widetilde m_{\mathsf Q}(y)\in\mathbb R$, for $y\in I_x$. Let $\varphi\in C_c(\mathbb R)$ satisfy
$\operatorname{supp}(\varphi)\subseteq I_x$.
Since $\widetilde m_{\mathsf Q}$ is continuous on $U$
and real-valued on $I_x$, its imaginary part converges
uniformly to zero on $\operatorname{supp}(\varphi)$
as the imaginary argument tends to zero.
Proposition~\ref{prop:multivariate-stieltjes-inversion},
applied in dimension one, therefore yields
\begin{align*}
    \int_{\mathbb R}
    \varphi(s)\,d\mathsf Q(s)
    &=
    \lim_{\eta\downarrow0}
    \frac1\pi
    \int_{\mathbb R}
    \varphi(s)
    \Im\left\{
        m_{\mathsf Q}(s+i\eta)
    \right\}
    \,ds
    \\
    &=
    \lim_{\eta\downarrow0}
    \frac1\pi
    \int_{I_x}
    \varphi(s)
    \Im\left\{
        \widetilde m_{\mathsf Q}(s+i\eta)
    \right\}
    \,ds
    =
    0.
\end{align*}
Since this holds for every such $\varphi$,
it follows that $\left.\mathsf Q\right|_{I_x}=0$.
Since $x\in(0,\infty)\setminus\mathcal S^+$
was arbitrary, and $\mathsf Q$ is supported on $[0,M]$,
we conclude that
\begin{equation}
\label{eq:GR-representation-Q-support}
    \operatorname{supp}(\mathsf Q)
    \subseteq
    \mathcal S^+\cup\{0\}.
\end{equation}
Together with \eqref{eq:no_singular_Q}, this implies
\begin{align*}\operatorname{supp}
    \bigl(
        \mathsf Q^{\mathrm{s}}
    \bigr)
    \subseteq
    \{0\}.
\end{align*}

We now consider the remaining mass of $\mathsf C$.
Set $D\coloneqq(0,M]\setminus\mathcal S^+$. By \eqref{eq:no_singular_Q} and \eqref{eq:no_singular_C}, $\left.\mu\right|_{\mathcal S^+}
    =
    \left.\mu^{\mathrm{ac}}\right|_{\mathcal S^+}$.
Moreover, \eqref{eq:GR-representation-Q-support} gives $\left.\mu\right|_D=\left.\mathsf C\right|_D$.
Thus
\begin{equation}
\label{eq:GR-representation-mu-decomposition}
    \mu=\left.\mu^{\mathrm{ac}}\right|_{\mathcal S^+}+\left.\mathsf C\right|_D+\mu(\{0\})\delta_0.
\end{equation}
We complete the definition of $\rho$ by setting
\begin{equation}
\label{eq:GR-representation-rho-outside}
    \rho(s)\coloneqq1,
    \qquad s\in D,
\end{equation}
and
\begin{equation}
\label{eq:GR-representation-rho-zero}
    \rho(0)\coloneqq
    \begin{cases}
        \displaystyle
        \frac{\mathsf C(\{0\})}{\mu(\{0\})},
        &\mu(\{0\})>0,\\[2mm]
        0,
        &\mu(\{0\})=0.
    \end{cases}
\end{equation}
Together with \eqref{eq:GR-representation-rho-explicit},
this defines a Borel function on $[0,M]$.
On $\mathcal S^+$,
\eqref{eq:GR-representation-rho-ac-RN},
\eqref{eq:no_singular_Q}, and
\eqref{eq:no_singular_C} imply
$
    \left.\mathsf C\right|_{\mathcal S^+}
    =
    \rho\left.\mu\right|_{\mathcal S^+}
$.
On $D$, we have $\left.\mathsf C\right|_D=\left.\mu\right|_D$
and $\rho=1$, whereas at zero the corresponding identity
follows from \eqref{eq:GR-representation-rho-zero}.
Consequently, $\mathsf C=\rho\mu$.
Furthermore,
$0\leq\rho(s)\leq1$, for $s\in[0,M].$.
On $\mathcal S^+$, this follows from
\eqref{eq:GR-representation-rho-explicit};
on $D$ it follows from
\eqref{eq:GR-representation-rho-outside};
and at zero it follows from
$0\leq\mathsf C(\{0\})\leq \mu(\{0\})$.
Thus, $\rho$ is a Radon--Nikodym derivative of $\mathsf C$
with respect to $\mu$. By \eqref{eq:GR-representation-neighborhood-gap},
we have $\mathsf C(V\setminus\mathcal S^+)=0$.
Since $\mathcal S^+$ is compact and contained
in the open set $V$,
\begin{equation}
\label{eq:GR-representation-outside-separated}
    d_V
    \coloneqq
    \operatorname{dist}
    \bigl(
        \mathcal S^+,
        \mathbb R\setminus V
    \bigr)
    >0.
\end{equation}

We next show that
\begin{equation}
\label{eq:GR-representation-C-inverse-moment}
    I_{\mathsf C}
    \coloneqq
    \int_{(0,M]}
    \frac1s\,d\mathsf C(s)
    <\infty.
\end{equation}
By \eqref{eq:GR-representation-support-separated}
and \eqref{eq:nonzero-support-identification}, $\operatorname{supp}(\overline F)
    \subseteq
    \{0\}\cup[\delta_0,M]$.
Set $\alpha\coloneqq\overline F(\{0\})$.
For $\lambda>0$, dominated convergence gives
\[
    v_0(-\lambda)
    \coloneqq
    \lim_{\eta\downarrow0}v(-\lambda+i\eta)
    =
    \frac{\alpha}{\lambda}
    +
    \int_{[\delta_0,M]}
    \frac1{s+\lambda}\,d\overline F(s)
    >0.
\]
Letting $z\to-\lambda$ through $\mathbb C^+$ in
\eqref{eq:GR-representation-C-transform},
dominated convergence gives
\begin{equation}
\label{eq:GR-representation-C-negative-axis}
    m_{\mathsf C}(-\lambda)
    \coloneqq
    \int_{[0,M]}
    \frac1{s+\lambda}\,d\mathsf C(s)
    =
    \frac1\lambda
    \int_{[0,\infty)}
    \frac{t}{1+t v_0(-\lambda)}
    \,dG(t).
\end{equation}
Suppose first that $\alpha>0$.
Since $v_0(-\lambda)\geq\alpha/\lambda$ and
\[
    0\leq
    \frac{t}{1+t v_0(-\lambda)}
    \leq
    \frac1{v_0(-\lambda)},
    \qquad t\geq0,
\]
\eqref{eq:GR-representation-C-negative-axis} yields
\[
    0
    \leq
    m_{\mathsf C}(-\lambda)
    \leq
    \frac1{\lambda v_0(-\lambda)}
    \leq
    \frac1\alpha,
    \qquad \lambda>0.
\]
Since $\mathsf C$ is a positive measure, we have
\[
    0
    \leq
    \frac{\mathsf C(\{0\})}{\lambda}
    \leq
    m_{\mathsf C}(-\lambda)
    \leq
    \frac1\alpha.
\]
Multiplying by $\lambda$ and letting $\lambda\downarrow0$
shows that $\mathsf C(\{0\})=0$.
Consequently,
\[
    m_{\mathsf C}(-\lambda)
    =
    \int_{(0,M]}
    \frac1{s+\lambda}\,d\mathsf C(s).
\]
Since $(s+\lambda)^{-1}\uparrow s^{-1}$ as
$\lambda\downarrow0$ for every $s>0$,
monotone convergence yields
\[
    I_{\mathsf C}
    =
    \int_{(0,M]}
    \frac1s\,d\mathsf C(s)
    =
    \lim_{\lambda\downarrow0}
    m_{\mathsf C}(-\lambda)
    \leq
    \frac1\alpha
    <\infty.
\]

Suppose next that $\alpha=0$. Then
$\overline F$ is supported on $[\delta_0,M]$, and
dominated convergence gives
\[
    v_0(-\lambda)
    \longrightarrow
    v_*
    \coloneqq
    \int_{[\delta_0,M]}
    \frac1s\,d\overline F(s)
    \in(0,\infty).
\]
Moreover,
\begin{align}
\label{eq:GR-representation-v-zero-difference}
    \frac{v_*-v_0(-\lambda)}{\lambda}
    =
    \int_{[\delta_0,M]}
    \frac1{s(s+\lambda)}\,d\overline F(s)
    \leq
    \frac1{\delta_0^2}.
\end{align}
By dominated convergence, letting $\lambda\downarrow0$, yields
\[
    \lambda m_{\mathsf C}(-\lambda)
    =
    \int_{[0,M]}
    \frac{\lambda}{s+\lambda}\,d\mathsf C(s)
    \longrightarrow  \int_{[0,M]}
    1_{\{0\}}(s)\,d\mathsf C(s)=
    \mathsf C(\{0\}).
\]
Together with \eqref{eq:GR-representation-C-negative-axis}
and $v_0(-\lambda)\to v_*$, dominated convergence gives
\begin{equation*}\label{eq:GR-representation-C-zero-mass}
    c_0
    \coloneqq
    \mathsf C(\{0\})
    =
    \lim_{\lambda\downarrow0}
    \lambda m_{\mathsf C}(-\lambda)
    =
    \int_{[0,\infty)}
    \frac{t}{1+t v_*}\,dG(t).
\end{equation*}
Hence, we obtain
\begin{align*}
    \int_{(0,M]}
    \frac1{s+\lambda}\,d\mathsf C(s)
    &=
    m_{\mathsf C}(-\lambda)-\frac{c_0}{\lambda}
    \\
    &=
    \frac1\lambda
    \int_{[0,\infty)}
    \left\{
        \frac{t}{1+t v_0(-\lambda)}
        -
        \frac{t}{1+t v_*}
    \right\}
    \,dG(t)
    \\
    &=
    \frac{v_*-v_0(-\lambda)}{\lambda}
    \int_{[0,\infty)}
    \frac{t^2}
         {\{1+t v_0(-\lambda)\}(1+t v_*)}
    \,dG(t)
    \\
    &\leq
    \frac1{\delta_0^2}
    \int_{[0,\infty)}t^2\,dG(t)
    <\infty,
\end{align*}
where we used \eqref{eq:GR-representation-v-zero-difference},
$v_0(-\lambda),v_*>0$, and the compact support of $G$.
Since $(s+\lambda)^{-1}\uparrow s^{-1}$ as
$\lambda\downarrow0$ for every $s>0$,
monotone convergence yields
\[
    I_{\mathsf C}
    =
    \lim_{\lambda\downarrow0}
    \int_{(0,M]}
    \frac1{s+\lambda}\,d\mathsf C(s)
    \leq
    \frac1{\delta_0^2}
    \int_{[0,\infty)}t^2\,dG(t)
    <\infty.
\]
This proves \eqref{eq:GR-representation-C-inverse-moment}
also in this case.

We finally combine the contributions to
\eqref{eq:GR-energy-finite} according to
\eqref{eq:GR-representation-mu-decomposition}.
Since $\rho=1$ on $D$, we have
$\kappa(s,u)=0$ for all $s,u\in D$.
Moreover, $\kappa(0,0)=0$ by definition.
Thus, the contributions from $D\times D$ and
from $(0,0)$ vanish.

We next consider the mixed contributions involving
$\left.\mu\right|_{\mathcal S^+}$ and $\left.\mu\right|_D$.
Recall that $\left.\mu\right|_D
    =
    \left.\mathsf C\right|_D$.
Since $D\cap V\subseteq V\setminus\mathcal S^+$,
the identity $\mathsf C(V\setminus\mathcal S^+)=0$
implies
\[
    \mu(D\cap V)
    =
    \mathsf C(D\cap V)
    =
    0.
\]
For every $s\in\mathcal S^+$ and $u\in D\setminus V$,
we have $u\in\mathbb R\setminus V$, so
\[
    |s-u|\geq d_V>0
\]
by \eqref{eq:GR-representation-outside-separated}.
Since $\rho(u)=1$ on $D$ and $0\leq\rho\leq1$,
we obtain
\begin{align*}
    \int_{[0,M]^2}
    |\kappa(s,u)|\,\left.d\mu\right|_{\mathcal S^+}(s)\,d\mu|_D(u)
    &=
    \int_{\mathcal S^+}
    \int_{D\setminus V}
    \frac{|\rho(s)-1|}{|s-u|}
    \,d\mathsf C(u)\,\left.d\mu\right|_{\mathcal S^+}(s)
    \\
    &\leq
    \frac1{d_V}
    \left.\mu\right|_{\mathcal S^+}([0,M])\mathsf C(D)
    <\infty.
\end{align*}
By the symmetry of $\kappa$, the other mixed contribution
satisfies the same bound.

It remains to bound the contributions to $\mathcal E$
from $\{0\}\times(0,M]$ and $(0,M]\times\{0\}$.
By the symmetry of $\kappa$ and since $0\leq\rho\leq1$,
their sum satisfies
\begin{align*}
    2\mu(\{0\})
    \int_{(0,M]}
    |\kappa(0,s)|\,d\mu(s)
    &=
    2\mu(\{0\})
    \int_{(0,M]}
    \frac{|\rho(0)-\rho(s)|}{s}\,d\mu(s)
    \\
    &\leq
    2\mu(\{0\})
    \int_{(0,M]}
    \frac1s\,d\mu(s)
    \\
    &\leq
    2\mu(\{0\})
    \left\{
        I_{\mathsf C}
        +
        \frac{\mathsf Q([0,M])}{\delta_0}
    \right\}
    <\infty.
\end{align*}

Combining the preceding estimates with
\eqref{eq:GR-representation-energy-ac}, we conclude that
\begin{align*}
    \mathcal E
    &\leq
    \int_{[0,M]^2}
    |\kappa(s,u)|\,\left.d\mu\right|_{\mathcal S^+}(s)\,\left.d\mu\right|_{\mathcal S^+}(u)
    +
    \frac{2}{d_V}
    \left.\mu\right|_{\mathcal S^+}([0,M])\mathsf C(D)
    \\
    &\quad+
    2\mu(\{0\})
    \left\{
        I_{\mathsf C}
        +
        \frac{\mathsf Q([0,M])}{\delta_0}
    \right\}
    <\infty.
\end{align*}
This proves \eqref{eq:GR-energy-finite}.
\end{proof}

\subsection{Proof of Lemma~\ref{lem:GR-candidate-measure}}

\begin{proof}[Proof of Lemma~\ref{lem:GR-candidate-measure}]
By Lemma~\ref{lem:GR-radon-nikodym-energy}, we have $\mathcal{E} = \|\kappa\|_{L^1(\mu\otimes \mu)} < \infty$. 
Fubini's theorem shows that $\ell$ from \eqref{eq:kappa_density_def_main} is well defined for
$\mu$-almost every $s\in[0,M]$.
We set $\ell(s)=0$ on the $\mu$-null set where
$\int_{[0,M]}|\kappa(s,t)|\,d\mu(t)=\infty$.
Then $\ell$ is a Borel function on $[0,M]$ satisfying
\begin{align}
\label{eq:ell_bound}
    \int_{[0,M]}
    |\ell(s)|
    \,d\mu(s)
    \leq
    \int_{[0,M]}
    \int_{[0,M]}
    |\kappa(s,t)|
    \,d\mu(t)d\mu(s)
    =
    \mathcal E
    <\infty.
\end{align}
Recall the notation introduced after
\eqref{eq:GR-candidate-main}. Using \eqref{eq:ell_bound},
$\kappa(\mu\otimes\mu)$ and $\ell\mu$ are finite signed
Borel measures on $[0,M]^2$ and $[0,M]$, respectively.
Since $\rho\,d\mu=d\mathsf C$ and $\mathsf C$ is a finite
positive Borel measure supported on $[0,M]$ by
Lemma~\ref{lem:GR-auxiliary-positive-measure},
\eqref{eq:GR-candidate-main} defines a finite signed
Borel measure supported on $[0,M]^2$.
Moreover, using the definition of $\mathsf K$ in \eqref{eq:GR-candidate-main},
the fact that
$\|D_\#\nu\|_{\mathrm{TV}}\leq\|\nu\|_{\mathrm{TV}}$
for every finite signed measure $\nu$,
and \eqref{eq:ell_bound}, we obtain
\begin{align*}
    \|\mathsf K\|_{\mathrm{TV}}
    &\leq
    \|D_\#\mathsf C\|_{\mathrm{TV}}
    +
    \gamma\|\kappa\,(\mu\otimes\mu)\|_{\mathrm{TV}}
    +
    \gamma\|D_\#(\ell\mu)\|_{\mathrm{TV}}
    \\
    &\leq
    \mathsf C([0,M])
    +
    \gamma
    \int_{[0,M]^2}
    |\kappa(s,t)|
    \,d\mu(s)d\mu(t)
    +
    \gamma
    \int_{[0,M]}
    |\ell(s)|
    \,d\mu(s)
    \\
    &\leq
    \mathsf C([0,M])
    +
    2\gamma\mathcal E.
\end{align*}
This proves \eqref{eq:GR-candidate-TV}.

It remains to compute the bivariate Stieltjes transform
of $\mathsf K$.
Fix $\zeta,\xi\in\mathbb C^+$ with $\zeta\neq\xi$.
For the diagonal pushforward of $\mathsf C$,
the following identity
\[
    \frac1{(s-\zeta)(s-\xi)}
    =
    \frac1{\zeta-\xi}
    \left(
        \frac1{s-\zeta}
        -
        \frac1{s-\xi}
    \right)
\]
gives
\begin{align}
\label{eq:GR-representation-diagonal-C-transform}
    m_{D_\#\mathsf C}(\zeta,\xi)
    =
    \int_{[0,M]}
    \frac1{(s-\zeta)(s-\xi)}
    \,d\mathsf C(s)
    =
    \frac{
        m_{\mathsf C}(\zeta)
        -
        m_{\mathsf C}(\xi)
    }{
        \zeta-\xi
    }.
\end{align}
The Stieltjes transform of $\ell\mu$ is given by
\begin{align}
\label{eq:stielt_L}
    m_{\ell\mu}(z)
    =
    \int_{[0,M]}
    \frac{\ell(s)}{s-z}
    \,d\mu(s),
    \qquad
    z\in\mathbb C^+.
\end{align}
Next, note that the integrals used below are absolutely
convergent. Indeed,
\[
    \int_{[0,M]^2}
    \left|
        \frac{\kappa(s,t)}
             {(s-\zeta)(t-\xi)}
    \right|
    \,d\mu(s)d\mu(t)
    \leq
    \frac{\mathcal E}
         {\Im\{\zeta\}\Im\{\xi\}}
    <
    \infty,
\]
and, by \eqref{eq:ell_bound},
\[
    \int_{[0,M]}
    \left|
        \frac{\ell(s)}{s-z}
    \right|
    \,d\mu(s)
    \leq
    \frac{\mathcal E}{\Im\{z\}}
    <
    \infty.
\]
The Stieltjes transform of $\kappa\,(\mu\otimes\mu)$ is given by
\[
    m_{\kappa\,(\mu\otimes\mu)}(\zeta,\xi)
    =
    \int_{[0,M]^2}
    \frac{\kappa(s,t)}
         {(s-\zeta)(t-\xi)}
    \,d\mu(s)d\mu(t).
\]
Using
$\zeta-\xi=(\zeta-s)+(s-t)+(t-\xi)$,
we obtain
\begin{align}
\label{eq:zwischenschritt}
\begin{split}
    (\zeta-\xi)m_{\kappa\,(\mu\otimes\mu)}(\zeta,\xi)
    &=
    -
    \int_{[0,M]^2}
    \frac{\kappa(s,t)}
         {t-\xi}
    \,d\mu(s)d\mu(t)
    \\
    &\quad+
    \int_{[0,M]^2}
    \frac{(s-t)\kappa(s,t)}
         {(s-\zeta)(t-\xi)}
    \,d\mu(s)d\mu(t)
    \\
    &\quad+
    \int_{[0,M]^2}
    \frac{\kappa(s,t)}
         {s-\zeta}
    \,d\mu(s)d\mu(t).
\end{split}
\end{align}
Since $\kappa(s,t)=\kappa(t,s)$,
Fubini's theorem and the definition of $\ell$ give
\[
    \int_{[0,M]}
    \kappa(s,t)
    \,d\mu(s)
    =
    \ell(t)
\]
for $\mu$-almost every $t$.
Consequently, the first and third terms in
\eqref{eq:zwischenschritt} are
$-m_{\ell\mu}(\xi)$ and $m_{\ell\mu}(\zeta)$
defined in \eqref{eq:stielt_L}, respectively.
Moreover, by the definition of $\kappa$,
\[
    (s-t)\kappa(s,t)=\rho(s)-\rho(t),
    \qquad s,t\in[0,M],
\]
including the case $s=t$, where both sides vanish.
Therefore,
\begin{align*}
    (\zeta-\xi)m_{\kappa\,(\mu\otimes\mu)}(\zeta,\xi)
    =
    -m_{\ell\mu}(\xi)
    +
    m_{\ell\mu}(\zeta)
    +
    \int_{[0,M]^2}
    \frac{\rho(s)-\rho(t)}
         {(s-\zeta)(t-\xi)}
    \,d\mu(s)d\mu(t).
\end{align*}
Since $\rho\,d\mu=d\mathsf C$ and
$\mu=\mathsf C+\mathsf Q$, Fubini's theorem gives
\begin{align*}
    \int_{[0,M]^2}
    \frac{\rho(s)-\rho(t)}
         {(s-\zeta)(t-\xi)}
    \,d\mu(s)d\mu(t)
    &=
    \left(
        \int_{[0,M]}
        \frac{\rho(s)}{s-\zeta}
        \,d\mu(s)
    \right)
    \left(
        \int_{[0,M]}
        \frac1{t-\xi}
        \,d\mu(t)
    \right)
    \\
    &\quad-
    \left(
        \int_{[0,M]}
        \frac1{s-\zeta}
        \,d\mu(s)
    \right)
    \left(
        \int_{[0,M]}
        \frac{\rho(t)}{t-\xi}
        \,d\mu(t)
    \right)
    \\
    &=
    m_{\mathsf C}(\zeta)
    \left\{
        m_{\mathsf C}(\xi)
        +
        m_{\mathsf Q}(\xi)
    \right\}
    -
    \left\{
        m_{\mathsf C}(\zeta)
        +
        m_{\mathsf Q}(\zeta)
    \right\}
    m_{\mathsf C}(\xi)
    \\
    &=
    m_{\mathsf C}(\zeta)m_{\mathsf Q}(\xi)
    -
    m_{\mathsf C}(\xi)m_{\mathsf Q}(\zeta).
\end{align*}
It follows that
\begin{align}
\label{eq:GR-representation-eta-L-identity}
    m_{\kappa\,(\mu\otimes\mu)}(\zeta,\xi)
    -
    \frac{
        m_{\ell\mu}(\zeta)
        -
        m_{\ell\mu}(\xi)
    }{
        \zeta-\xi
    }
    =
    \frac{
        m_{\mathsf C}(\zeta)m_{\mathsf Q}(\xi)
        -
        m_{\mathsf C}(\xi)m_{\mathsf Q}(\zeta)
    }{
        \zeta-\xi
    }.
\end{align}
Similarly,
\begin{align*}
    m_{D_\#(\ell\mu)}(\zeta,\xi)
    =
    \int_{[0,M]}
    \frac{\ell(s)}{(s-\zeta)(s-\xi)}
    \,d\mu(s)
    =
    \frac{
        m_{\ell\mu}(\zeta)
        -
        m_{\ell\mu}(\xi)
    }{
        \zeta-\xi
    }.
\end{align*}
Combining this identity with
\eqref{eq:GR-representation-diagonal-C-transform},
\eqref{eq:GR-representation-eta-L-identity},
and \eqref{eq:GR-candidate-main}, using
$\rho\,d\mu=d\mathsf C$, we obtain
\begin{align}
\label{eq:GR-representation-K-transform-cq}
    m_{\mathsf K}(\zeta,\xi)
    =
    \frac{
        m_{\mathsf C}(\zeta)
        -
        m_{\mathsf C}(\xi)
    }{
        \zeta-\xi
    }
    +
    \gamma
    \frac{
        m_{\mathsf C}(\zeta)m_{\mathsf Q}(\xi)
        -
        m_{\mathsf C}(\xi)m_{\mathsf Q}(\zeta)
    }{
        \zeta-\xi
    }.
\end{align}
Recall \eqref{eq:GR-C-transform} and \eqref{eq:Q:stielt:def} and note that
\[
    1+\gamma m_{\mathsf Q}(z)
    =
    -\frac1{zv(z)}.
\]
Consequently,
\begin{align*}
    &
    m_{\mathsf C}(\zeta)
    -
    m_{\mathsf C}(\xi)
    +
    \gamma
    \left\{
        m_{\mathsf C}(\zeta)m_{\mathsf Q}(\xi)
        -
        m_{\mathsf C}(\xi)m_{\mathsf Q}(\zeta)
    \right\}
    \\
    &\quad=
    m_{\mathsf C}(\zeta)
    \left\{
        1+\gamma m_{\mathsf Q}(\xi)
    \right\}
    -
    m_{\mathsf C}(\xi)
    \left\{
        1+\gamma m_{\mathsf Q}(\zeta)
    \right\}
    \\
    &\quad=
    \frac1{\zeta\xi}
    \left\{
        \frac1{v(\xi)}
        \int_{[0,\infty)}
        \frac{t}{1+t v(\zeta)}
        \,dG(t)
        -
        \frac1{v(\zeta)}
        \int_{[0,\infty)}
        \frac{t}{1+t v(\xi)}
        \,dG(t)
    \right\}
    \\
    &\quad=
    \frac1{\zeta\xi v(\zeta)v(\xi)}
    \int_{[0,\infty)}
    t
    \left\{
        \frac{v(\zeta)}{1+t v(\zeta)}
        -
        \frac{v(\xi)}{1+t v(\xi)}
    \right\}
    \,dG(t)
    \\
    &\quad=
    \frac{
        v(\zeta)-v(\xi)
    }{
        \zeta\xi v(\zeta)v(\xi)
    }
    \int_{[0,\infty)}
    \frac{t}
         {\{1+t v(\zeta)\}\{1+t v(\xi)\}}
    \,dG(t).
\end{align*}
Substituting this into
\eqref{eq:GR-representation-K-transform-cq} yields
\begin{align}
\label{eq:rep:m_K_neq_points}
    m_{\mathsf K}(\zeta,\xi)
    =
    \frac1{\zeta\xi}
    \frac{
        v(\zeta)-v(\xi)
    }{
        (\zeta-\xi)v(\zeta)v(\xi)
    }
    \int_{[0,\infty)}
    \frac{t}
         {\{1+t v(\zeta)\}\{1+t v(\xi)\}}
    \,dG(t)
    =
    L(\zeta,\xi)
\end{align}
for $\zeta,\xi\in\mathbb C^+$ satisfying $\zeta\neq\xi$. Since $\mathsf K$ is finite, its bivariate Stieltjes
transform is continuous on $(\mathbb C^+)^2$.
Moreover, since $v$ is holomorphic on $\mathbb C^+$,
we have
\[
    \frac{v(\zeta)-v(\xi)}{\zeta-\xi}
    \longrightarrow
    \frac{d}{d\zeta}v(\zeta)
\]
as $\xi\to\zeta$.
Since $G$ is compactly supported by
Assumptions~\ref{ass:GR-basic-3} and \ref{ass:GR-basic-5},
and $\Im\{v(z)\}>0$ for $z\in\mathbb C^+$,
the factors $1+t v(\zeta)$ and $1+t v(\xi)$
in \eqref{eq:rep:m_K_neq_points} are bounded away
from zero uniformly for $t\in\operatorname{supp}(G)$
whenever $\xi$ is sufficiently close to $\zeta$.
Hence, by dominated convergence,
\[
    \int_{[0,\infty)}
    \frac{t}
    {\{1+t v(\zeta)\}\{1+t v(\xi)\}}
    \,dG(t)
    \longrightarrow
    \int_{[0,\infty)}
    \frac{t}{\{1+t v(\zeta)\}^2}
    \,dG(t)
\]
as $\xi\to\zeta$.
Together with the continuity of $m_{\mathsf K}$,
this proves \eqref{eq:GR-candidate-transform}
also when $\zeta=\xi$. 

Finally, define
\[
    \widetilde B(g,h)
    \coloneqq
    \int_{[0,M]^2}
    g(s)h(t)
    \,d\mathsf K(s,t),
    \qquad
    g,h\in C([0,M];\mathbb C).
\]
Since $\mathsf K$ is a finite signed Borel measure,
$\widetilde B$ is a bounded bilinear form satisfying
\[
    |\widetilde B(g,h)|
    \leq
    \|\mathsf K\|_{\mathrm{TV}}
    \|g\|_{\infty}
    \|h\|_{\infty}.
\]
Furthermore, for every $\zeta,\xi\in\mathbb C^+$,
\begin{align*}
    \widetilde B
    \bigl(
        r_\zeta,
        r_\xi
    \bigr)
    =
    m_{\mathsf K}(\zeta,\xi)
    =
    L(\zeta,\xi)
    =
    B
    \bigl(
        r_\zeta,
        r_\xi
    \bigr).
\end{align*}
The uniqueness assertion in
Theorem~\ref{lem:GR-bilinear-convergence}
therefore gives $\widetilde B=B$.
This proves
\eqref{eq:GR-bilinear-measure-representation}
and completes the proof.
\end{proof}

\subsection{Proof of Corollary~\ref{cor:ridge-optimality2}}

\begin{proof}[Proof of Corollary~\ref{cor:ridge-optimality2}]
Write $r_f(s)=1-sf(s)$ and define
\[
A_n(f)
\coloneqq
\frac1n Xf(S_n)\Sigma_n f(S_n)X^\top,
\qquad
b_n(f)
\coloneqq
\frac1n Xf(S_n)\Sigma_n r_f(S_n)\beta_n^\star.
\]
Note that, by the definitions in \eqref{eq:def_R_n} and \eqref{eq:def-Ln},
\begin{align*}
    \mathrm R_{\Sigma_n}(\hat\beta(f)) = \mathbb E\left[
\mathcal L_n(f)\mid X
\right] .
\end{align*}
Then, using \eqref{eq:formula_diff_hat_beta_beta_star}, we obtain
\[
\mathcal L_n(f)
-\mathrm R_{\Sigma_n}(\hat\beta(f))
=
\frac1n\left\{
\varepsilon^\top A_n(f)\varepsilon
-\sigma^2\tr A_n(f)
\right\}
-2\varepsilon^\top b_n(f).
\]
Using the constant $M$ from
Theorem~\ref{lem:GR-bilinear-convergence},
we set $E_n=\{\|S_n\|_2\leq M\}$.
Hence, for every fixed $f\in C([0,\infty))$,
we have on $E_n$,
\[
\|A_n(f)\|_2\leq C_f,
\qquad
\|b_n(f)\|_2^2\leq \frac{C_f}{n}.
\]
Set $q=2+\delta/2$. Note that $E_n$ is measurable with respect to the $\sigma$-field generated by $X$. 
Since $\varepsilon$ is independent of $X$,
\cite[Lemma~B.26]{BaiSilverstein2010},
applied conditionally on $X$, gives
\[
\mathbb E\left[
\left|
\frac1n\left\{
\varepsilon^\top A_n(f)\varepsilon
-\sigma^2\tr A_n(f)
\right\}
\right|^q
\mathbf1_{E_n}
\right]
\leq C_f n^{-q/2}.
\]
Moreover, independence and the finite fourth
moment of the errors yield
\[
\mathbb E\left[
|\varepsilon^\top b_n(f)|^4\mid X
\right]
=
3\sigma^4\|b_n(f)\|_2^4+
\bigl(\mathbb E[\varepsilon_1^4]-3\sigma^4\bigr)
\sum_{i=1}^n b_{n,i}(f)^4\leq C\|b_n(f)\|_2^4,
\]
where $b_{n,i}(f)$ are the components of
$b_n(f)$. Consequently,
\[
\mathbb E\left[
|\varepsilon^\top b_n(f)|^4
\mathbf1_{E_n}
\right]
\leq C_f n^{-2}.
\]
Since $q>2$, Markov's inequality and
Borel--Cantelli, together with the definition of
$E_n$ and Theorem~\ref{lem:GR-bilinear-convergence},
imply \eqref{eq:risk_cond_XY}
for every fixed $f\in C([0,\infty))$. The asserted
characterization follows directly from
Theorem~\ref{cor:ridge-optimality}.
\end{proof}

\section{Further technical results for Section~\ref{sec:proof-bilinear}} \label{appendix_D}

\subsection{Proposition~\ref{prop:bilinear-convergence} and its proof}

\begin{proposition}\label{prop:bilinear-convergence}
Let $X$ be a normed space over
$\mathbb F\in\{\mathbb R,\mathbb C\}$, and let
$B_n:X^2\to\mathbb F$, $n\geq1$, be continuous bilinear forms
satisfying
\[
    \sup_{n\geq1}\|B_n\|<\infty.
\]
Let $D\subseteq X$ be dense and suppose that
$B_n(d_1,d_2)$ converges for every $d_1,d_2\in D$.
Then there exists a unique continuous bilinear form
$B:X^2\to\mathbb F$ satisfying
\[
    B_n(x,y)\longrightarrow B(x,y),
\]
for all $x,y\in X$, and
\[
    \|B\|
    \leq
    \sup_{n\geq1}\|B_n\|.
\]
\end{proposition}

\begin{proof}
Set $C\coloneqq\sup_{n\geq1}\|B_n\|<\infty$.
Fix $x,y\in X$. Since $D$ is dense in $X$, there exist
sequences $x_k,y_k\in D$ such that $x_k\to x$ and $y_k\to y$.
By bilinearity and the uniform bound on $\|B_n\|$,
\begin{align*}
    |B_n(x,y)-B_n(x_k,y_k)|
    &\leq
    |B_n(x-x_k,y)|+|B_n(x_k,y-y_k)|
    \\
    &\leq
    C\bigl(
        \|x-x_k\|\|y\|
        +
        \|x_k\|\|y-y_k\|
    \bigr).
\end{align*}
Consequently, for every $n,m,k\geq1$,
\begin{align*}
    |B_n(x,y)-B_m(x,y)|
    &\leq
    |B_n(d_k,e_k)-B_m(d_k,e_k)|
    \\
    &\quad+
    2C\bigl(
        \|x-d_k\|\|y\|
        +
        \|d_k\|\|y-e_k\|
    \bigr).
\end{align*}
For every fixed $k$, the first term tends to zero as
$n,m\to\infty$ by the assumed convergence on $D\times D$.
The remaining terms tend to zero as $k\to\infty$.
Thus $\{B_n(x,y)\}_{n\geq1}$ is a Cauchy sequence in
$\mathbb F$, and we may define
\[
    B(x,y)\coloneqq\lim_{n\to\infty}B_n(x,y).
\]
Passing to the limit in the bilinearity identities for $B_n$
shows that $B$ is bilinear. Moreover,
\[
    |B(x,y)|
    =
    \lim_{n\to\infty}|B_n(x,y)|
    \leq
    C\|x\|\|y\|,
\]
so $B$ is continuous and $\|B\|\leq C$.
Uniqueness follows from the uniqueness of the pointwise limits.
\end{proof}

\subsection{Proof of Lemma~\ref{lem:abstract-bilinear-extension}}

\begin{proof}[Proof of Lemma~\ref{lem:abstract-bilinear-extension}]
Define
\[
    \mathcal R_M
    \coloneqq
    \operatorname{span}
    \left\{
        r_z:z\in\mathbb C^+
    \right\},
\]
and let $\overline{\mathcal R_M}$ denote its closure in the supremum
norm. We first show that
$\overline{\mathcal R_M}=C([0,M];\mathbb C)$.
For $t>0$ and $x\in[0,M]$,
\[
    \left|-it\,r_{it}(x)-1\right|
    =
    \frac{x}{|x-it|}
    \leq
    \frac{M}{t}.
\]
Thus $-it\,r_{it}$ converges uniformly to the constant
function $1$ as $t\to\infty$, and hence,
$1\in\overline{\mathcal R_M}$. Moreover,
\[
    x r_z(x)
    =
    1+z r_z(x),
    \qquad z\in\mathbb C^+,\quad x\in[0,M].
\]
Since $\overline{\mathcal R_M}$ is a linear space, it follows that $xf\in\overline{\mathcal R_M}$, for $ f\in\mathcal R_M$.
Now let $f\in\overline{\mathcal R_M}$ and choose
$f_k\in\mathcal R_M$ such that
\[
    \|f_k-f\|_{\infty}\longrightarrow0.
\]
By the preceding observation,
$xf_k\in\overline{\mathcal R_M}$ for every $k$. Moreover,
\[
    \|xf_k-xf\|_{\infty}
    \leq
    M\|f_k-f\|_{\infty}
    \longrightarrow0.
\]
Since $\overline{\mathcal R_M}$ is closed, it follows that
$xf\in\overline{\mathcal R_M}$. Thus,
$\overline{\mathcal R_M}$ is invariant under multiplication by the
identity function $x\mapsto x$.
Since $1\in\overline{\mathcal R_M}$, induction shows that
$\overline{\mathcal R_M}$ contains every complex polynomial.
The Stone--Weierstrass theorem therefore gives
\begin{equation*}
\label{eq:resolvent-span-dense}
    \overline{\mathcal R_M}
    =
    C([0,M];\mathbb C).
\end{equation*}
For
$g=\sum_{i=1}^k a_i r_{\zeta_i}$ and
$h=\sum_{j=1}^{\ell} b_j r_{\xi_j}$
in $\mathcal R_M$, bilinearity gives
\begin{align*}
    B_n(g,h)
    =
    \sum_{i=1}^k\sum_{j=1}^{\ell}
    a_i b_j
    B_n\left(r_{\zeta_i},r_{\xi_j}\right)
    \longrightarrow
    \sum_{i=1}^k\sum_{j=1}^{\ell}
    a_i b_j L(\zeta_i,\xi_j),
\end{align*}
where $L(\zeta,\xi)$ denotes the limit of $B_n(r_{\zeta}, r_{\xi})$.
Thus $B_n(g,h)$ converges for every $g,h\in\mathcal R_M$. Proposition~\ref{prop:bilinear-convergence}, applied with
$X=C([0,M];\mathbb C)$ and $D=\mathcal R_M$, yields a unique continuous
bilinear form $B$ with $\|B\|\leq C_B$ such that
$B_n(g,h)\to B(g,h)$ for every $g,h\in C([0,M];\mathbb C)$.
In particular, since $r_\zeta, r_\xi \in\mathcal R_M$,
\[
B\left(r_\zeta,r_\xi \right)
    =
    \lim_{n\to\infty}
    B_n\left(r_\zeta,r_\xi\right)
    =
    L(\zeta,\xi),
    \qquad
    \zeta,\xi\in\mathbb C^+.
\]
This proves the assertion.
\end{proof}

\stopcontents[supplement]

\end{document}